\documentclass[10pt]{article}
\usepackage{multirow}
\usepackage{comment}
\usepackage{graphicx}
\usepackage{amsmath,bm}
\usepackage{amssymb}
\usepackage{amsfonts,bbm}
\usepackage{amsthm}
\usepackage{xcolor}
\usepackage{enumitem}

\numberwithin{equation}{section}
\def\R{\mathbb{R}}
\def\Z{\mathbb{Z}}
\def\To{\mathbb{T}}
\def\P{\mathbb{P}}

\def\kk{{\mathbf k}}
\def\cc{{\mathbf c}}

\def\DD{\mathbb{ D}}

\def\La{\Lambda}
\def\T{\mathcal{T}}
\def\P{\mathcal{P}}
\def\f{\textbf{f}}
\def\uu{\textbf{u}}
\def\vv{\textbf{v}}
\def\w{\textbf{w}}
\def\varphib{\bm{\varphi}}
\def\omegab{\bm{\omega}}

\def\neweq#1{\begin{equation}\label{#1}}
\def\endeq{\end{equation}}

\def\eps{\varepsilon}
\newtheorem{theorem}{Theorem}[section]
\newtheorem{corollary}[theorem]{Corollary}
\newtheorem{lemma}[theorem]{Lemma}
\newtheorem{proposition}[theorem]{Proposition}
\newtheorem{computation}[theorem]{Computation}

\newtheorem{remark}[theorem]{Remark}

\usepackage[linktocpage,colorlinks=true,raiselinks]{hyperref}
\hypersetup{urlcolor=blue, citecolor=red, linkcolor=blue}

\title{ \bf On a 3D Stokes eigenvalue problem  under Navier
  slip-with-friction boundary conditions and applications to Navier-Stokes equations}

\author{Luigi C. Berselli$^\dag$ - Alessio Falocchi$^\ddag$ -
  Rossano Sannipoli$^\dag$\\ {\small $\dag$ Dipartimento di Matematica, Dipartimento di Eccellenza 2023-2027 
    - Università di Pisa, Italy}\\
    {\small $\ddag$ Dipartimento di
    Matematica, Dipartimento
di Eccellenza 2023-2027 - Politecnico di Milano, Italy}
}
\date{}

\begin{document}
\maketitle
\begin{abstract}
    In this paper we consider, by means of a precise spectral analysis,
    the 3D Navier-Stokes equations endowed with Navier
    slip-with-friction boundary conditions. We study the problem in a
    very simple geometric situation as the region between two parallel
    planes, with periodicity along the two planes. This setting, which
    is often used in the theory of boundary layers,  requires some
    special treatment for what concerns the functional setting and
    allows us to characterize in a rather explicit manner eigenvalues
    and eigenfunctions of the associated Stokes problem. These, will
    be then used in order to identify infinite dimensional  classes of data
    leading to  global strong solutions for the corresponding
    evolution Navier-Stokes equations. 

\end{abstract}
\setcounter{tocdepth}{1}
\tableofcontents
\section{Introduction}
\label{sec:introduction}
In this paper we study the incompressible 3D Navier-Stokes equations with Navier slip-with-friction boundary conditions;  our aim is to characterize the eigen-values/functions in order to prove later on results of global existence for large classes of nontrivial strong solutions.  More precisely, the basic building block is the following eigenvalue problem 
\begin{equation}
\label{stokes}
\left\{
\begin{array}{ll}
-\Delta \uu+\nabla q=\lambda \uu\quad & \mbox{in }\Omega\,
\\
\nabla\cdot \uu=0\quad & \mbox{in }\Omega\,
\\
\uu\cdot\nu=0\quad & \mbox{on }\Gamma\,
\\
2(\DD\uu\cdot\nu)\cdot\tau+\beta \uu\cdot\tau=0\quad & \mbox{on }\Gamma\,,
\\
\end{array}\right.
\end{equation}
where  $\uu=(u,v,w)$  is the fluid velocity vector, $q$ is the fluid pressure, and $\lambda\in \R$ is an eigenvalue. 
Moreover, the unit vectors $\nu$ and $\tau$ are respectively the outward normal and the tangential vectors associated to $\Gamma$, $\DD\uu$ is the strain tensor 
\begin{equation*}
\DD\uu:=\frac{\ \nabla \uu+\nabla^\top \uu\ }{2},
\end{equation*}
$\beta \ge 0$ is the friction coefficient (a given parameter), and
--without loss of generality-- we consider the kinematic viscosity
equal to 1. The non-negative number $\beta$ determines the friction
and depends on the fluid viscosity and the roughness of the boundary; it is the inverse of the slip length, see G{\'e}rard-Varet  and
Masmoudi~\cite{GVM2010}.  These boundary conditions are also at the
basis of several Large Eddy Simulation (LES) models for turbulent flows, see  Chac{\'o}n Rebollo and  Lewandowski~\cite{CL2014} and~\cite{BIL2006}, where the problem becomes  nonlinear, since $\beta$ could be also linked through a nonlinear law to the local tangential velocity.  

Problem~\eqref{stokes} is the starting point to analyze evolution 3D Navier-Stokes equations  with Navier slip-with-friction boundary conditions: for $T>0$, we set $Q_T:= \Omega\times (0,T)$ and we consider the problem
\begin{equation}
\label{ns}
\left\{
\begin{array}{ll}
\vv_t-\Delta \vv+(\vv\cdot\nabla)\vv+\nabla q=\f\quad &\mbox{in }Q_T
\\
\nabla\cdot \vv=0\quad &\mbox{in }Q_T
\\
\vv\cdot\nu=0\quad &\mbox{on } \Gamma\times (0,T)
\\
2(\DD\vv\cdot\nu)\cdot\tau+\beta \vv\cdot\tau=0\quad &\mbox{on } \Gamma\times (0,T)
\\
\vv(0,x,y,z)=\vv_0(x,y,z)\quad & \mbox{in }\Omega\,,
\end{array}\right.
\end{equation}
where in the sequel we put $\f =\mathbf{0}$ (but this assumption is unessential).
%
The Navier condition is often written as $2(\DD\vv\cdot\nu)\cdot\tau_j+\beta \vv\cdot\tau_j=0$, where $\tau_j$, for $j=1,2$, are any two linearly independent tangential unit vectors.

By restricting to a simple (but physically meaningful) domain, the  spectral
analysis of problem~\eqref{stokes} allows us to study in an explicit way the existence and regularity of solutions to~\eqref{ns}, with special emphasis on the role of $\beta$. We point out that
for $\beta=0$ we fall in the pure Navier (frictionless) case, while in the limit $\beta\to\infty$, we
recover the Dirichlet case.

The analytical study of problem~\eqref{ns} and of the corresponding
linearized steady version (which is a basic building block of the
theory) dates back to Solonnikov and \v{S}\v{c}adilov~\cite{SS1973},
Beir\~ao da Veiga~\cite{BDV2004,BDV2005} for the Hilbertian theory;
see also the recent results by Amrouche \textit{et
  al.}~\cite{ATACG2021,ABA2019,AR2014}. The most relevant results in
the steady case are the existence and uniqueness ($\beta>0)$ of weak
and strong solutions. This can be used to formulate a variational
theory also for the time-evolution problem; more precisely,
restricting to small data/small times, a fixed point argument can be
used to tackle the nonlinear evolution problem. Particular care is
given in the above papers to determine if there exist nontrivial
solutions of the homogeneous problem, in particular geometries (with
rotational invariance) and in dependence on the values of
$\beta\geq 0$. Note that we focus on the 3D case, while the two
dimensional one has been studied by several authors, also to consider
the vanishing viscosity limit, employing tools typical of the 2D
problem (cf.~Clopeau, Mikeli{\'c}, and Robert~\cite{CMR1998} and Lopes
Filho, Nussenzveig Lopes, and Planas~\cite{LFNLP2005}); their
technique are not replicable in the 3D case.

Here, we consider a simple geometric setting, which corresponds to the space between two parallel planes, with periodicity along the axes on the two planes; the
Navier boundary conditions involve only the vertical direction. This can be
seen as a sort of ``doubly-periodic channel" and is a test case used
in several computations of turbulent flows or in the study of boundary
layers, cf.~Pinier~\textit{et al.}~\cite{PLMC2021}, and coupling
fluid-fluid/fluid-air, cf.~\cite{BLL2024}; this is the simplest setting in
which  appears the effect of the boundary. The analytical
results we prove, can be used for a priori testing  experimental data
related with channel flows.

More precisely, we set 
$$
\Omega:=\To^2\times(-1,1),
$$
where $\To:=\R/(2\pi\Z)$ and we denote by $\Gamma:=\To^2\times\{\pm1\}$ a portion of 
the boundary of $\Omega$.  The fact that the domain is a periodic torus in the $x$ and $y$
variables means that functions are $2\pi$-periodic in these two directions, i.e., %
\begin{equation}
\label{per3dNSE}
\vv(t,x,y,z)=\vv(t,x+2k\pi,y,z)=\vv(t,x,y+2k\pi,z)\qquad \forall\,k\in\mathbb{Z},\ \text{ and a.e. }t\in(0,T),
\end{equation}
completing the definition of the boundary conditions, which are of Navier slip-with-friction only at $z=\pm1$. Obviously in the sequel we will  use as tangent vectors $\tau_j$, $j=1,2$,  the unit vectors in the direction of the axis $x$ and $y$, respectively.

Preliminary, we recall some classical theorems of existence and regularity of solutions to the steady and to the evolution case, since they are not generally stated in this geometric setting. The local existence result for strong solutions in the time evolution case is based on abstract knowledge of the eigenfunctions.

Then we pass to study the spectrum of \eqref{stokes}; in this geometric setting, we adapt the precise spectral analysis performed
by Rummler~\cite{Rum1997a,Rum1997b,Rum2000} and in~\cite{FG2023} in order to explicitly evaluate the eigenvalues and the eigenfunctions; see also 
the results by Hu~\cite{Hu2007} in the case of a very thin domain. 
Several topics, as for instance: the determination of the smallest eigenvalue; the asymptotic limit in terms of $\beta$; the characterization of
eigenfunctions associated with trivial and nontrivial pressure are significantly different from the previous works. For instance, although the domain is flat, we cannot use the reflection technique used
in the case of Navier boundary conditions without friction~\cite{Ber2010a}, making the problem much closer to the Dirichlet one. In particular, if $\beta>0$, passing to
the curl system drives towards a problem which cannot be handled, as in the usual Dirichlet case, due to the presence of boundary terms which cannot be controlled. The spectral analysis requires very precise and long computations in order to extract the main properties of the spectrum; due to the presence of $\beta$ the eigenvalues are obtained solving transcendental equations.  The spectral information presented in~\cite{FG2023} on a cubic domain having frictionless Navier boundary conditions ($\beta=0$) can be viewed as a particular case of our analysis.

After the full characterization of the eigen-values/functions for the Stokes problem, we pass to the main aim of the work: to build nontrivial families of global smooth solutions, starting from the spectral decomposition.
They are not variations of 2D problem or families corresponding to ``rarified'' sequences of eigenfunctions~\cite{FG2023}; instead, we identify infinite dimensional families of ``2+1/2'' solutions, which are globally smooth, without being two dimensional. These results are in the same spirit of the  3D extension of classical bi-periodic vortex solution by Taylor, found by Ethier and Steinman~\cite{RES94} and Antuono~\cite{Ant2020} (which are all analytical and ``Tri-periodic''). Such solutions
can be used as benchmarks for
checking the accuracy/debugging of 3D numerical codes for Direct Numerical Simulations (DNS) as well as for Large Eddy Simulations (LES). Moreover, the physics of the problem is related to study the onset of turbulence in the case of Navier slip-with-friction condition on a solid planar surface, hence in a situation typical of the turbulent channel flow.

Since many results are technically demanding (even to be stated), we first specify the notation and the functional setting, reporting 
all the results in a dedicated section. The various proofs, which are based on a few
cumbersome explicit computations, will be given in the subsequent
sections. 
\bigskip

\textbf{Plan of the paper:} In Section~\ref{sec:preliminary} we specify the functional setting proper to deal with our problem endowed by periodic and Navier boundary conditions. Moreover, we report some preliminary results which are known in classical Navier-Stokes framework and we adapt them
for this specific situation.
In Section~\ref{sec:main-results} we collect the main results concerning the full characterization of the spectrum of the Stokes operator and the identification of infinite dimensional families of initial data leading to smooth global solutions for the corresponding time dependent Navier-Stokes equations.
In Section~\ref{sec:4}  we prove the preliminary results, while  the proofs of the main results are collected in Section~\ref{sec:proofs}.
In Appendix~\ref{eigenfunctions1} we give some further details on the Stokes spectrum, e.g. the full (and lengthy!)  explicit expression of the eigenfunctions, while in Appendix~\ref{nonlin_proj} we report the Helmholtz-Weyl decomposition, needed to construct  the solutions of the related nonlinear time dependent problem.
\section{Preliminaries}
\label{sec:preliminary}
\subsection{Functional setting and notations}
We denote by $L^p(\Omega)$, for $1\leq p \leq\infty$, the usual Lebesgue spaces, equipped
with norm $\|\cdot\|_{L^p(\Omega)}$, while we denote by $H^s(\Omega)$,
for $s\geq0$, the classical Sobolev spaces. We often use
$\|\cdot\|_{L^p(\Gamma)}$ referring to the $L^p$ norm on the boundary
$\Gamma$. We shall use the same symbol for both scalar and vector
function spaces. In general we omit $dx$ and $dS$ for the integrals
over $\Omega$ and $\Gamma$, respectively.

Let us recall the definition of the usual spaces used in the treatment of the Navier-Stokes equations. The set $H$ is obtained by taking the closure of $\mathcal{V}:=\{\vv\in C^\infty_0(\Omega):\, \nabla\cdot \vv=0\}$ in $L^2(\Omega)$; being $\Omega$ Lipschitz and bounded, by~\cite[Theorem III.2.3]{galdi} we have
\begin{equation*}
H:=\{\vv\in L^2(\Omega);\, \nabla\cdot \vv=0,\, \vv\cdot\nu=0\mbox{ on }\Gamma\}\, ,
\end{equation*}
in which, the divergence is taken in distributional sense, while $\vv\cdot\nu$ denotes the normal trace of $\vv$.  We define $V_{\tau}$ as
\begin{equation*}
V_{\tau}:=H^1(\Omega)\cap H.
\end{equation*}
 By construction, $H$ is a closed subspace of $L^2(\Omega)$; therefore, $V_{\tau}$ is a closed subspace
of $H^1(\Omega)$. Moreover, we  denote the tangential component on $\Gamma$ of a vector-valued function $\vv$ as follows 
\begin{equation*}
    \vv_\tau := \vv-(\vv\cdot\nu)\,\nu.
\end{equation*}

Since we are in a special setting in which the functions are periodic in $x$ and $y$ directions and the boundary of the domain is given by $z=\pm1$, we use the Fourier series expansion in the $x$ and $y$ variables. More precisely, we decompose a field $\vv=(u,v,w)$ as follows 
\begin{equation*}
\vv(x,y,z)=\sum_{\kk=(k,h)\in\mathbb{Z}^2}\cc_{\kk}(z)e^{i(kx+hy)}\qquad\text{with}\quad  \overline{\cc}_{\kk}(z)=-\cc_{\kk}(z),
\end{equation*}
where $\cc_{\kk}(z):=\frac{1}{(2\pi)^2}\int_{\To}\vv(x,y,z) e^{-i(kx+hy)}$ are the vector-valued Fourier coefficients and $\overline{\cc}_{\kk}$ is the complex conjugate of $\cc_{\kk}$, see also~\cite{temam}. Hence, we introduce the following (partially) periodic functional spaces
\begin{equation*}
\begin{aligned}
    L^2_{\mathcal{P}}(\Omega)&=\Big\{\vv\in L^2(\Omega):\  \|\vv\|^2_{L^2_\P}:=\sum_{\kk\in\mathbb{Z}^2}\|\cc_{\kk}\|_{L^2(-1,1)}^2<\infty\Big\}\quad
\\
    H^1_{\mathcal{P}}(\Omega)&=\Big\{\vv\in H^1(\Omega):\ \|\vv\|^2_{H^1_\P}:= \sum_{\kk\in\mathbb{Z}^2}\big[(1+|\kk|^{2})\|\cc_{\kk}\|_{L^2(-1,1)}^2+\|\cc'_{\kk}\|_{L^2(-1,1)}^2\big]<\infty\Big\},
\end{aligned}
\end{equation*}
and the spaces
\begin{equation}
  \label{HP}
H_\P:=H\cap L^2_\P(\Omega)\qquad\text{and}\qquad V_\P:=V_\tau\cap H^1_\P(\Omega),
\end{equation}
on which we consider respectively the bilinear forms
\begin{equation*}
(\vv,\varphib):=\int_\Omega \vv\cdot \varphib\, \qquad \text{and}\qquad
\quad (\DD \vv,\DD \varphib):=\int_\Omega \DD \vv:\DD \varphib\,.
\end{equation*}
 Here ``$\,.\,$'' denotes the scalar product between vectors, while  ``$\,:\,$'' indicates the scalar product between matrices.  We denote by $\langle\cdot,\cdot\rangle$ the duality pairing between $V_{\P}$ and its dual space $V_{\P}^*$, and by $\|\cdot\|_{V_\P^*}$ the corresponding dual norm.
 
It is worth mentioning that, while the bilinear form $(\cdot,\cdot)$ defines a scalar product over $H_\P$, with corresponding norm $\|\vv\|_{L^2(\Omega)
}^2:=\int_\Omega |\vv|^2\,$ equivalent to $\|\vv\|_{H_\P}^2$,  the
situation could be different for $(\DD \vv,\DD \varphib)$. More
precisely, if $\beta>0$  $\|\DD \vv\|^2_{L^2(\Omega)}:=(\DD \vv,\DD
\vv)$  is an equivalent  norm on $V_\P$, while if $\beta=0$ then  $\|\DD \vv\|_{L^2(\Omega)}$ is a semi-norm, being the possibility of zero eigenvalue in the Stokes operator. In order to study the limit case $\beta\rightarrow\infty$ we introduce the subspace 
\begin{equation*}
    V^D_\mathcal{P}= V_\mathcal{P}\cap \{\vv \in H^1(\Omega) : \vv \cdot \tau = 0\},
\end{equation*}
corresponding to Dirichlet boundary conditions on $\Gamma$. We point out that the functional space to study the limit case $\beta=0$ is always $V_\mathcal{P}$. In the paper we put the indexes $D$ and $N$ to underline respectively when we are considering the Dirichlet boundary conditions ($\beta\rightarrow\infty$) and the Navier frictionless boundary conditions ($\beta=0$).

In the following proposition we emphasize an important advantage of dealing with the domain $\Omega$, in terms of equivalence between norms. The proof is given in Section~\ref{lem:equiv}.
\begin{proposition}
\label{lem:gradsym}
Let $\beta \ge 0$ then $\|\DD \cdot\|_{L^2(\Omega)}$ and $\|\nabla\cdot\|_{L^2(\Omega)}$ are equivalent semi-norms (norms if $\beta>0$) in $V_\P$, i.e. it holds
\begin{equation*}
      2  \|\DD \textbf{v}\|_{L^2(\Omega)}^2= \|\nabla\textbf{v}\|_{L^2(\Omega)}^2, \qquad \forall\, \vv \in V_{\mathcal P}.
\end{equation*}
\end{proposition}

Throughout the paper when we write $\widetilde\Delta\vv$ we mean  the projection of the $\Delta\vv$ onto the space $H_\P$, i.e. the Stokes operator; moreover, when we use the symbol $\mathbb{N}$ we refer to the natural numbers including $0$, while $\mathbb{N}_+$ refers to the strictly positive integers.

In the next section we extend to our framework some classical regularity results on the Stokes and Navier-Stokes equations under Navier boundary conditions. These preliminary results are useful in order to state the main results of the paper.
\subsection{Preliminary results}
First of all we consider the Stokes problem in $\Omega$ under Navier boundary conditions with a source term  $\f\in V_\P^*$    
\begin{equation}
\label{stokesf}
\left\{
\begin{array}{ll}
-\Delta \uu+\nabla q=\textbf{f}\quad & \mbox{in }\Omega\,\\
\nabla\cdot \uu=0\quad & \mbox{in }\Omega\,\\
\uu\cdot\nu=2(\DD\uu\cdot\nu)\cdot\tau+\beta \uu\cdot\tau=0\quad & \mbox{on }\Gamma\,,
\end{array}\right.
\end{equation}
whose weak formulation is 
\begin{equation*}
    2(\DD \uu,\DD \varphib)+\beta\int_{\Gamma}\uu_\tau\cdot\varphib_\tau=\langle\f, \varphi\rangle\qquad \forall\,\varphib\in V_\P.
\end{equation*}

 In the next proposition we state existence and uniqueness results for the Stokes problem.
We omit the proof since it easily follows adapting the results in the half-space~\cite{BDV2005} to our setting; see also~\cite{BDV2004} for the analysis of the problem in the case with non zero divergence and rotational symmetries.
\begin{proposition}\label{prop:Stokes-existence}
    Let $\beta\geq0$, $\f\in V_{\P}^*$ then there exists a weak solution $\uu\in V_{\P}$  of~\eqref{stokesf} such that

\begin{equation*}
    \|\DD\uu\|_{L^2(\Omega)}\leq C(\Omega)\|\f\|_{V_\P^*},
\end{equation*}
with $C(\Omega)>0$. The solution is unique if $\beta>0$, while if $\beta=0$ is
  unique up to constant vectors of the type $(a,b,0)$, $a,b\in\R$.
\end{proposition}

In the following statement, proved in Section~\ref{appendix:regularityevolutive}, we extend classical regularity results to the weak solution of~\eqref{stokesf}. They are fundamental in order to prove the compactness of the inverse of the Stokes operator and to ensure the existence of an Hilbert basis of eigenfunctions for such operator. We denote by $H^1_{\P,\#}(\Omega)$ the subspace of $H^1_{\P}(\Omega)$ with zero mean value over $\Omega$ and we identify the pressure, as in usual Stokes problems, by its zero mean value.
\begin{proposition}\label{lemma:regularity}
 Let $\beta\geq0$ and $\textbf{f}\in L^2(\Omega)$. Then, there exists a strong solution $(\textbf{u},q)\in \big(H^2(\Omega)\cap V_\P\big)\times H^1_{\P,\#}(\Omega)$ of~\eqref{stokesf} such that
\begin{equation}\label{estimategrauf}
    \| \textbf{u}\|_{H^2(\Omega)}+\|q\|_{H^1(\Omega)}\le C \|\textbf{f}\|_{L^2(\Omega)},
\end{equation}
with $C>0$. If $\beta>0$ the velocity $\uu$ is unique, while if $\beta=0$ is unique up to constant vectors of the type $(a,b,0)$, $a,b\in\R$.
 Moreover, $\|\widetilde\Delta\uu\|_{L^2(\Omega)}$ is a norm on $H^2(\Omega)\cap H_{\P}$, equivalent to the norm induced by $H^2(\Omega)$.
\end{proposition}
We give now some preliminary results on the Stokes eigenvalue problem~\eqref{stokes}, which
in weak form reads 
\begin{equation}\label{weak}
    2(\DD \uu,\DD \varphib)+\beta\int_{\Gamma}\uu_\tau\cdot\varphib_\tau=\lambda(\uu, \varphib) \qquad \forall\,\,\varphib\in V_\P;
\end{equation}
taking $\varphib=\uu$, we notice that the first eigenvalue
of~\eqref{stokes} is non-negative, i.e. $\lambda_0\geq 0$ for all
$\beta\geq0$. In the following proposition, proved in Section~\ref{proof:lambda0}, we give a more precise statement on the first eigenvalue of~\eqref{stokes} with respect to the choice of the friction coefficient $\beta$. When $\beta=0$ the presence of the zero eigenvalue in the spectrum of~\eqref{stokes} has important consequences in terms of Poincaré inequalities, see also~\cite{falocchi2022remarks} for a similar problem.
\begin{proposition}\label{lemma0}
   Let $\lambda_0$ be the first eigenvalue of~\eqref{stokes}. If $\beta>0$ then $\lambda_0>0$ and if $\beta=0$ then $\lambda_0=0$, corresponding to the two linear independent constant eigenfunctions $(1,0,0)$ and $(0,1,0)$. Moreover, there exist  $C_0=C_0(\Omega,\beta), C_0^D(\Omega), C_0^N(\Omega)>0$ such that
   \begin{equation}
   	\label{poincare}
   		\|\vv\|^2_{L^2(\Omega)}\leq\begin{cases}
   			C_0\big(2\|\DD \vv\|^2_{L^2(\Omega)}+\beta\|\vv_\tau\|^2_{L^2(\Gamma)}\big)\qquad\, &\text{if }\beta>0,\quad\hspace{2mm}\forall\,\vv\in V_\P
   			\\  
   			2 C^D_0\
   			\|\DD \vv\|^2_{L^2(\Omega)}\qquad \hspace{25mm}&\text{if }\beta\rightarrow\infty,\quad\forall\,\vv\in V^D_\P
   			\\
   			2C_0^N
   			\|\DD \vv\|^2_{L^2(\Omega)}\qquad \hspace{5mm}&\text{if }\beta=0,\quad\hspace{2.5mm}\forall\,\vv\in V_\P\cap Z^{\perp}_\P,
   		\end{cases} 
   \end{equation}
where the space $Z^{\perp}_\P$ is the orthogonal complement of the kernel of the Stokes operator.
\end{proposition}
 Since in the next sections we study the full spectrum of~\eqref{stokes} we will provide more precise information on the constants in~\eqref{poincare}, see Remark~\ref{poincare_constants}.

We pass now to state some preliminary results on  the  evolution Navier-Stokes problem~\eqref{ns}, that in weak form reads 
$$
\begin{aligned}
    &\int_0^T(\vv(t),\varphib)\, \phi'(t)\,dt+(\vv_0,\varphib)\phi(0)=\int_0^T\left\{2(\DD \vv(t),\DD \varphib) +\beta\int_{\Gamma}\vv_\tau(t)\cdot\varphib_\tau+( (\vv(t)\cdot \nabla )\,\vv(t)\cdot\varphib) \right\}\phi(t)\,dt,
\end{aligned}
$$
for all $\varphib\in V_\P$ and for all $\phi\in C^\infty_0[0,T)$.

The global existence of weak solutions as well as local results on strong solutions can be found in Beir\~ao da Veiga~\cite{BdV2007}, even if  in the slightly different case of a smooth bounded domain.  The following result is a direct consequence of the variational formulation and is similar to the standard theorem for Leray-Hopf weak solutions under Dirichlet boundary conditions.  
\begin{lemma}
Let $T>0$ and $\vv_0\in H_\P$. Then, there exists at least a weak solution $\vv\in L^\infty(0,T;H_\P)\cap L^2(0,T;V_\P)$ to~\eqref{ns} in the time interval $[0,T]$, satisfying the energy  inequality
\begin{equation*}
    \frac{1}{2}\|\vv(T)\|^2_{L^2(\Omega)}+\int_0^T 2\|\DD\vv(t)\|^2_{L^2(\Omega)}\,dt+\beta\int_{0}^{T}\|\vv_{\tau}(t)\|_{L^{2}(\Gamma)}^{2}\,dt\leq \frac{1}{2}\|\vv_0\|^2_{L^2(\Omega)}.
\end{equation*}
\end{lemma}

One of the main features of the Navier boundary conditions is the
possibility to infer properties on $\omegab=\text{curl }\vv$ at the
boundary. It is well-known~\cite{Ber2010a} that the condition $(\DD\vv\cdot\nu)\cdot\tau=0$ on a flat
boundary  is equivalent to
$\omegab\times\nu=0$. Moreover, in the case of Navier frictionless
($\beta=0$) boundary conditions on $\mathbb{R}^{3}_{+}$, it is possible to show that $-\int_{ \mathbb{R}^{3}_{+}} \Delta\omegab\cdot\omegab
 =\int_{\mathbb{R}^{3}_{+}}|\nabla\omegab|^{2}$, since the boundary
 terms vanish. In the case of a non-flat boundary a
 lower order term appears and, for a smooth domain
 $D\subset\mathbb{R}^{3}$ if we have $\omegab\times\nu=0$ on $\partial D$ (instead of Navier), then 
 \begin{equation}\label{bound}
   -\int_{D} \Delta\omegab\cdot\omegab
 \geq\int_{D}|\nabla\omegab|^{2}-c\int_{\partial D}|\omegab|^{2}, 
\end{equation} 
with $c>0$ depending on the curvature of $\partial D$, see e.g.~\cite{BB2009}.
Although in our case the domain is flat, being $\beta$ in general different from zero, we are not able to obtain an estimate like~\eqref{bound}, hence we face similar problems as in the Dirichlet case when trying to use the vorticity field. 

For these reasons in the next proposition we prove classical regularity results on~\eqref{ns}, without passing through the curl of the solution. Although the techniques used are standard, for sake of completeness in Section~\ref{sec:existence} we give the full proof; similar problems are treated e.g. in~\cite{BdV2007} and~\cite{falgaz}, but here we focus more closely on the data dependence.
\begin{proposition}
\label{existence}
	Let $T>0$, $\beta>0$, and let $\textbf{f}=\textbf{0}$, $\vv_0\in H_\P$. Then, there exists  a (global) weak solution $\vv\in L^\infty(0,T;H_\P)\cap L^2(0,T;V_\P)$ to the b.i.v.p.~\eqref{ns}  and, moreover, there is a constant $K_\Omega>0$ such that:\\
	$\bullet$  If $\vv_0\in V_\P$, then there exists 
	\begin{equation*}
	0<\frac{K_{\Omega}}{\bigg(2\|\DD \vv_0\|^2_{L^2(\Omega)}+\beta\|(\vv_0)_\tau\|^2_{L^2(\Gamma)}\bigg)^2}\le T^*\le T,
	\end{equation*}
such that the weak solution $\vv$ of~\eqref{ns} is unique in $[0,T^*)$ and
	\begin{equation}\label{regularitya1}
	\vv\in L^\infty(0,T^*;V_\P)\qquad \vv_t, \Delta \vv, \nabla q\in L^2(Q_{T^*});
	\end{equation}
 $\bullet$ If $\vv_0\in V_\P$ and \begin{equation}\label{datopiccolo}\|\vv_0\|_{L^2(\Omega)}\big(2\|\DD \vv_0\|^2_{L^2(\Omega)}+\beta\|(\vv_{0})_\tau\|_{L^2(\Gamma)}^2\big)<2 K_\Omega,
 \end{equation}
 then the solution $\vv$ of~\eqref{ns} satisfies $\vv\in L^\infty(\R^+;V_\P)$, i.e. it is unique and global in time.\\
Furthermore, for any global weak solution $\vv$ of~\eqref{ns}, there exists $\T=\T(\vv)>0$ such that
	\begin{equation}\label{regularitya2}
	\vv\in L^\infty(\T,\infty;V_\P)\qquad \qquad \vv_t, \Delta \vv, \nabla q\in L^2(\T,\infty;L^2(\Omega)).
	\end{equation}
\end{proposition}

In the sequel we are mainly interested in results of global regularity, while for the decay at infinity of the solution the reader can refer to~\cite{BFP2024} and~\cite{Kelliher2024}.
\section{Main results}
\label{sec:main-results}
In this section we present at first the results on the Stokes eigenvalue problem~\eqref{stokes}, characterizing its eigenvalues and eigenfunctions.  Then, advantaged from the spectral analysis, we exhibit some nontrivial smooth global solutions to the Navier-Stokes equations~\eqref{ns} with $\textbf{f}=\textbf{0}$.
\subsection{Spectral analysis}
The simple geometry of the domain $\Omega$ allows us to compute explicitly the eigenvalues and the eigenfunctions of the Stokes operator associated to~\eqref{stokes}. 
For $m,n\in\mathbb N$ we introduce $\mu:=\sqrt{m^2+n^2}$ and the periodic functions satisfying the boundary conditions~\eqref{per3dNSE}
\begin{equation}\label{periodic}
\begin{split}
       \mathcal{P}^u_{m,n}(x,y):=&a\cos(mx)\sin(ny)-b\sin(mx)\sin(ny)-c\sin(mx)\cos(ny)+d\cos(mx)\cos(ny)\\
        \mathcal{P}^v_{m,n}(x,y):=&a\sin(mx)\cos(ny)+b\cos(mx)\cos(ny)-c\cos(mx)\sin(ny)-d\sin(mx)\sin(ny)\\
         \mathcal{P}_{m,n}(x,y):=&a\sin(mx)\sin(ny)+b\cos(mx)\sin(ny)+c\cos(mx)\cos(ny)+d\sin(mx)\cos(ny),
\end{split}
\end{equation}
for some $a,b,c,d\in \R$, not all contemporary zero. In the spectral analysis the pressure assumes a relevant role, hence we distinguish two classes of eigenvalues with respect to the case where the pressure is constant or not.
In Section~\ref{prooftheorem3d} we prove the following theorem. 
\begin{theorem}\label{eigcube}
Let $\beta>0$. For any $m,n\in \mathbb{N}$ there exist two positive and increasing sequences of diverging eigenvalues of~\eqref{stokes} $\lambda:=\lambda_{m,n,p}(\beta)$  and $\Lambda:=\Lambda_{m,n,p}(\beta)$ ($p\in\mathbb{N}$), counted with their multiplicity, such that:
\begin{itemize}
 \item[i)] if $\nabla q\equiv0$ the eigenvalues $\lambda\in \big(\mu^2+\frac{\pi^2}{4}p^2,\mu^2+\frac{\pi^2}{4}(1+p)^2\big)$ satisfy
\begin{equation}\label{eig3d2}
  \begin{cases}
    \displaystyle\tan(\sqrt{\lambda-\mu^2}) =-\frac{\sqrt{\lambda-\mu^2}}{\beta} &\text{if } p\,\text{odd}\\
        \displaystyle\cot(\sqrt{\lambda-\mu^2}) =\frac{\sqrt{\lambda-\mu^2}}{\beta} & \text{if } p=0\,\, \text{or even};
    \end{cases}
\end{equation}
   \item[ii)] if $\nabla q\not\equiv 0$ the eigenvalues $\La\in \big(\mu^2+\frac{\pi^2}{4}(1+p)^2,\mu^2+\frac{\pi^2}{4}(2+p)^2\big)$ with $\mu>0$ satisfy
 \begin{equation}\label{eig3d}
    \begin{cases}
    \displaystyle\sqrt{\La-\mu^2}\tan(\sqrt{\La-\mu^2}) = -\frac{\La}{\beta}-\mu\tanh (\mu) &\text{if }  p=0\,\, \text{or even}\\
        \displaystyle\sqrt{\La-\mu^2}\cot(\sqrt{\La-\mu^2}) = \frac{\La}{\beta}+\mu\coth (\mu) & \text{if } p\, \text{odd}.
    \end{cases}
\end{equation}
\end{itemize}
 In both the cases the eigenfunctions are $(u_{m,n,p},v_{m,n,p},w_{m,n,p})$ where
 \begin{equation}\label{eig_strutt}
u_{m,n,p}=\mathcal U_{m,n,p}(z)\mathcal{P}^u_{m,n}(x,y), \quad v_{m,n,p}=\mathcal V_{m,n,p}(z)\mathcal{P}^v_{m,n}(x,y), \quad w_{m,n,p}= \mathcal W_{m,n,p}(z)\mathcal{P}_{m,n}(x,y),
\end{equation}
    with $\P^u_{m,n}(x,y), \P^v_{m,n}(x,y)$, and $\P_{m,n}(x,y)$ as in~\eqref{periodic}. The set of  eigenfunctions forms a basis of $H_\P$.
\end{theorem}

This theorem is a reduced version of the spectral theorem proved in Section~\ref{prooftheorem3d}; in the Appendix~\ref{eigenfunctions1} we provide the full (long) version of it, see Theorem~\ref{eigcube1}, in which many more information are given as the explicit expression of all the eigenfunctions. We point out that if we restrict to the case $m=0$ (or $n=0$) we obtain the Stokes spectrum of  the strip $\To\times (-1,1)$, i.e. the corresponding 2D domain.

    We observe that for every fixed $p$ and $\mu$, each of the equations in~\eqref{eig3d2} or in~\eqref{eig3d} admits a unique solution. Moreover, from direct computations it is easy to see that the solutions of~\eqref{eig3d2} are always different from those of~\eqref{eig3d}. This is important in order to get information on the multiplicity of the eigenvalues.   
\begin{remark}\label{rmk_mult}
   The geometric multiplicity of the eigenvalues depends on the four free parameters $a,b,c,d \in \mathbb R$ in~\eqref{periodic}; to each not vanishing parameter it corresponds at least a linear independent eigenfunction. 
    
    In particular, for each fixed $p\in \mathbb N$ only one of the following cases can occur: 
    \begin{itemize}
        \item If $m=n=0$ the eigenvalue has multiplicity $2$ and
          corresponds to the eigenfunctions
          $$
          (d\,\mathcal U_{0,0,p}(z),b\,\mathcal V_{0,0,p}(z),0),
          $$
          for some $b,d\in\mathbb{R}$ not contemporary zero; moreover, $\mathcal U_{0,0,p}(z)=\mathcal V_{0,0,p}(z)$, see Theorem~\ref{eigcube1}.
         \item If $m=0$ and $n\neq 0$ the eigenvalue has multiplicity
           $4$ and corresponds to the eigenfunctions
           $$
           \big([a\sin(ny)+d\cos(ny)]\,\mathcal
           U_{0,n,p}(z)\,,\,[b\cos(ny)-c\sin(ny)]\,\mathcal
           V_{0,n,p}(z)\,,\,[b\sin(ny)+c\cos(ny)]\mathcal
           W_{0,n,p}(z)\big),
           $$
           for some $a,b,c,d\in\mathbb{R}$ not contemporary zero. Similarly if $n=0$ and $m\neq 0$.
        \item Otherwise the periodic functions in~\eqref{periodic} will always have $4$ free parameters.
        Thus, the multiplicity of the eigenvalue  $m_\lambda$ is
        \begin{equation*}
        m_\lambda = 4\cdot\# A,
        \end{equation*}
       where $A=\{(m,n)\in \mathbb N_+\times \mathbb N_+ :\ m^2+n^2=\mu^2\}$ and $\#A$ denotes its cardinality.  
       
    By the Fermat theorem on sum of squares integers, if $\mu$ is large enough, $\#A$ can be arbitrarily large.
    \end{itemize}
\end{remark}
As a corollaries of Theorem~\ref{eigcube}, we study the spectrum of the Stokes operator in the two extremal cases corresponding to $\beta=0$, i.e. Navier frictionless boundary conditions, and $\beta\to \infty$, i.e. Dirichlet boundary conditions, see the notations in Section~\ref{sec:preliminary}. First, we state the case $\beta=0$, in which the structure of the spectrum is similar to~\cite{FG2023} and the pressure is always constant. The proof is in Section~\ref{proof:corollaryN}.
\begin{corollary}[Navier frictionless]
\label{corollary}
Let $\beta=0$. For any $m,n\in \mathbb{N}$ there exists a non-negative sequence of diverging eigenvalues  of~\eqref{stokes} $\lambda^N:=\lambda_{m,n,p}^N$ ($p\in\mathbb{N}$), counted with their multiplicity, given by
\begin{equation}\label{lambdaN}
\lambda_{m,n,p}^N=\mu^2+\frac{\pi^2}{4}p^2\qquad p\in\mathbb{N}.
\end{equation}
The corresponding pressure is always constant and the eigenfunctions are as in~\eqref{eig_strutt} with $\mathcal{U}_{m,n,p}(z)$, $\mathcal{V}_{m,n,p}(z)$, and $\mathcal W_{m,n,p}(z)$ trigonometric.
\end{corollary}
On the other hand the Dirichlet case, proved in Section~\ref{proof:corollaryD}, retraces the structure of Theorem~\ref{eigcube} where the constant or not pressure gives rise to different cases.
\begin{corollary}[Dirichlet]
\label{corollaryD}
Let $\beta \to \infty$. 
For any $m,n\in \mathbb{N}$ there exist two positive and increasing sequences of diverging eigenvalues of~\eqref{stokes} $\lambda^D:=\lambda^D_{m,n,p}$ ($p\in\mathbb{N}_+$)  and $\Lambda^D:=\Lambda^D_{m,n,p}$ ($p\in\mathbb{N}$), counted with their multiplicity, such that:
\begin{itemize}

\item[i)] if $\nabla q^D\equiv 0$ the eigenvalues are 
\begin{equation}\label{lambdaD}
\lambda_{m,n,p}^D=\mu^2+\frac{\pi^2}{4}p^2 \qquad p\in \mathbb N_+.
\end{equation}
    \item[ii)] if $\nabla q^D\not\equiv 0$ the eigenvalues $\La^D\in \big(\mu^2+\frac{\pi^2}{4}(1+p)^2,\mu^2+\frac{\pi^2}{4}(2+p)^2\big)$ with $\mu>0$ satisfy
    \begin{equation}
\label{eig3dbetainfty}
 \begin{cases}
    \displaystyle\sqrt{\La^D-\mu^2}\tan(\sqrt{\La^D-\mu^2}) = -\mu\tanh (\mu) & \text{if}\,\, p=0\,\, \text{or even}\\
        \displaystyle\sqrt{\La^D-\mu^2}\cot(\sqrt{\La^D-\mu^2}) = \mu\coth (\mu) &\text{if}\,\, p\, \text{odd}.
    \end{cases}
\end{equation}
\end{itemize}
In both the cases  the eigenfunctions are as in~\eqref{eig_strutt}.
\end{corollary}

An extended version of these corollaries is given in the Appendix~\ref{eigenfunctions1} where all the details about the eigenfunctions are given, see  respectively Corollary~\ref{corollaryN1} and~\ref{corollaryD1}.

We consider now the eigenvalues belonging to the same class with respect to the parameter $\beta$. The following result, proved in Section~\ref{proof:corollarymono}, highlights the monotone behaviour of the eigenvalues with respect to $\beta$.
\begin{corollary}
\label{cor:monotonicitybeta}
Let $m,n,p \in \mathbb N$ be fixed and let $\beta>0$. Then, the eigenvalues $\lambda_{m,n,p}= \lambda(\beta)$ and $\La_{m,n,p}= \La(\beta)$ are strictly increasing with respect to $\beta$. In particular, we have
\begin{equation*}
    0<\lambda(\beta)<\lim_{\beta\to \infty}\lambda(\beta)=\lambda^{D}\qquad\text{and}\qquad 0<\La(\beta)<\lim_{\beta\to \infty}\La(\beta)=\La^{D},
\end{equation*}
where $\lambda^D$ and $\La^D$ are given in Corollary~\ref{corollaryD}.
\end{corollary}
\begin{table}[h!]\centering
	\begin{tabular}{ |p{0.3cm}||p{0.9cm}|p{0.9cm}|p{0.4cm}||p{0.9cm}|p{0.9cm}|p{0.4cm}||p{0.9cm}|p{0.9cm}|p{0.4cm}||p{0.9cm}|p{0.95cm}|p{0.4cm}|}
		\hline
		\multicolumn{1}{|c||}{} & \multicolumn{3}{|c||}{$\beta=0$} & \multicolumn{3}{|c||}{$\beta=1$} & \multicolumn{3}{|c||}{$\beta=10$} & \multicolumn{3}{|c|}{$\beta\to \infty$} \\
		\hline
		$j$ & $\lambda^N_{m,n,p}$ &\,\, $\lambda_j^0$ & $m_\lambda$ & $\lambda_{m,n,p}$ & \,\,$\lambda_j^0$  & $m_\lambda$ & $\lambda_{m,n,p}$ &\,\,$\lambda_j^0$&$m_\lambda$  & $\lambda^D_{m,n,p}$ &\,\, $\lambda_j^0$& $m_\lambda$ 
		\\
		\hline\hline
		$1$ & $\lambda_{0,0,0}$ & \,\,\,\,\,$0$ & \,\,$2$ &$\lambda_{0,0,0}$ & \,$0.74$ & \,\,$2$ &$\lambda_{0,0,0}$ & \,$2.04$ & \,\,$2$& $\lambda_{0,0,1}$ & \,$2.47$ &\,\,$2$ 
		\\
		\hline
		$2$ & $\lambda^{\,\frown}_{1,0,0}$& \,\,\,\,\,$1$ & \,\,$8$ & $\lambda^{\,\frown}_{1,0,0}$ & \,$1.74$ & \,\,$8$ &  $\lambda^{\,\frown}_{1,0,0}$ & \,$3.04$ & \,\,$8$ & $\lambda^{\,\frown}_{1,0,1}$& \,$3.47$ & \,\,$8$ 
		\\
		\hline
		$3$ &$\lambda_{1,1,0}$ & \,\,\,\,\,$2$ & \,\,$4$ & $\lambda_{1,1,0}$ & \,$2.74$ & \,\,$4$ &  $\lambda_{1,1,0}$ & \,$4.04$ & \,\,$4$ & $\lambda_{1,1,1}$ & \,$4.47$ & \,\,$4$ 
		\\
		\hline
		$4$ & $\lambda_{0,0,1}$ & \,$2.47$ & \,\,$2$ & $\lambda_{0,0,1}$ & \,$4.12$ & \,\,$2$ &  $\lambda^{\,\frown}_{2,0,0}$ & \,$6.04$ & \,\,$8$ & $\lambda^{\,\frown}_{2,0,1}$  & \,$6.47$ & \,\,$8$
		\\
		\hline
		$5$ & \,$\lambda^{\,\frown}_{1,0,1}$ & \,$3.47$ & \,\,$8$ & $\lambda^{\,\frown}_{2,0,0}$ & \,$4.74$ & \,\,$8$ &  $\lambda^{\,\frown}_{2,1,0}$ & \,$7.04$ & \,\,$8$ & 
		$\lambda^{\,\frown}_{2,1,1}$ & \,$7.47$ &\,\,$8$ 
		\\
		\hline
		$6$ & $\lambda^{\,\frown}_{2,0,0}$ &\,\,$4$ & \,\,$8$ & $\lambda^{\,\frown}_{1,0,1}$  & \,$5.12$ & \,\,$8$ &  $\lambda_{0,0,1}$ & \,$8.20$ & \,\,$2$ & \,$\lambda_{0,0,2}$ & \,$9.87$ & \,\,$2$
		\\
		\hline
		$7$ & \,$\lambda_{1,1,1}$ & \,$4.47$ &\,\,$4$ & $\lambda^{\,\frown}_{2,1,0}$ & \,$5.74$ & \,\,$8$ &  $\lambda^{\,\frown}_{1,0,1}$  & \,$9.20$ &\,\,$8$ & $\lambda_{2,2,1}$ & \,$10.47$ & \,\,$4$ 
		\\
		\hline
		$8$ & $\lambda^{\,\frown}_{2,1,0}$ & \,\,\,\,$5$ & \,\,$8$ & \,$\lambda_{1,1,1}$ & \,$6.12$ & \,\,$4$ &  $\lambda_{2,2,0}$ & \,$10.04$ &\,\,$4$ & 
		$\lambda^{\,\frown}_{1,0,2}$ & \,$10.87$ & \,\,$8$
		\\
		\hline
		$9$ & $\lambda^{\,\frown}_{2,0,1}$  & \,$6.47$ &\,\,$8$ & $\lambda^{\,\frown}_{2,0,1}$ & \,$8.12$ & \,\,$8$ &  $\lambda_{1,1,1}$ & \,$10.20$ & \,\,$4$ & 
		$\lambda^{\,\frown}_{3,0,1}$  & \,$11.47$ & \,\,$8$ 
		\\
		\hline
		$10$ & $\lambda^{\,\frown}_{2,1,1}$  & \,$7.47$ & \,\,$8$ & $\lambda_{2,2,0}$ & \,$8.74$ & \,\,$4$ & $\lambda^{\,\frown}_{3,0,0}$  & \,$11.04$ & \,\,$8$  & $\lambda_{1,1,2}$  & \,$11.87$ &\,\,$4$ 
		\\
		\hline
	\end{tabular}
	\caption{The first 10 values of eigenvalues $\lambda_j^0(\beta)$ when $\nabla q\equiv 0$; for each $\beta$ we fix the indexes $m,n,p\in\mathbb{N}$, the eigenvalue and the geometric multiplicity $m_\lambda$. $\lambda^{\,\frown}_{m,n,p}$ means that $m$ and $n$ permute.}
	\label{tab1}\end{table}
\begin{table}[h!]\centering
	\begin{tabular}{ |p{0.3cm}||p{0.9cm}|p{0.9cm}|p{0.4cm}||p{0.9cm}|p{0.9cm}|p{0.4cm}||p{0.9cm}|p{0.95cm}|p{0.4cm}|}
		\hline
		\multicolumn{1}{|c||}{} & \multicolumn{3}{|c||}{$\beta=1$} & \multicolumn{3}{|c||}{$\beta=10$} & \multicolumn{3}{|c|}{$\beta\to \infty$} \\
		\hline
		$j$ &  $\Lambda_{m,n,p}$ & \,\,$\Lambda_j^0$  & $m_\Lambda$ & $\Lambda_{m,n,p}$ &\,\,$\Lambda_j^0$&$m_\Lambda$  & $\Lambda^D_{m,n,p}$ &\,\, $\Lambda_j^0$& $m_\Lambda$ 
		\\
		\hline\hline
		$1$  &$\Lambda^{\,\frown}_{1,0,0}$ & \,$4.65$ & \,\,$8$ &$\Lambda^{\,\frown}_{1,0,0}$ & \,$7.80$ & \,\,$8$& $\Lambda^{\,\frown}_{1,0,0}$ & \,$9.31$ &\,\,$8$ 
		\\
		\hline
		$2$ &  $\Lambda_{1,1,0}$ & \,$5.39$ & \,\,$4$ &  $\Lambda_{1,1,0}$ & \,$7.97$ & \,\,$4$ & $\Lambda_{1,1,0}$& \,$9.33$ & \,\,$4$ 
		\\
		\hline
		$3$ &$\Lambda^{\,\frown}_{2,0,0}$ & \,$7.11$ & \,\,$8$ &  $\Lambda^{\,\frown}_{2,0,0}$ & \,$9.02$ & \,\,$8$ & $\Lambda^{\,\frown}_{2,0,0}$ & \,$10.16$ & \,\,$8$ 
		\\
		\hline
		$4$ & $\Lambda^{\,\frown}_{2,1,0}$ & \,$8.03$ & \,\,$8$ &  $\Lambda^{\,\frown}_{2,1,0}$ & \,$9.72$ & \,\,$8$ & $\Lambda^{\,\frown}_{2,1,0}$  & \,$10.77$ & \,\,$8$
		\\
		\hline
		$5$ & $\Lambda_{2,2,0}$ & \,$10.87$ & \,\,$4$ &  $\Lambda_{2,2,0}$ & \,$12.16$ & \,\,$4$ & 
		$\Lambda_{2,2,0}$ & \,$13.04$ &\,\,$4$ 
		\\
		\hline
		$6$ &  $\Lambda^{\,\frown}_{3,0,0}$  & \,$11.84$ & \,\,$8$ &  $\Lambda^{\,\frown}_{3,0,0}$ & \,$13.04$ & \,\,$8$ & $\Lambda^{\,\frown}_{3,0,0}$ & \,$13.87$ & \,\,$8$
		\\
		\hline
		$7$ & $\Lambda^{\,\frown}_{1,0,1}$& \,$12.44$ & \,\,$8$ &  $\Lambda^{\,\frown}_{3,1,0}$  & \,$13.93$ &\,\,$8$ & $\Lambda^{\,\frown}_{3,1,0}$ & \,$14.73$ & \,\,$8$ 
		\\
		\hline
		$8$ &  $\Lambda^{\,\frown}_{3,1,0}$ & \,$12.81$ & \,\,$8$ &  $\Lambda^{\,\frown}_{2,3,0}$ & \,$16.69$ &\,\,$8$ & 
		$\Lambda^{\,\frown}_{3,2,0}$ & \,$17.40$ & \,\,$8$
		\\
		\hline
		$9$ & $\Lambda_{1,1,1}$ & \,$13.31$ & \,\,$4$ &  $\Lambda^{\,\frown}_{1,0,1}$ & \,$17.53$ & \,\,$8$ & 
		$\Lambda^{\,\frown}_{4,0,0}$  & \,$20.18$ & \,\,$8$ 
		\\
		\hline
		$10$ &  $\Lambda^{\,\frown}_{2,0,1}$ & \,$15.11$ & \,\,$8$ & $\Lambda_{1,1,1}$  & \,$18.07$ & \,\,$4$  & $\Lambda^{\,\frown}_{1,0,1}$  & \,$20.57$ &\,\,$8$ 
		\\
		\hline
	\end{tabular}
	\caption{The first 10 values of eigenvalues $\Lambda_j^0(\beta)$ when $\nabla q\not\equiv 0$; for each $\beta$ we fix the indexes $m,n,p\in\mathbb{N}$, the eigenvalue and the geometric multiplicity $m_\Lambda$. $\Lambda^{\,\frown}_{m,n,p}$ means that $m$ and $n$ permute.}\label{tab:tab2}
\end{table}
 \begin{figure}[h!]
	\centering
	\includegraphics[width=83mm]{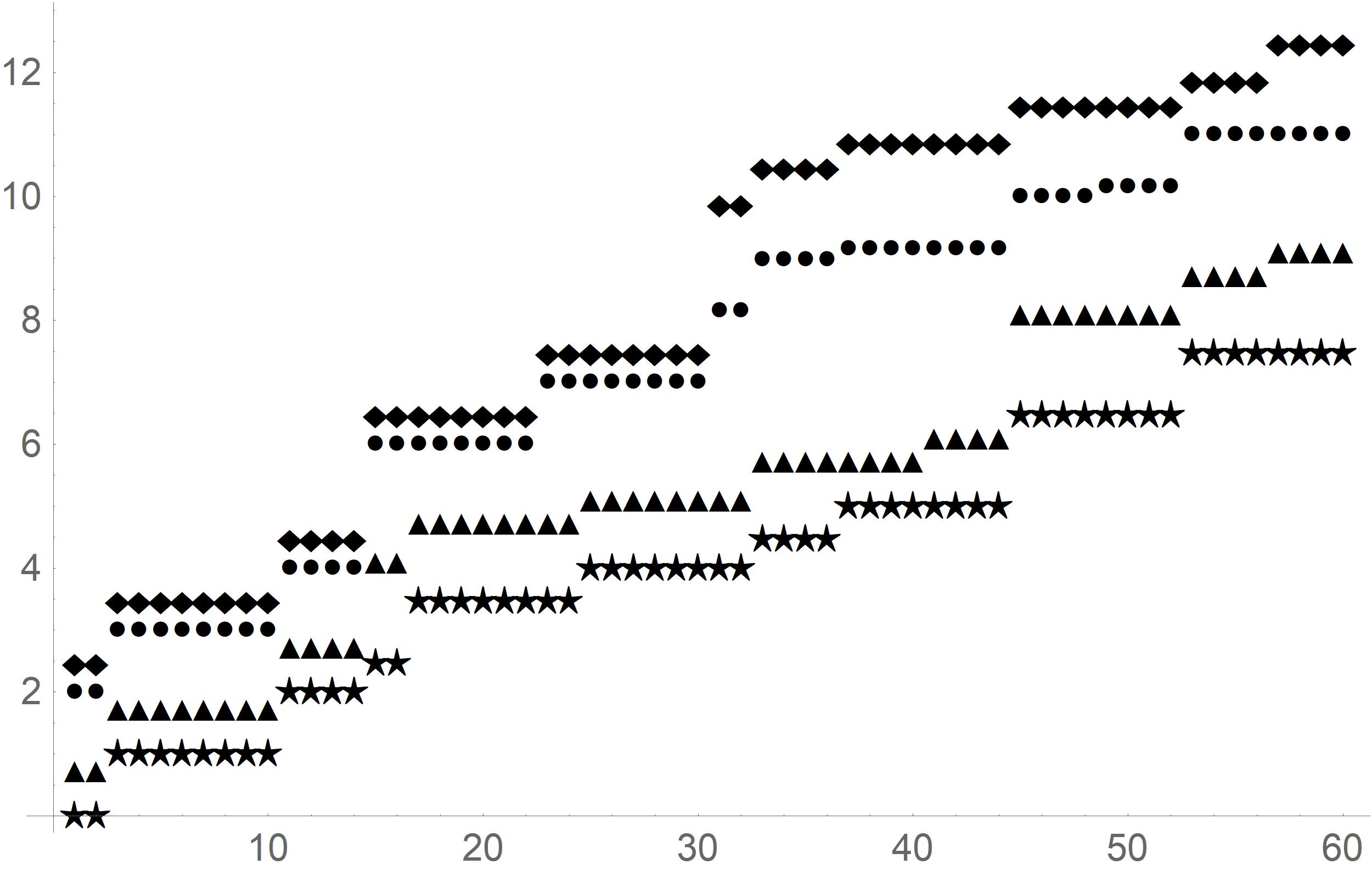}\qquad\includegraphics[width=83mm]{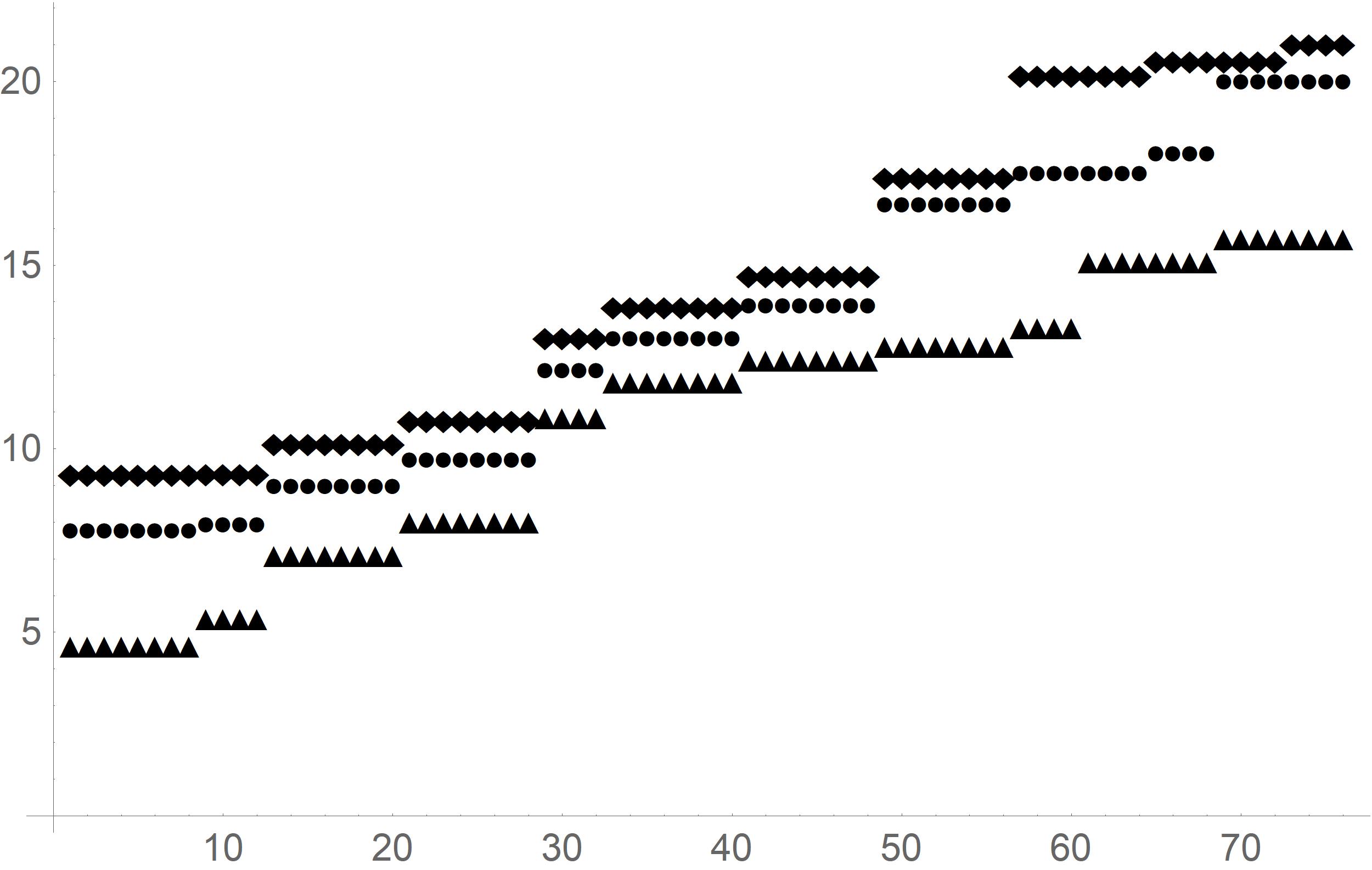}
	\caption{On the left the first eigenvalues $\lambda_k(\beta)$ when $\nabla q\equiv 0$ for 4 different $\beta$, on the right the first eigenvalues $\Lambda_k(\beta)$ when $\nabla q\not\equiv 0$ for 3 $\beta$; in the order from below stars $\beta=0$, triangles $\beta=1$, circles $\beta=10$ and squares $\beta\rightarrow\infty$.}
	\label{fig:enter-label}
\end{figure}
We complete this section computing the first few eigenvalues
of~\eqref{stokes}, represented for fixed $m,\,n,\,p\in\mathbb{N}$, and  varying
$\beta$. From Theorem~\ref{eigcube} and
Corollaries~\ref{corollary}-\ref{corollaryD} we see that for all
$\beta\geq 0$ the eigenvalues corresponding to constant pressure are
the first ones appearing in the spectrum; hence, let  begin with summarizing the information we have on the first eigenvalue, for details see Proposition~\ref{primo} in Appendix~\ref{eigenfunctions1}.
\begin{remark}\label{poincare_constants}
    For $\beta>0$ the first eigenvalue is $\lambda_0=\lambda_{0,0,0}(\beta)\in (0,\frac{\pi^2}{4})$, implying that the Poincaré constant in~\eqref{poincare}$_1$ is $C_0>\frac{4}{\pi^2}$. If $\beta\rightarrow\infty$ the constant is exactly $C_0^D=\frac{4}{\pi^2}$, while by Corollary~\ref{cor:monotonicitybeta} the constant increases as $\beta\searrow 0$. Since in the case $\beta=0$ the second eigenvalue is $1$, the constant in~\eqref{poincare}$_3$ is $C_0^N=1$.
\end{remark}
In order to consider eigenvalues with large multiplicity we need some notations.
We denote by $\lambda_k$ the \textit{eigenvalues} repeated according to their multiplicity and by $\lambda_j^0$ the {\em values of the eigenvalues}. Clearly the eigenvalues $\lambda_k$ are ordered along a non-decreasing sequence, while the values  $\lambda_j^0$ are ordered along a strictly increasing sequence.  
 In Table~\ref{tab1} (respectively Table~\ref{tab:tab2}) we report the first 10 values $\lambda_j^0$ (resp. $\Lambda_j^0$) of the eigenvalues of~\eqref{stokes} referred to the case with constant (resp. non constant) pressure and their multiplicity $m_\lambda$. According to the notations of Theorem~\ref{eigcube} and Corollaries~\ref{corollary}-\ref{corollaryD} we put in evidence the indexes $m,n,p\in\mathbb{N}$ for which the eigenvalues are obtained.
 
 When we write $\lambda^{\,\frown}_{m,n,p}$ we mean that the same
 value is obtained permuting indexes $m$ and $n$, so that the multiplicity is double.
Considering the first 10 values we find multiplicity at most equal to
8, but going further it is possible to find eigenvalues with larger
multiplicity, e.g. for $\beta=0$,  the eigenvalue corresponding to $\mu^2=325=1^2+18^2=6^2+17^2=10^2+15^2$ has multiplicity 24, see Remark~\ref{rmk_mult}.
In Figure~\ref{fig:enter-label} on the right (resp. on the left) we plot the eigenvalues $\lambda_k$ (resp. $\Lambda_k$) of the Table~\ref{tab1} (resp. Table~\ref{tab:tab2}) repeated according to their multiplicities. The plot exhibits the increasing behaviour with respect to $\beta$ proved in Corollary~\ref{cor:monotonicitybeta}.
\subsection{Global smooth solutions to the evolution problem}
In this section we exploit the spectral information obtained on the Stokes operator in order to exhibit global solutions to the corresponding Navier-Stokes problem~\eqref{ns}. To reduce the computations we focus on the case where~\eqref{stokes} has eigenfunctions corresponding to constant pressure, i.e.
\begin{equation}\label{press}
    q=Q_0\in\R;
\end{equation}
thus, when we say unique solution, we mean unique up to the constant $Q_0$.
For reader convenience we rewrite  Theorem~\ref{eigcube} (and Theorem~\ref{eigcube1}) when~\eqref{press} holds, in the following corollary.
\begin{corollary}
\label{corollary:press_cost}
    Let $\beta>0$ and assume~\eqref{press}. For any $m\in\mathbb{N}$, $n\in \mathbb{N}_+$ there exists a sequence of diverging eigenvalues $\lambda:=\lambda_{m,n,p}(\beta)\in \big(\mu^2+\frac{\pi^2}{4}p^2,\mu^2+\frac{\pi^2}{4}(1+p)^2\big)$ of~\eqref{stokes}, $p\in\mathbb{N}$,  satisfying \eqref{eig3d2}.
The corresponding eigenfunctions are 
\begin{equation}\label{eig222}
\uu_{m,n,p}(x,y,z)=\mathcal Z_{m,n,p}(z)\big(\mathcal{P}^u_{m,n}(x,y)\,,\,-\tfrac{m}{n}\mathcal{P}^v_{m,n}(x,y)\,,\,0\big),
\end{equation}
where
\begin{equation*}
    \mathcal Z_{m,n,p}(z)=\begin{cases}
        \cos(\sqrt{\lambda-\mu^2}z)\hspace{3mm} \text{if } p=0\,\, \text{or even},\vspace{3mm}\\
        \sin(\sqrt{\lambda-\mu^2}z)\hspace{3mm} \text{if } p \text{ odd},
    \end{cases}
\end{equation*}
with $\mathcal{P}^u_{m,n}$ and $\mathcal{P}^v_{m,n}$ as in~\eqref{periodic}, for some $a,b,c,d\in\R$, chosen so that $\|\uu_{m,n,p}\|_{L^2(\Omega)}=1$.
\end{corollary}
The eigenfunctions in the previous corollary deserve some comments:
first of all they always have vanishing third component, but the
dependence on $(x,y,z)$ of the first two components produces
difficulties similar to the full 3D problem. We do not include in the
analysis the case $n=0$, because we should modify the eigenfunctions,
for details see Theorem ~\ref{eigcube1}; hence, we consider
$n\in\mathbb{N}_+$ in order to not burden the notations of this
section, but (up to slight changes) it is possible to include also $n=0$.

In the next results we observe that there are particular choices of the initial datum in~\eqref{ns} for which the nonlinearity plays no role and the solution is global and smooth. We may say that such solutions are rarefied, according to the definition given in~\cite{FG2023}.

Let us begin with a simple two-dimensional result, obtained taking $m=0$ in the eigenfunctions we considered in Corollary~\ref{corollary:press_cost}. 
\begin{proposition}\label{prop-mono}
    Let $n_k\in\mathbb{N}_+$, $p_k\in\mathbb{N}$, and let
    $\uu_{0,n_k,p_k}(y,z)=\big(u_{0,n_k,p_k}(y,z),0,0\big)$ be the
    eigenfunctions in~\eqref{eig222} corresponding to the eigenvalue
    $\lambda_k:=\lambda_{0,n_k,p_k}$. Let $\textbf{f}=\textbf{0}$,
    $\vv_0(y,z)=\sum_{k=1}^\infty\gamma_k \uu_{0,n_k,p_k}(y,z)\in
    V_\P$ be the initial datum in~\eqref{ns} for some $\gamma_k\in \R$ such that $\sum_{k=1}^\infty \lambda_k\gamma_k^2<\infty$. Then, the solution $(\vv,q)$ of~\eqref{ns} is global and unique for all $t\in\R^+$, $\vv\in L^\infty(\R^+,V_\P)$ and is given by
    \begin{equation}\label{sol-mono}
        \begin{split}
            &\vv(t,y,z)=\sum_{k=1}^\infty\gamma_k\, e^{-\lambda_{k}t}\,\uu_{0,n_k,p_k}(y,z),\qquad q(t,y,z)=Q_0\in\R.
        \end{split}
    \end{equation}
\end{proposition}
The proof is given in Section~\ref{proof:prop-mono} and Proposition~\ref{prop-mono} can be properly adapted choosing in the series other combinations of eigenfunctions of~\eqref{stokes}, e.g. 
$$\uu_{m_k,0,p_k}(x,z)=\big(0,v_{m_k,0,p_k}(x,z),0\big)\quad\text{ or also }\quad \uu_{0,0,p_k}(z)=\big(u_{0,0,p_k}(z),v_{0,0,p_k}(z),0\big).$$
We point out that the vanishing of the nonlinearity in~\eqref{ns} is
peculiar of such eigenfunctions that are not fully 3D both in the space
variables and in the vector components.   
\medskip

In the next proposition, proved in
Section~\ref{proof:prop-single}, we show that for specific choices of
the coefficients in~\eqref{periodic}, it is possible to construct
global solutions starting from initial data being a single
eigenfunction, which depends on all three spatial variables.
\begin{proposition}\label{prop-single}
    Let $n,m\in\mathbb{N}_+$, $p\in\mathbb{N}$ and $\uu_{m,n,p}(x,y,z)$ be an eigenfunction in~\eqref{eig222} corresponding to the eigenvalue $\lambda:=\lambda_{m,n,p}$. Let $\textbf{f}=\textbf{0}$, $\vv_0(x,y,z)=\gamma \uu_{m,n,p}(x,y,z)$  the initial datum in~\eqref{ns} for some $\gamma\in \R$.
    If the coefficients in $\mathcal{P}^u_{m,n}$ and $\mathcal{P}^v_{m,n}$, see~\eqref{periodic}, are
\begin{equation}\label{hp}
    a=-c\,\,\text{ and }\,\, b=d \qquad \text{or alternatively}\qquad a=c\,\,\text{ and }\,\, b=-d,
\end{equation}
 then the solution $(\vv,q)$ of~\eqref{ns} is global and unique for all
 $t\in \R^+$ and is given by
    \begin{equation}\label{sol-single}
        \begin{split}
            &\vv(t,x,y,z)=\gamma\, e^{-\lambda t}\,\uu_{m,n,p}(x,y,z),\qquad q(t,x,y,z)=Q_0\in\R.
        \end{split}
    \end{equation}
\end{proposition}
In the previous propositions we were able to explicitly determine the solution (since the problem reduces to uncoupled heat equations) and, in turn, also their regularity. In the next result we assume an initial datum more general depending on 3D variables; we cannot write the explicit solution, but we show that it is globally regular, since the projection of the nonlinearity onto the eigenspace spanned by the eigenfunction vanishes. 
\begin{theorem}
\label{teo_p_nulla}
Let $m_k,n_k\in\mathbb{N}_+$, $p_k\in\mathbb{N}$, and let $
\uu_{m_k,n_k,p_k}(x,y,z)$ be the eigenfunctions in~\eqref{eig222} corresponding to the eigenvalue $\lambda_k:=\lambda_{m_k,n_k,p_k}$. Let $\textbf{f}=\textbf{0}$,
$\vv_0(x,y,z)=\sum_{k=1}^\infty\gamma_k \uu_{m_k,n_k,p_k}(x,y,z)\in
V_\P$ be the initial datum in~\eqref{ns} for some $\gamma_k\in \R$ such that $\sum_{k=1}^\infty \lambda_k\gamma_k^2<\infty$. 
If the coefficients in $\mathcal{P}^u_{m_k,n_k}$ and
$\mathcal{P}^v_{m_k,n_k}$, see~\eqref{periodic}, satisfy for all $k$
\begin{equation}
\label{hp3}
    a,b\in\mathbb{R}\,\,\text{ and }\,\,c=d=0 \qquad \text{or}\qquad   a,d\in\mathbb{R}\,\,\text{ and }\,\,c=b=0 \qquad \text{or}\qquad   b,d\in\mathbb{R}\,\,\text{ and }\,\,c=a=0,
\end{equation}
then the corresponding solution $(\vv,q)$ of~\eqref{ns}  is global,
unique for all $t\in\mathbb{R}^+$ and $\vv\in L^\infty(\mathbb{R}^+;V_\P)$.
\end{theorem}
The proof of this result will be given in Section~\ref{proof:theorem_p_nulla} and is based on particular choices of the coefficients in the eigenfunctions of Corollary~\ref{corollary:press_cost}, see \eqref{hp3}. Observe that in Theorem~\ref{teo_p_nulla} to ensure global uniqueness we can play with the coefficient $a,b,c,d$ according to~\eqref{hp3}, but in all the cases $c$ must be zero. We point out that the solution to our problem with $a=b=d=0$ and $c\neq 0$ coincides with the solution of the physical problem on the parallelepiped $(-\pi,\pi)^2\times(-1,1)$ where, instead of periodic conditions in the $x$ and $y$ directions, we put Navier frictionless boundary conditions, see~\cite{FG2023}. In this case Theorem~\ref{teo_p_nulla} fails, suggesting that in real physical situations (but in a domain with corners) the nonlinearity may have a role, that it is determinant, see also~\cite{afg} for specific analysis on this topic.

\section{Proofs of the preliminary results}\label{sec:4}
\subsection{Proof of Proposition~\ref{lem:gradsym}}\label{lem:equiv}
Let us denote $\textbf{v}=(u,v,w)$ the components of $\textbf{v}$ and since $C^\infty(\Omega)$ is dense in $H^1(\Omega)$, we can suppose these components smooth and the general result will follow by approximation arguments.  Since $|\nabla \textbf{v}|^2=|\nabla^T\textbf{v}|^2$, we get
\begin{equation*}
        |\DD\textbf{v}|^2 = \frac{1}{4}(\nabla \textbf{v}+\nabla^T\textbf{v})^2= \frac{1}{4}|\nabla \textbf{v}|^2+\frac{1}{4}|\nabla^T\textbf{v}|^2+\frac{1}{2}\nabla \textbf{v}:\nabla^T \textbf{v}= \frac{1}{2}|\nabla \textbf{v}|^2+\frac{1}{2}\nabla \textbf{v}:\nabla^T \textbf{v}.
    \end{equation*}
We have
\begin{equation*}
( \nabla\cdot \textbf{v})^2 = (u_x)^2+(v_y)^2+(w_z)^2+2u_xv_y+2u_xw_z+2v_yw_z,
\end{equation*}
so that
\begin{equation*}
    \begin{split}
   \nabla \textbf{v}:\nabla^T \textbf{v}&= (u_x)^2+(v_y)^2+(w_z)^2+ 2 u_yv_x+2u_zw_x+2w_yv_z\\
      &= (\nabla\cdot\textbf{v})^2+2(u_yv_x-u_xv_y)+2(u_zw_x-u_xw_z)+ 2(w_yv_z-v_yw_z).
            \end{split}
\end{equation*}
Now integrating over $\Omega$, we obtain
\begin{equation}
\label{Dw12nablaw}
    \|\DD \textbf{v}\|_{L^2(\Omega)}^2= \frac{1}{2}\|\nabla\textbf{v}\|_{L^2(\Omega)}^2 + \frac{1}{2}\|\nabla\cdot\textbf{v}\|_{L^2(\Omega)}^2+\int_\Omega(u_yv_x-u_xv_y)+\int_\Omega(u_zw_x-u_xw_z)+ \int_\Omega(w_yv_z-v_yw_z).
\end{equation}
Since $\vv$ is divergence-free, the thesis follows showing that the last three integrals in the previous equation are zero. 
If we denote by $\nu_k$ the $k-$th component of the outer unit normal, integrating by parts we get
\begin{equation*}
    \int_\Omega(u_yv_x-u_xv_y)= \int_{\partial \Omega}(u_yv\nu_1-u_xv\nu_2)-\int_\Omega (u_{yx}-u_{xy})v=\int_{\Gamma}(u_yv\nu_1-u_xv\nu_2),
\end{equation*}
due to Schwarz theorem and the periodic conditions in $x$ and $y$ directions. Since on $\Gamma$ we have $\nu_1=\nu_2=0$,  also the integral on $\Gamma$ vanishes.

About the remaining two integrals in~\eqref{Dw12nablaw} we integrate by parts arguing as before. Thus, since $\vv\in V_\P$
$$
\int_\Gamma(u_zw\nu_1-u_xw\nu_3)+\int_\Gamma(v_zw\nu_2-v_yw\nu_3)=-\int_{\Gamma}(u_x+v_y)\textbf{v}\cdot\nu = 0.
$$
\subsection{Proof of Proposition~\ref{lemma:regularity}}
\label{appendix:regularityevolutive}
We prove at first~\eqref{estimategrauf}. By adapting the argument
in~\cite{BDV2004} we use an approximation by artificial
compressibility. Note that the simple existence of $\uu$ weak solution
as in Proposition~\ref{prop:Stokes-existence} need the external force
$\f\in V_{P}^{*}$. Here, since we are interested in estimates which
involve also the pressure, we need to work within usual negative
Sobolev spaces (without the divergence-free constraint). Hence we
define the following space
$$
  H^1_{\mathcal{P},0}(\Omega):=\{\vv\in H^1_{\mathcal{P}}(\Omega):\, \vv\cdot\nu=0\mbox{ on }\Gamma\}\, ,
    $$
    to ensure that functions are tangential to the boundary. We equip
    $H^{1}_{\P,0}$ with  the $L^2$-norm of the gradient and we set
    $H^{-1}_{\mathcal{P}}(\Omega):=(H^1_{\mathcal{P},0}(\Omega))^{*}$.

    We consider with full details the case $\beta>0$, since the
    treatment of the case $\beta=0$ is similar and the
    difference is mainly in the determination of the non-null kernel.  Let
    $\varepsilon > 0$, and let us consider the following auxiliary
    problem \begin{equation}\label{stokesflambda} \left\{
\begin{array}{ll}
-\Delta \uu+\nabla q=\textbf{f}\quad & \mbox{in }\Omega\,\\
\varepsilon q+\nabla\cdot \uu=0\quad & \mbox{in }\Omega\,\\
\uu\cdot\nu=2(\DD\uu\cdot\nu)\cdot\tau+\beta \uu\cdot\tau=0\quad & \mbox{on }\Gamma\,.
\end{array}\right.
\end{equation}
We introduce the bilinear form $B_\varepsilon(\cdot,\cdot)$ and the linear form $\langle L,\cdot\rangle$, respectively given by
\begin{equation*}
  B_\varepsilon(\textbf{V},\Phi)= 2\int_\Omega \DD \uu : \DD \varphib
-\int_\Omega (\nabla\cdot\uu)(\nabla\cdot
    \varphib)
  -\int_{\Omega}q\,\nabla\cdot\varphib+\beta\int_{\Gamma}\uu_\tau\cdot\varphib_\tau+\varepsilon\int_{\Omega}q\,\psi+\int_{\Omega}\psi\,\nabla\cdot \uu\quad\text{and}\quad \langle L,\Phi\rangle = \langle\textbf{f},\varphib\rangle,
\end{equation*}
where $\langle\,.\,\rangle$ denotes the duality pairing between
$H^{1}_{\mathcal{P},0}(\Omega)$ and its dual, while
$\textbf{V}:=(\uu,q)$ and $\Phi:=(\varphib,\psi)$.
Since the pressure $q$ is determined up to constants, as usual we fix it asking for a zero mean value and $L^{2}_{\P,\#}(\Omega)$
denotes the space of square summable scalar functions,
$2\pi$-periodic in the first two-variables and with zero mean value. 
The weak formulation of \eqref{stokesflambda} reads as: find  $\textbf{V}\in H^1_{\P,0}(\Omega)\times L^2_{\P,\#}(\Omega)$ such that
\begin{equation}\label{weakbilinear}
    B_\varepsilon(\textbf{V},\Phi)= \langle L,\Phi\rangle, \qquad \forall\, \Phi\in H^1_{\P,0}(\Omega)\times L^2_{\P,\#}(\Omega),\ \forall\,\varepsilon >0,
\end{equation}
where now we set
$\textbf{V}=(\textbf{u}(\varepsilon),q(\varepsilon))$, to highlight
the dependence on $\varepsilon$.

Following~\cite[Theorem 1.1]{BDV2004}, it is possible to prove that,
for every fixed $\varepsilon > 0$, the bilinear form
$B_\varepsilon(\cdot,\cdot)$ is continuous and coercive on
$H^1_{\P,0}(\Omega)\times L^2_{\P,\#}(\Omega)$  and $\langle
L,\cdot\rangle$ is continuous on $H^1_{\P,0}(\Omega)$. Note in particular
that the following equality 
\begin{equation*}
  2\|\DD\uu\|_{L^{2}(\Omega)}^{2}=\|\nabla\uu\|_{L^{2}(\Omega)}^{2}+\|\nabla\cdot\uu\|_{L^{2}(\Omega)}^{2}\qquad\forall\,\uu\in H^{1}_{\mathcal{P},0}(\Omega),
\end{equation*}
follows by integration by parts since the boundary terms vanish,
cf.~Proposition~\ref{lem:gradsym}.  
 This implies the existence and uniqueness of a weak solution $(\textbf{u}(\varepsilon), q(\varepsilon))\in H^1_{\P,0}(\Omega)\times L^2_{\P,\#}(\Omega)$ to the problem~\eqref{stokesflambda}, satisfying (for some positive constant $C>0$ independent of $\varepsilon$) the following inequality
 \begin{equation}\label{weakestimatelambda0}
   \|\textbf{u}(\varepsilon)\|_{H^1(\Omega)}^{2}+
   \varepsilon\|q(\varepsilon)\|_{L^2(\Omega)}^{2}\le C
   \|\textbf{f}\|_{H^{-1}_{\P}(\Omega)}^{2}. 
\end{equation}
In order to take the limit as $\varepsilon\to0^{+}$ and to construct
solutions of the original problem, we need estimates independent of
$\varepsilon>0$. To this end, let us take $\Phi=(\varphib,0)$ with $\varphib$ of
class $C^\infty$, $2\pi$-periodic in the first two variables and
vanishing near $\Gamma$. In such a way we get 
\begin{equation*}
 \int_{\Omega} \nabla q(\varepsilon)\cdot \varphib =\int_{\Omega}(\nabla\cdot\uu) (\nabla\cdot \varphib)-2\int_{\Omega}\DD\uu:\DD
 \varphib+\langle \f, \varphib\rangle
\end{equation*}
which implies, using also the estimate on $\uu$ coming
from~\eqref{weakestimatelambda0}, that the pressure
$q(\varepsilon)\in L^{2}_{\P}(\Omega)\subset H^{-1}_{\P}(\Omega)$
satisfies as an element of $H^{-1}_{\P}(\Omega)$ the following estimate
\begin{equation*}
\|q(\varepsilon)\|_{L^{2}(\Omega)}\leq C \| \nabla q(\varepsilon)\|_{H^{-1}_{\P}(\Omega)}\leq
 C\left(\|\f\|_{H^{-1}_{\P}(\Omega)}+\|\nabla \uu\|_{L^{2}(\Omega)}\right)\leq
 C \|\f\|_{H^{-1}_{\P}(\Omega)}.
\end{equation*}
Hence, adding the square of the latter to~\eqref{weakestimatelambda0}, we get the
following estimate
 \begin{equation}\label{weakestimatelambda}
\|\textbf{u}(\varepsilon)\|_{H^1(\Omega)}^{2}+(1+\varepsilon)\|q(\varepsilon)\|_{L^2(\Omega)}^{2}\le
C \|\textbf{f}\|_{H^{-1}_{\P}(\Omega)}^{2}.
\end{equation}
Since estimate~\eqref{weakestimatelambda} is independent on
$\varepsilon$, for $\varepsilon\to 0^+$ (up to a subsequence) the
family $(\textbf{u}(\varepsilon),q(\varepsilon))$ converges weakly in
$H^{1}_{\mathcal{P},0}(\Omega)\times L^2_{\mathcal{P},\#}(\Omega)$ to some
$(\textbf{u},q)$, solving the initial Stokes problem.
The solution $(\textbf{u},q)$ is unique (recall that also the mean value of $q$ is zero), satisfies~\eqref{weakbilinear} with $\varepsilon=0$, and the estimate
\begin{equation}\label{weakestimate}
  \|\textbf{u}\|_{H^1(\Omega)}^{2}+\|q\|_{L^2(\Omega)}^{2}\le C
  \|\textbf{f}\|_{H^{-1}_{\P}(\Omega)}^{2}.
\end{equation}
Since we have a flat boundary and the estimates~\eqref{weakestimatelambda}-\eqref{weakestimate} do not depend directly on $\beta$, we proceed as in~\cite[Theorem 1.2]{BDV2005}, where the problem in the half-space with Dirichlet boundary condition is studied. The idea is to consider the incremental ratios in all the directions orthogonal to the $z$-axis and to plug them into~\eqref{weakbilinear} in order to estimate all the second order derivatives of $\textbf{u}$ (except for the last component and the derivatives in $z$-direction)  and all the first order derivatives of $q$ (except for the one in $z$-direction). Finally, by using the equation~\eqref{stokesflambda}, we bound the remaining derivatives.
 
 To be more precise, let $h\neq 0$ small enough and let us introduce
 the translation operator $\tau_h$ that translates one of the first
 two components of the argument of a generic function $w$, i.e. $w_h\equiv\tau_h
 w(x,y,z)=w(x+h,y,z)$ or alternatively $w_h\equiv\tau_h w(x,y,z)=w(x,y+h,z)$. Let us
 denote by~\eqref{weakbilinear}$_{-h}$, equation~\eqref{weakbilinear},
 where we substitute $\Phi$, by
 $\Phi_{-h}=(\tau_{-h}\varphib,\tau_{-h}\psi)$. Since  $\langle
 w,\Phi_{-h} \rangle= \langle w_h,\Phi\rangle$, subtracting
 equations~\eqref{weakbilinear} and~\eqref{weakbilinear}$_{-h}$, and
 dividing by the scalar $h>0$, we get after integration by parts
\begin{equation*}    B_\varepsilon\bigg(\frac{\textbf{V}_h-\textbf{V}}{h}, \Phi \bigg)=
    \bigg\langle L,\frac{\Phi_{-h}-\Phi}{h}\bigg\rangle, \qquad
    \forall\,\Phi\in H^1_{\P,0}(\Omega)\times L^2_{\P,\#}(\Omega), \ \varepsilon>0.
\end{equation*}
We choose the admissibile $\Phi:=(\uu(\varepsilon),q(\varepsilon))$
and observe that (cf.~\cite[App.~I]{BDV2005}, but here the proof is
more direct due to the $x$-$y$-periodicity)
$\|\f_{h}-\f\|_{H^{-1}_{\P}(\Omega)}\leq
|h|\,\|\f\|_{L^{2}_{\P}(\Omega)}$. %
We can use the estimate
\eqref{weakestimate} applied to the differences to infer that 
$\|\textbf{u}_{h}-\uu\|_{H^1(\Omega)}^{2}+\|q_{h}-q\|_{L^2(\Omega)}^{2}\le
C|h|^{2}\, \|\textbf{f}\|^{2}_{L^{2}_{\P}(\Omega)}$. Hence, using the
definition of Sobolev norm, we proved the following estimate 
\begin{equation}\label{weakestimatestar}
    \|\nabla^2_*\textbf{u}(\varepsilon)\|_{L^2(\Omega)}^{2}+\|\nabla_*q(\varepsilon)\|_{L^2(\Omega)}^{2}\le C \|\mathbf{f}\|_{L^2_{\P}(\Omega)}^{2},
\end{equation}
where $\nabla^2_*$ denotes the second order derivatives except for
$\partial^2/\partial z^2$, $\nabla_*$ denotes the gradient except
for the last component, and the constant $C>0$ does not depend on
$\varepsilon>0$. By weak convergence and uniqueness of the limit this
allows to infer the same estimates on  $(\uu,q)$.

To estimate the remaining derivatives we consider the third scalar
component of~\eqref{stokesflambda}$_1$ and we
differentiate~\eqref{stokesflambda}$_2$ by $z$. Note that in the
interior of $\Omega$ after localization (which introduces several
terms of lower 
order) a classical argument allows to show that all second order
derivatives of $\uu(\varepsilon)$ and first order derivatives of
$q(\varepsilon)$  can be estimated by means of the differential
quotients and so they belong also to $L^2(\Omega')$ for each
strip $\Omega'\subset \Omega$, with positive distance from $\Gamma$. Hence, for  a.e. $x\in\Omega$ we can
write the following pointwise identities for $\uu=(u,v,w)$ and
$\mathbf{f}=(f_{1},f_{2},f_{3})$:
\begin{equation*}
    \begin{cases}
        -\frac{\partial^2w(\varepsilon)}{\partial z^2} +
        \ \frac{\partial q(\varepsilon)}{\partial z}= \Delta_* w(\varepsilon)+{f}_3=: A\\
        \ \ \frac{\partial^2w(\varepsilon)}{\partial z^2} + \varepsilon\frac{\partial q(\varepsilon)}{\partial z}= -\frac{\partial}{\partial z}(\nabla_*\cdot \textbf{u}_*(\varepsilon))=: B,
    \end{cases}
\end{equation*}
where again the symbol ``$_{*}$'' means that we are not considering the last component and the $z$ variable in differentiating. Solving the $2\times 2$ system in the unknowns $\partial^2w(\varepsilon)/\partial z^2$ and $\partial q(\varepsilon)/\partial z$ (the determinant is $-1-\varepsilon \neq 0$, for $\varepsilon$ small enough), we have the following estimates
\begin{equation*}
    \begin{cases}
        |\frac{\partial^2w(\varepsilon)}{\partial z^2}|\le
        |B-\varepsilon A|
        \\
                |\frac{\partial q(\varepsilon)}{\partial z}| \le |A+B|.
    \end{cases}
\end{equation*}
It is clear that by~\eqref{weakestimatelambda} and~\eqref{weakestimatestar}
\begin{equation*}
    \begin{cases}
        \|A\|^2_{L^2(\Omega)}\le \|\nabla^2_* w(\varepsilon)\|_{L^2(\Omega)} +\|{f}_3\|_{L^2(\Omega)}\le C \|\textbf{f}\|_{L^2(\Omega)}\\
        \|B\|^2_{L^2(\Omega)}\le \|\nabla_*^2\textbf{u}_*(\varepsilon)\|_{L^2(\Omega)}\le C \|\textbf{f}\|_{L^2(\Omega)},
    \end{cases}
\end{equation*}
and so
\begin{equation*}
    \bigg\|\frac{\partial^2w(\varepsilon)}{\partial z^2}\bigg\|^2_{L^2(\Omega)} + \bigg\| \frac{\partial q(\varepsilon)}{\partial z}\bigg\|^2_{L^2(\Omega)}\le \Bar{C} \|\textbf{f}\|_{L^2(\Omega)}^{2},
\end{equation*}
where $\Bar{C}>0$ is again independent of $\varepsilon>0$. These bounds, together with the first two scalar equations in~\eqref{stokesflambda}$_1$, allow us to estimate the remaining derivatives. Summing up we prove \eqref{estimategrauf} letting $\varepsilon \to 0^+$.

We show the equivalence between norms. Let $\uu\in H^2(\Omega)\cap V_\P$, by definition $\widetilde\Delta \uu \in L^2(\Omega)$ if there exists $\textbf{f}\in L^2(\Omega)$ such that 
    \begin{equation}\label{weak_stokes}
         2(\DD \uu, \DD \varphib) +\beta \int_\Gamma \uu_\tau\cdot\varphib_\tau=(\textbf{f},\varphib)\, \qquad \forall\,\varphib\in V_\P\,. 
    \end{equation}
   By interpreting~\eqref{weak_stokes} as the weak formulation
   of~\eqref{stokesf} we deduce that $\uu\in H^2(\Omega)$  and that
   $$
\|\uu\|_{H^2(\Omega)}\le C \|\textbf{f}\|_{L^2(\Omega)}= C\|\widetilde\Delta \uu\|_{L^2(\Omega)} \,,
   $$
for some $C=C(\Omega)>0$. The inverse inequality is trivial, then it
follows that $\|\widetilde\Delta \uu\|_{L^2(\Omega)}$ is a norm on
$H^2(\Omega)\cap V_\P$.

To conclude, the case $\beta=0$ can be treated similarly working in $V_\P\cap Z_\P^\perp$, where $Z_\P^\perp$ is the orthogonal complement of the kernel of the Stokes operator.
\begin{remark}
    This proposition does not depend on the dimension of the space and it is valid for any $\Omega= \mathbb T^{n-1}\times(-\ell,\ell)\subset \R^n$, $n>1$, and  $\ell>0$; clearly the constants will depend on $A$ or on the size of the torus, if not $2\pi$.
\end{remark}
\subsection{Proof of Proposition~\ref{lemma0}}\label{proof:lambda0}
If $\lambda=0$ we should find $\uu\in V_\P$ nontrivial solution of 
\begin{equation}\label{weak3dzero}
2\|\DD \uu\|_{L^2(\Omega)}^2+\beta\int_{\Gamma}|\uu_\tau|^2=0.
\end{equation}
If $\DD \uu \equiv 0$, it means that
\begin{equation*}
    \begin{cases}
        u=a+dy+ez \\
        v=b-dx+fz\\
        w= c-ex-fy,
    \end{cases}
\end{equation*}
for some $a,b,c,d,e,f \in \mathbb R$, for computational details see~\cite{falocchi2022remarks}.
Since the condition $\uu\cdot \nu=0$ holds on $\Gamma$, we have
\begin{equation*}
    (a+dy+ez,b-dx+fz,c-ex-fy)\cdot (0,0,\pm 1) = \pm (c-ex-fy)= 0 \qquad \forall\,(x,y) \in \To^2.
\end{equation*}
This means that $c=e=f=0$ must be satisfied. In this way the velocity is 
\begin{equation}\label{firsteig}
    \textbf{u}= (a+dy,b-dx,0) \qquad a,b,d \in \mathbb R;
\end{equation}
since such solutions have to be periodic we infer $d=0$. If $\beta >0$, from~\eqref{weak3dzero} we have $\uu\cdot \tau=0$ on $\Gamma$, which implies $c_1a=-c_2b$, not holding for all $c_1,c_2\in\R$.
 If $\beta = 0$, then the velocity is of the kind~\eqref{firsteig} with $d=0$
 and, therefore, $\uu=(a,b,0)$ for some $a,b\in \R$ not contemporary zero is a corresponding eigenfunction.
 
 To prove \eqref{poincare} we distinguish two cases. If $\beta>0$ by using \eqref{weak} we obtain the variational characterization of the first eigenvalue of the Stokes operator as  
 \begin{equation*}
 	\lambda_0 = \min_{\textbf{0}\neq\vv \in V_\mathcal{P}} \frac{2\|\DD \vv\|^2_{L^2(\Omega)}+\beta\|\vv_\tau\|^2_{L^2(\Gamma)}}{\|\vv\|^2_{L^2(\Omega)}},
 \end{equation*}
 from which \eqref{poincare}$_1$ is proved with $C_0=\frac{1}{\lambda_0}$; similarly \eqref{poincare}$_2$ being $\vv_\tau=0$ on $\Gamma$. If $\beta=0$ we know that $\lambda_0=0$, hence we write the variational characterization of the first (non-zero) eigenvalue $\lambda_1>0$ of the Stokes operator as  
 \begin{equation*}
 	\lambda_1 = \min_{\textbf{0}\neq\vv \in V_\mathcal{P}\cap Z^{\perp}_\P} \frac{2\|\DD \vv\|^2_{L^2(\Omega)}}{\|\vv\|^2_{L^2(\Omega)}},
 \end{equation*}
getting \eqref{poincare}$_3$ with $C_0^N=\frac{1}{\lambda_1}$. Observe that the minimum is here taken among the nontrivial functions in $V_\P$, belonging to the orthogonal complement of the Stokes operator kernel.
\subsection{Proof of Proposition~\ref{existence}}\label{sec:existence}
The proof of this proposition is standard, hence we give here the main steps; for the complete version see e.g.~\cite{falgaz} or~\cite{galdi,heywood,temam} for the corresponding problem under Dirichlet boundary conditions. Throughout the proof we denote by $C$ as a positive constant, possibly changing line by line. 

In the space $V_\P$, we consider the Stokes eigenvalue problem~\eqref{stokes}.
Since $V_\P$ is a separable Hilbert space and the Stokes operator is linear, with compact inverse, self-adjoint, and positive, all the eigenvalues of~\eqref{stokes} have finite multiplicity and can be ordered in an increasing divergent sequence $\{\lambda_k\}_{k\in \mathbb{N}_+}$, in which the eigenvalues are repeated according to their multiplicity. 
Up to normalization, the set of eigenfunctions $\{\uu_k\}_{k\in \mathbb{N}_+}$ is a complete orthonormal system in $H_\P$, and complete orthogonal in $V_\P$.

We consider the eigenvectors $\{\uu_k\}_{k=1}^\infty\subset V_\P$ of~\eqref{stokes}
and the Faedo-Galerkin $n^{th}$-order approximation of~\eqref{ns}, that is, the system of ordinary differential equations
\begin{equation}\label{ode}
\begin{cases}
(\vv^n_t(t),\uu_k)-(\widetilde\Delta \vv^n(t),\uu_k)=-\big((\vv^n(t)\cdot\nabla) \vv^n(t), \uu_k\big)
\qquad k=1,\dots,n,
\\
\vv^n(0)=\vv_0^n,
\end{cases}
\end{equation}
where  $\vv_0^n:=\sum_{k=1}^n (\vv_0,\uu_k)\,\uu_k$ is the projection
in $H_\P$ of $\vv_0$ onto the space spanned by $\uu_1,\dots, \uu_n$ and $(\widetilde\Delta \vv^n,\uu_k)=-2(\DD \vv^n,\DD \uu_k)-\beta\int_{\Gamma}(\vv^n\cdot\tau)(\uu_k\cdot\tau) $.
By the theory of ode's systems,~\eqref{ode} admits a unique (local in time) solution
\begin{equation}\label{approx}
\vv^n(t,x,y,z):=\sum_{k=1}^n c^n_{k}(t)\,\uu_k(x,y,z),
\end{equation}
with $c^n_k(t)$ being smooth (scalar) coefficients. We multiply~\eqref{ode} by $c^n_k(t)$
and we sum for $k$ from 1 to $n$, obtaining
\begin{equation}\label{eq000}
\frac{d}{dt}\|\vv^n(t)\|^2_{L^2(\Omega)}+4\|\DD \vv^n(t)\|^2_{L^2(\Omega)}+2\beta\|\vv^n_\tau(t)\|^2_{L^2(\Gamma)}=0.
\end{equation}
Integrating with respect to the time, we find that for all $t\in[0,T]$
\begin{eqnarray}
\!\|\vv^n(t)\|^2_{L^2(\Omega)}+4\!\int_0^t\!\|\DD \vv^n(s)\|^2_{L^2(\Omega)}ds+2\beta\!\int_0^t\!\|\vv^n_\tau(s)\|^2_{L^2(\Gamma)}ds\leq \| v_0\|_{L^2(\Omega)}^2,\label{eq21}
\end{eqnarray}
giving a \textit{uniform bound} on $\|\vv^n(t)\|_{L^2(\Omega)}$, hence proving that the interval of existence of $c_{k}^{n}$ is arbitrary. In particular, with the a priori bound in $L^\infty(0,T;H_\P)\cap L^2(0,T;V_\P)$,
derived from~\eqref{eq21}, one can obtain (by using standard compactness arguments) the existence of a weak solution $u\in L^\infty(0,T;H_\P)\cap L^2(0,T;V_\P)$.

If the initial datum is smoother, say in $V_\P$  we can  multiply the equations in~\eqref{ode} by $\lambda_kc^n_k(t)$ and, summing for $k$ from 1 to $n$, we get
\begin{equation}\label{stima2}
\frac{1}{2}\frac{d}{dt}\big(2\|\DD \vv^n(t)\|^2_{L^2(\Omega)}+\beta\|\vv^n_\tau\|^2_{L^2(\Gamma)}\big)+ \|\widetilde\Delta \vv^n(t)\|^2_{L^2(\Omega)}=\big((\vv^n(t)\cdot\nabla) \vv^n(t),\widetilde\Delta \vv^n(t) \big),
\end{equation}
since  $\big(\vv^n_t(t),-\widetilde\Delta \vv^n(t)\big)=\frac{1}{2}\frac{d}{dt}\big(2\|\DD \vv^n(t)\|^2_{L^2(\Omega)}+\beta \|\vv^n_\tau\|^2_{L^2(\Gamma)}\big)$. 

Hence, we bound the nonlinear term in~\eqref{stima2}, using Sobolev
inequality, Poincaré inequality \eqref{poincare}$_1$, the equivalence between the norms $\|\nabla
\w\|_{L^2(\Omega)}$ and $\|\DD\w\|_{L^2(\Omega)}$,
see Proposition~\ref{lem:gradsym}, and the equivalence between the norms
$\|\w\|_{H^2(\Omega)}$ and $\|\widetilde\Delta\w\|_{L^2(\Omega)}$, see
Proposition~\ref{lemma:regularity}. More precisely, there exists $C>0$ such that
	\begin{subequations}
		\begin{align}
		&\|\w\|_{L^6(\Omega)}\leq C\big(2\|\DD \w\|^2_{L^2(\Omega)}+\beta\|\w_\tau\|^2_{L^2(\Gamma)}\big)^{1/2} \quad\hspace{21mm}  \forall\,\w\in V_\P, \label{sob1}
  \\&\|\nabla \w\|_{L^3(\Omega)}\leq  C\|\DD \w\|^{1/2}_{L^2(\Omega)}\|\widetilde\Delta \w\|^{1/2}_{L^2(\Omega)}\notag\\&\hspace{17mm}\leq  C\big(2\|\DD \w\|^2_{L^2(\Omega)}+\beta\|\w_\tau\|^2_{L^2(\Gamma)}\big)^{1/4} \|\widetilde\Delta \w\|^{1/2}_{L^2(\Omega)}\quad \hspace{1.5mm} \forall\,\w\in H^2(\Omega)\cap V_\P.\label{sob2}
		\end{align}
	\end{subequations}

Since $\|(\uu\cdot\nabla) \w\|_{L^2(\Omega)}\leq \|\uu\|_{L^6(\Omega)}\|\nabla \w\|_{L^3(\Omega)}$ for all $\uu,\w\in H^2(\Omega)\cap V_\P$, we then infer
	\begin{equation}\label{stimanonlin}
	\begin{split}
	|\big((\vv^n\cdot\nabla) \vv^n,\widetilde\Delta \vv^n\big )|&\leq \|(\vv^n\cdot\nabla) \vv^n\|_{L^2(\Omega)}\|\widetilde\Delta \vv^n\|_{L^2(\Omega)}\\&\leq C\big(2\|\DD \vv^n\|^2_{L^2(\Omega)}+\beta\|\vv^n_\tau\|^2_{L^2(\Gamma)}\big)^{3/4}\|\widetilde\Delta \vv^n\|^{3/2}_{L^2(\Omega)}\\&\leq C\big(2\|\DD \vv^n\|^2_{L^2(\Omega)}+\beta\|\vv^n_\tau\|^2_{L^2(\Gamma)}\big)^3+\dfrac{1}{2}\|\widetilde\Delta \vv^n\|^2_{L^2(\Omega)},
	\end{split}
	\end{equation}
in which we used the H\"{o}lder inequality,~\eqref{sob1}-\eqref{sob2} and the Young inequality $ab\leq \eps\frac{a^4}{4}+\frac{3}{4\eps}b^{4/3}$ for $a,b>0$ and $\eps=\tfrac{3}{2}$.\par
	
From~\eqref{stima2} we obtain
\begin{equation}\label{eqdiff}
	\frac{d}{dt}\big(2\|\DD \vv^n(t)\|^2_{L^2(\Omega)}+\beta\|\vv^n_\tau(t)\|^2_{L^2(\Gamma)}\big)\leq C\big(2\|\DD \vv^n(t)\|^2_{L^2(\Omega)}+\beta\|\vv^n_\tau(t)\|^2_{L^2(\Gamma)}\big)^3.
\end{equation}
 Setting $y(t):=2\|\DD \vv^n(t)\|^2_{L^2(\Omega)}+\beta\|\vv^n_\tau(t)\|^2_{L^2(\Gamma)}$ and, solving the differential inequality $\dot{y}(t)\leq y(t)^3$, with initial datum $y(0)=2\|\DD \vv^n_0\|^2_{L^2(\Omega)}+\beta\|(\vv_0^n)_\tau\|^2_{L^2(\Gamma)}$ we get the following inequality
\begin{equation}
\label{bd1}
2\|\DD\vv^n(t)\|^2_{L^2(\Omega)}
+\beta\|\vv^n_\tau(t)\|^2_{L^2(\Gamma)}
\leq \dfrac{1}{\sqrt{\big(2\|\DD \vv_0\|^2_{L^2(\Omega)}+\beta\|(\vv_0)_\tau\|^2_{L^2(\Gamma)}\big)^{-2}-2C t}}
:=F(t)\qquad \forall\,t\in [0,T^*),
\end{equation}
for some
\begin{equation*}
	T^*\geq\frac{K_{\Omega}}{\bigg(2\|\DD \vv_0\|^2_{L^2(\Omega)}+\beta\|(\vv_0)_\tau\|^2_{L^2(\Gamma)}\bigg)^2},
\end{equation*}
where $K_{\Omega}>0$ depends on $\Omega$.

We integrate~\eqref{stima2} from $0$ to $t\in [0,T^*)$ and, through~\eqref{stimanonlin}, we find $G(t)>0$ on $[0,T^*)$ such that
	\begin{equation}\label{bd2}
	\begin{split}
	&2\|\DD \vv^n(t)\|^2_{L^2(\Omega)}+\beta\|\vv^n_\tau(t)\|^2_{L^2(\Gamma)}+\int_0^t\|\widetilde\Delta \vv^n(s)\|^2_{L^2(\Omega)}ds\\\leq &2\|\DD \vv_0\|^2_{L^2(\Omega)}+\beta\|(\vv_0)_\tau\|^2_{L^2(\Gamma)}+C \int_0^t\big(2\|\DD \vv^n(s)\|^2_{L^2(\Omega)}+\beta\|\vv^n(s)\|_{L^2(\Gamma)}^2\big)^3ds\\
	\Rightarrow&\int_0^t\|\widetilde\Delta \vv^n(s)\|^2_{L^2(\Omega)}ds\leq 2\|\DD \vv_0\|^2_{L^2(\Omega)}+\beta\|(\vv_0)_\tau\|^2_{L^2(\Gamma)}+C \int_0^tF(s)^3ds:=G(t)\qquad \forall\,t\in [0,T^*).
	\end{split}
	\end{equation}
Then, we multiply the first equation in~\eqref{ode} by $\frac{d}{dt}c^n_k(t)$ and, summing for $k$ from 1 to $n$, we obtain
	\begin{equation*}
	\begin{split}
	\|\vv^n_t(t)\|^2_{L^2(\Omega)}=& \big(\widetilde\Delta \vv^n(t),\vv^n_t(t)\big)-\big((\vv^n(t)\cdot\nabla) \vv^n(t),\vv^n_t(t)\big).
	\end{split}
	\end{equation*}
By proceeding as for~\eqref{stimanonlin}, through H\"{o}lder and Young inequalities we have
	\begin{equation*}
	\begin{split}
	\|\vv^n_t(t)\|_{L^2(\Omega)}&\leq \|\widetilde\Delta \vv^n(t)\|_{L^2(\Omega)}+\|(\vv^n(t)\cdot\nabla) \vv^n(t)\|_{L^2(\Omega)}\\&\leq C\Big[\|\widetilde\Delta \vv^n(t)\|_{L^2(\Omega)}+\big(2\|\DD \vv^n\|^2_{L^2(\Omega)}+\beta\|\vv^n_\tau\|^2_{L^2(\Gamma)}\big)^{3/2}\Big].
	\end{split}
	\end{equation*}
	After squaring and integrating from $0$ to $t\in [0,T^*)$, we obtain
	\begin{equation}\label{bd3}
	\begin{split}
	\int_0^t\|\vv^n_t(s)\|^2_{L^2(\Omega)}ds\leq& C\bigg(\int_{0}^t\|\widetilde\Delta \vv^n(s)\|^2_{L^2(\Omega)}ds+\int_{0}^t\big(2\|\DD \vv^n\|^2_{L^2(\Omega)}+\beta\|\vv^n_\tau\|^2_{L^2(\Gamma)}\big)^{3}\bigg)\\\leq& C\bigg(G(t)+\int_{0}^t F(s)^3ds\bigg)\qquad\qquad  \forall\,t\in [0,T^*).
	\end{split}
	\end{equation}

From~\eqref{bd1}-\eqref{bd2}-\eqref{bd3} we infer the boundedness of $\vv^n$ in $L^\infty(0,T^*;V_\P)$ and the boundedness of $\widetilde\Delta \vv^n$, $\vv^n_t$ in $L^2(Q_{T^*})$. Hence, up to a subsequence, we have weak convergence in $L^2(0,T^*;H_\P)$ of $\vv_t^n$ and $\widetilde\Delta \vv^n$, respectively, to $\vv_t$ and $\widetilde\Delta \vv$.  
In standard way, see e.g.~\cite{heywood}, is possible to infer that $\vv_t-\Delta \vv+(\vv\cdot \nabla)\, \vv\in L^2(Q_{T^*})$ and by density that there exists a function $q$ with $\nabla q\in L^2(Q_{T^*})$ such that $\vv_t- \Delta \vv+(\vv\cdot \nabla) \vv=-\nabla q$; this proves that $(\vv,q)$  satisfy~\eqref{regularitya1}.

 Moreover, we can show the exponential decay of weak solution; we consider the approximate solution in~\eqref{approx} and applying~\eqref{poincare}$_1$ to~\eqref{eq000}, we obtain
\begin{equation*}
		\dfrac{d}{dt}\|\vv^n(t)\|_{L^2(\Omega)}^2+\frac{2}{C_0}\|\vv^n(t)\|_{L^2(\Omega)}^2\leq 0\quad \Rightarrow\quad \|\vv^n(t)\|_{L^2(\Omega)}^2\leq \|\vv_0\|_{L^2(\Omega)}^2\,e^{-\frac{2}{C_0}t}.
\end{equation*}
We then integrate~\eqref{eq000} over $[0,T]$, so that
	\begin{equation*}
	\|\vv^n(T)\|^2_{L^2(\Omega)}+2\int_0^T\big(2\|\DD \vv^n(t)\|^2_{L^2(\Omega)}+\beta\|\vv^n_\tau(t)\|_{L^2(\Gamma)}^2\big)\,dt=\|\vv^n(0)\|^2_{L^2(\Omega)}.
	\end{equation*}
Letting $T\rightarrow\infty$ we obtain
\begin{equation}\label{eq012}
\int_0^\infty \big(2\|\DD \vv^n(t)\|^2_{L^2(\Omega)}+\beta\|\vv^n_\tau(t)\|_{L^2(\Gamma)}^2\big)\,dt=\frac{\|\vv^n(0)\|^2_{L^2(\Omega)}}{2}.
\end{equation}
	As before, we set $y(t)=\big(2\|\DD
        \vv^n(t)\|^2_{L^2(\Omega)}+\beta\|\vv^n_\tau(t)\|_{L^2(\Gamma)}^2\big)$,
        $C=\frac{1}{K_\Omega}$, $E=\int_0^\infty y(t)\,dt$,  and we apply~\eqref{eqdiff}; then through a classical ode
        result (see e.g.~\cite[Lemma 5]{heywood}), if $E<1/\big(C
        y(0)\big)$ and $\dot{y}(t)\leq C y(t)^3$, then there exists $K>0$ such that $y(t)\leq K$ for all $t\geq0$. Thus if
	 $$E=\dfrac{\|\vv^n(0)\|^2_{L^2(\Omega)}}{2}<\dfrac{ K_\Omega}{ 2\|\DD \vv^n(0)\|^2_{L^2(\Omega)}+\beta\|\vv^n_\tau(0)\|_{L^2(\Gamma)}^2},$$
  i.e. if and only if
  $$\|\vv^n(0)\|_{L^2(\Omega)}\big(2\|\DD \vv^n(0)\|^2_{L^2(\Omega)}+\beta\|\vv^n_\tau(0)\|_{L^2(\Gamma)}^2\big)<2 K_\Omega,$$
  we infer a uniform bound for the $L^\infty(\R^+,V_\P)$ norm of $\vv^n$ (that also holds for the limit $\vv$). Since~\eqref{datopiccolo} ensures the previous inequality, the third statement is proved.

Any local weak solution $\vv$ can be globally extended to $\vv\in L^\infty(\R^+;H_\P)\cap L^2(\R^+;V_\P) $.
By the lower semicontinuity of the norm with respect to weak convergence, from~\eqref{eq012}
we infer
$$
\int_0^\infty \big(2\|\DD \vv(t)\|^2_{L^2(\Omega)}+\beta\|\vv_\tau(t)\|_{L^2(\Gamma)}^2\big)\,dt<\infty\, ,
$$
which yields the existence of $\T>0$ such that $\big(2\|\DD \vv(t)\|^2_{L^2(\Omega)}+\beta\|\vv_\tau(t)\|_{L^2(\Gamma)}^2\big)< C(\Omega,\|\vv(\T)\|_{L^2(\Omega)})$, that is,
\eqref{datopiccolo} translated at initial time $\T$; therefore,~\eqref{regularitya2} holds.

\section{Proofs of the main results}
\label{sec:proofs}
In this section we give the proofs of the main results of the paper;
further details will be  given in the Appendices.
\subsection{Proof of Theorem~\ref{eigcube}}\label{prooftheorem3d}
 We prove here Theorem~\ref{eigcube} and its corresponding expanded version, Theorem~\ref{eigcube1} in Appendix~\ref{eigenfunctions1}.
We look for  solutions to~\eqref{stokes} separating variables as follows
\begin{equation*}
\begin{split}
&u_{m,n,p}(x,y,z)=U_{m,n,p}(z)e^{i(mx+ny)}\\
&v_{m,n,p}(x,y,z)=V_{m,n,p}(z)e^{i(mx+ny)}\\
&w_{m,n,p}(x,y,z)=W_{m,n,p}(z)e^{i(mx+ny)}\\
&q_{m,n,p}(x,y,z)=Q_{m,n,p}(z)e^{i(mx+ny)}.
\end{split}
\end{equation*}

In any case the $z$-components of the solution have to solve the scalar problems
\begin{equation}\label{stokes23d}
\left\{
\begin{array}{ll}
U''(z)+(\lambda-\mu^2) U(z)=imQ(z)\quad & z\in(-1,1)\,\\
V''(z)+(\lambda-\mu^2) V(z)=inQ(z)\quad & z\in(-1,1)\,\\
W''(z)+(\lambda-\mu^2) W(z)=Q'(z)\quad & z\in(-1,1)\,\\
imU(z)+inV(z)+W'(z)=0\quad & z\in(-1,1)\,\\
Q''(z)-\mu^2Q(z)=0\quad & z\in(-1,1)\,\\
\pm U'(\pm 1)+\beta U(\pm 1)=0\,\\
\pm V'(\pm 1)+\beta V(\pm 1)=0\,\\
W(\pm 1)=0,
\end{array}\right.
\end{equation}
where we added $\Delta q=0$ in $\Omega$ and, for sake of simplicity, we omit the dependencies of the previous functions on $m, n$, and $p$.

We distinguish different cases. We stress the fact that the only condition that couples $U, V$, and $W$ comes from the null divergence $\eqref{stokes23d}_4$. In particular if $m=n=0$, then $U,V$, and $W$ are independent from each other.
\\
\\
\textbf{\large{$\star$ The case $m=n=0$ ($\mu=0$).}}\\
\noindent
\\
In this case~\eqref{stokes23d} becomes the real valued system
\begin{equation}
\label{stokesmnzero}
\left\{
\begin{array}{ll}
U''(z)+\lambda U(z)=0\quad & z\in(-1,1)\,\\
V''(z)+\lambda V(z)=0\quad & z\in(-1,1)\,\\
W''(z)+\lambda W(z)=Q'(z)\quad & z\in(-1,1)\,\\
W'(z)=0\quad & z\in(-1,1)\,\\
Q''(z)=0\quad & z\in(-1,1)\,\\
\pm U'(\pm 1)+\beta U(\pm 1)=0\,\\
\pm V'(\pm 1)+\beta V(\pm 1)=0\,\\
W(\pm 1)=0.
\end{array}\right.
\end{equation}
From~\eqref{stokesmnzero}$_4$-\eqref{stokesmnzero}$_8$ it is clear
that $W(z)\equiv 0$, while from~\eqref{stokesmnzero}$_3$
and~\eqref{stokesmnzero}$_5$ the pressure turns out to be $Q(z) \equiv Q_0$, with $Q_0\in\R$. So in this case we have the two following decoupled problems 
\begin{equation*}
\begin{cases}
U''(z)+\lambda U(z)=0\quad & z\in(-1,1)\,\\
\pm U'(\pm 1)+\beta U(\pm 1)=0  
\end{cases} \qquad\text{and} \qquad
\begin{cases}
V''(z)+\lambda V(z)=0\quad & z\in(-1,1)\,\\
\pm V'(\pm 1)+\beta V(\pm 1)=0.
\end{cases}
\end{equation*}
Solving the first problem, we obtain $U(z)=a_1\sin(\sqrt{\lambda} z)+b_1 \cos(\sqrt{\lambda}z)$  for some $a_1,b_1\in\R$ satisfying the boundary conditions, i.e.
\begin{equation*}
\begin{cases}
+a_1[\sqrt{\lambda}\cos(\sqrt{\lambda})+\beta\sin(\sqrt{\lambda})] +b_1[-\sqrt{\lambda}\sin(\sqrt{\lambda})+\beta\cos(\sqrt{\lambda})]=0 \vspace{1.5mm}\\
-a_1[\sqrt{\lambda}\cos(\sqrt{\lambda})+\beta\sin(\sqrt{\lambda})] +b_1[-\sqrt{\lambda}\sin(\sqrt{\lambda})+\beta\cos(\sqrt{\lambda})]=0.
\end{cases}
\end{equation*}
The system admits nontrivial solutions if and only if
the eigenvalues satisfy 
\begin{equation}\label{eq:m=0}
[-\sqrt{\lambda}\sin(\sqrt{\lambda})+\beta\cos(\sqrt{\lambda})][\sqrt{\lambda}\cos(\sqrt{\lambda})+\beta\sin(\sqrt{\lambda})]=0.
\end{equation}
Since $\beta>0$ we have $\sqrt{\lambda}\neq p\pi/2$, $p\in\mathbb{N}$, hence the eigenvalues solve
$$
\tan(\sqrt\lambda)=-\dfrac{\sqrt{\lambda}}{\beta} \quad\vee\quad \cot(\sqrt{\lambda})=\dfrac{\sqrt{\lambda}}{\beta},
$$
i.e.~\eqref{eig3d2} with $\mu=0$. This implies the existence of an increasing and diverging sequence of eigenvalues such that 
$$
\lambda_{0,0,p}\in\bigg(\dfrac{\pi^2}{4}p^2,\dfrac{\pi^2}{4}(p+1)^2\bigg)\qquad\text{with } p\in\mathbb N.
$$
The corresponding eigenfunctions $(u_{0,0,p},v_{0,0,p},0)$ are written explicitly in~\eqref{eig03d}.
\\ \\
\textbf{\large{$\star$ The case $m\in\mathbb N_+ \vee n\in\mathbb N_+$ ($i.e.,\ \mu>0$).}}\\
\noindent
\\
From~\eqref{stokes23d} the pressure  is of the form
\begin{equation}\label{pressure3d}
    Q(z) = c_1 \sinh (\mu z) +c_2 \cosh (\mu z),
\end{equation}
for some $c_1,c_2 \in \mathbb R$, possibly both zero (to be discussed separately).  Consequently $W(z)$ solves 
\begin{equation}
\begin{cases}\label{Wmnp}
W''(z)+(\lambda-\mu^2) W(z)=\mu [c_1 \cosh{(\mu z)}+c_2\sinh{ (\mu z)}] & \qquad z \in (-1,1)\\
W(\pm 1)=0.
\end{cases}
\end{equation}
We observe that $Q(z)$ and $W(z)$ are real functions, being $c_1$, $c_2\in\mathbb R$ and $\lambda\in\mathbb R$. On the other hand the functions $U(z)$ and $V(z)$ may be complex.

From now on we distinguish three different subcases, from which we will see that only the third case $\lambda>\mu^2$ gives nontrivial solutions.\\
\noindent
\\
\textbf{\normalsize{1) Case $\lambda= \mu^2$.}} The solution to~\eqref{Wmnp} is
\begin{equation*}
    W(z)= a_0z+b_0 + \frac{c_1}{\mu} \cosh{(\mu z)}+\frac{c_2}{\mu}\sinh{(\mu z)},
\end{equation*}
with $a_0$, $b_0\in \mathbb R$; imposing the boundary conditions~\eqref{stokes23d}$_8$ we find
\begin{equation*}
\begin{cases}
a_0 = -\frac{c_2}{\mu}\sinh{(\mu)}\\
    b_0=-\frac{c_1}{\mu}\cosh{(\mu)}.
    \end{cases}
\end{equation*}
From~\eqref{stokes23d}$_1$-\eqref{stokes23d}$_2$ we obtain
\begin{equation*}
\begin{split}
        U(z) = a_1z+b_1+\frac{c_1im}{\mu^2}\sinh{(\mu z)} +\dfrac{c_2im}{\mu^2} \cosh{(\mu z)}\quad\quad
        V(z) = a_2z+b_2+\frac{c_1in}{\mu^2}\sinh{(\mu z)} +\dfrac{c_2in}{\mu^2} \cosh{(\mu z)},
\end{split}
\end{equation*}
with $a_1$, $b_1$, $a_2$, $b_2\in\mathbb C$.
Applying the divergence condition~\eqref{stokes23d}$_4$ we find
\begin{equation}\label{divlambamu2}
\begin{cases}
na_2=-ma_1\\
nb_2=-mb_1-ic_2\dfrac{\sinh(\mu)}{\mu},
\end{cases}
\end{equation}
hence, from~\eqref{stokes23d}$_6$ we obtain the system 
\begin{equation}\label{sist3d}
    \begin{cases}
        \displaystyle +a_1(1+\beta)+b_1\beta+ c_1\frac{im}{\mu}\bigg(\cosh{(\mu)}+\frac{\beta}{\mu}\sinh(\mu)\bigg)+c_2\frac{im}{\mu}\bigg(\sinh{(\mu)}+\frac{\beta}{\mu}\cosh(\mu)\bigg)=0\vspace{1.5mm}\\
        \displaystyle -a_1(1+\beta)+b_1\beta- c_1\frac{im}{\mu}\bigg(\cosh{(\mu)}+\frac{\beta}{\mu}\sinh(\mu)\bigg)+c_2\frac{im}{\mu}\bigg(\sinh{(\mu)}+\frac{\beta}{\mu}\cosh(\mu)\bigg)=0.
    \end{cases}
\end{equation}

If $c_1=c_2=0$ there are no $\beta>0$ such that 
\begin{equation}\label{eq010}
\beta(\beta+1)=0,
\end{equation}
then we find the trivial solution.

If $c_1\neq 0 \vee c_2\neq 0$ we have to distinguish further subcases.

$\bullet$ \textit{Case} $m,n\in\mathbb N_+.$
We obtain
\begin{equation*}
    \begin{cases}
    \displaystyle a_1=- c_1\frac{im}{(1+\beta)\mu}\bigg(\cosh{(\mu)}+\frac{\beta}{\mu}\sinh(\mu)\bigg)\vspace{1.5mm}\\
        \displaystyle b_1=-c_2\frac{im}{\beta\mu}\bigg(\sinh{(\mu)}+\frac{\beta}{\mu}\cosh(\mu)\bigg),
    \end{cases}
\end{equation*}
and
\begin{equation*}
\begin{split}
        V(z) =\dfrac{i}{n\mu}\bigg\{ &c_1\bigg[\dfrac{m^2}{(1+\beta)}\bigg(\cosh(\mu)+\frac{\beta}{\mu}\sinh(\mu)\bigg)z+\frac{n^2}{\mu}\sinh(\mu z)\bigg]+\\&c_2\bigg[\dfrac{m^2}{\beta}\bigg(\sinh(\mu)+\frac{\beta}{\mu}\cosh(\mu)\bigg)-\sinh(\mu) +\dfrac{n^2}{\mu} \cosh(\mu z)\bigg]\bigg\}.
\end{split}
\end{equation*}
From the boundary condition~\eqref{stokes23d}$_7$ we get
\begin{equation}\label{sist3}
    \begin{cases}
    \displaystyle +c_1[\mu^2\cosh(\mu)+\beta\mu\sinh(\mu)]+c_2[(\mu^2-\beta)\sinh(\mu)+\beta\mu\cosh(\mu)]=0
     \vspace{.2cm}\\
        \displaystyle -c_1[\mu^2\cosh(\mu)+\beta\mu\sinh(\mu)]+c_2[(\mu^2-\beta)\sinh(\mu)+\beta\mu\cosh(\mu)]=0.
    \end{cases}
\end{equation}
Nontrivial solutions of~\eqref{sist3} can be found for  $m,n\in\mathbb{N}_+$ and $\beta>0$ if and only if
\begin{equation}\label{eq2}
    \big[\mu\cosh(\mu)+\beta\sinh(\mu)\big]\big[(\mu^2-\beta)\sinh(\mu)+\beta\mu\cosh(\mu)\big]=0;
\end{equation}
we observe that $\mu\cosh(\mu) +\beta\sinh(\mu)>0$ for all $\mu,\beta>0$ and $h(\mu):=(\mu^2-\beta)\sinh(\mu)+\beta\mu\cosh(\mu)$, $h'(\mu)=\mu^2\cosh(\mu)+(2+\beta)\mu\sinh(\mu)>0$ for all $\mu,\beta>0$, giving $h(\mu)>0$ for all $\mu,\beta>0$ since $h(0)=0$. This implies that there are no solutions of~\eqref{eq2}.

We observe that the cases $c_1=0$ and $c_2\neq 0$ or $c_1\neq0$ and $c_2=0$ lead to discuss only one of the factors in~\eqref{eq2}, yielding to the same conclusion.

$\bullet$ \textit{Case} $m=0$, $ n\in\mathbb{N}_+.$ The
system~\eqref{sist3d} gives nontrivial solutions if only
if~\eqref{eq010} holds; this is not possible for $\beta>0$. Hence $U(z)\equiv 0$ and we obtain~\eqref{sist3} with $\mu=n$, giving trivial solution.

$\bullet$ \textit{Case} $n=0$, $ m\in\mathbb{N}_+.$  We find the same result as before as long as the first and the second components of the velocity field are exchanged.\\ \\
\noindent
\textbf{2) Case $0<\lambda < \mu^2$.} Solving~\eqref{Wmnp} we  get
\begin{equation*}
    W(z) = a_0 \sinh (\sqrt{\mu^2-\lambda} z) + b_0\cosh (\sqrt{\mu^2-\lambda} z)+\frac{c_1\mu}{\lambda} \cosh(\mu z)+\frac{c_2\mu}{\lambda}\sinh(\mu z),
\end{equation*}
with $a_0$, $b_0\in\mathbb R$.
The boundary conditions~\eqref{stokes23d}$_8$ yield
\begin{equation*}
  \begin{cases}
a_0 = -\frac{c_2\mu}{\lambda}\frac{\sinh(\mu)}{\sinh(\sqrt{\mu^2-\lambda})}\\
    b_0=-\frac{c_1\mu}{\lambda}\frac{\cosh(\mu)}{\cosh(\sqrt{\mu^2-\lambda})},
    \end{cases}
\end{equation*}
being $\lambda\neq 0$ and $\lambda\neq\mu^2$.
From~\eqref{stokes23d}$_1$-\eqref{stokes23d}$_2$ we get
\begin{equation*}
\begin{split}
        U(z) = a_1\sinh(\sqrt{\mu^2-\lambda}z)+b_1\cosh(\sqrt{\mu^2-\lambda}z)+\frac{c_1im}{\lambda}\sinh(\mu z) +\dfrac{c_2im}{\lambda} \cosh(\mu z)\\
        V(z) = a_2\sinh(\sqrt{\mu^2-\lambda}z)+b_2\cosh(\sqrt{\mu^2-\lambda}z)+\frac{c_1in}{\lambda}\sinh(\mu z) +\dfrac{c_2in}{\lambda} \cosh(\mu z),
\end{split}
\end{equation*}
with $a_1$, $b_1$, $a_2$, $b_2\in\mathbb C$.
Imposing the divergence-free condition~\eqref{stokes23d}$_4$ we find
\begin{equation}\label{div2}
\begin{cases}
na_2=-ma_1-i\dfrac{c_1\mu}{\lambda}\dfrac{\cosh(\mu)}{\cosh(\sqrt{\mu^2-\lambda})}\sqrt{\mu^2-\lambda}\\
nb_2=-mb_1-i\dfrac{c_2\mu}{\lambda}\dfrac{\sinh(\mu)}{\sinh(\sqrt{\mu^2-\lambda})}\sqrt{\mu^2-\lambda},
\end{cases}
\end{equation}
hence, from~\eqref{stokes23d}$_6$ we obtain the system 
\begin{equation}\label{sist23d}
    \begin{cases}
        \displaystyle a_1[\sqrt{\mu^2-\lambda}\cosh(\sqrt{\mu^2-\lambda})+\beta\sinh(\sqrt{\mu^2-\lambda})]+b_1[\sqrt{\mu^2-\lambda}\sinh(\sqrt{\mu^2-\lambda})+\beta\cosh(\sqrt{\mu^2-\lambda})]+
        \\ c_1\frac{im}{\lambda}\big[\mu\cosh(\mu)+\beta\sinh(\mu)\big]+c_2\frac{im}{\lambda}\big[\mu\sinh(\mu)+\beta\cosh(\mu)\big]=0
        \vspace{.3cm}
        \\
        \displaystyle -a_1[\sqrt{\mu^2-\lambda}\cosh(\sqrt{\mu^2-\lambda}+\beta\sinh(\sqrt{\mu^2-\lambda})]+b_1[\sqrt{\mu^2-\lambda}\sinh(\sqrt{\mu^2-\lambda})+\beta\cosh(\sqrt{\mu^2-\lambda})]+
        \\ -c_1\frac{im}{\lambda}\big[\mu\cosh(\mu)+\beta\sinh(\mu)\big]+c_2\frac{im}{\lambda}\big[\mu\sinh(\mu)+\beta\cosh(\mu)\big]=0.
    \end{cases}
\end{equation}

If $c_1=c_2=0$ there are no $\lambda\in(0,\mu^2)$ such that 
\begin{equation}\label{eq00}
    \big[\sqrt{\mu^2-\lambda}+\beta\tanh(\sqrt{\mu^2-\lambda})\big]\big[\sqrt{\mu^2-\lambda}\tanh(\sqrt{\mu^2-\lambda})+\beta\big]=0,
\end{equation}
being both factors strictly positive. 

If $c_1\neq 0 \vee c_2\neq 0$ we have to distinguish further subcases.

$\bullet$ \textit{Case} $m,n\in\mathbb N_+.$
We obtain
\begin{equation*}
    \begin{cases}
    \displaystyle a_1=- c_1\frac{im\big(\mu\cosh(\mu)+\beta\sinh(\mu)\big)}{\lambda[\sqrt{\mu^2-\lambda}\cosh(\sqrt{\mu^2-\lambda})+\beta\sinh(\sqrt{\mu^2-\lambda})]}
    \vspace{.2cm}
    \\
        \displaystyle b_1=-c_2\frac{im\big(\mu\sinh(\mu)+\beta\cosh(\mu)\big)}{\lambda[\sqrt{\mu^2-\lambda}\sinh(\sqrt{\mu^2-\lambda})+\beta\cosh(\sqrt{\mu^2-\lambda})]},
    \end{cases}
\end{equation*}
where the denominators are strictly positive.
From the boundary condition~\eqref{stokes23d}$_7$ we get
\begin{equation}\label{sist3d2}
    \begin{cases}
        +c_1[\cosh(\mu) (\lambda - \beta \sqrt{\mu^2-\lambda}\tanh(\sqrt{\mu^2-\lambda}))+\beta \mu \sinh(\mu)]+\\ \,\,\,\,\,c_2[\sinh(\mu) (\lambda - \beta \sqrt{\mu^2-\lambda}\coth(\sqrt{\mu^2-\lambda}))+\beta \mu \cosh(\mu)]=0
        \vspace{.2cm}
        \\
        -c_1[\cosh(\mu) (\lambda - \beta \sqrt{\mu^2-\lambda}\tanh(\sqrt{\mu^2-\lambda}))+\beta \mu \sinh(\mu)]+ \\\,\,\,\,\,c_2[\sinh(\mu) (\lambda - \beta \sqrt{\mu^2-\lambda}\coth(\sqrt{\mu^2-\lambda}))+\beta \mu \cosh(\mu)]=0,
    \end{cases}
\end{equation}
\normalsize
having nontrivial solutions if and only if
\begin{equation}\label{eqcaso2}
\begin{split}
    &\big[\cosh(\mu) (\lambda - \beta\sqrt{\mu^2-\lambda}\tanh(\sqrt{\mu^2-\lambda}))+\beta \mu \sinh(\mu)\big]\cdot\\&\cdot\big[\sinh(\mu) (\lambda - \beta \sqrt{\mu^2-\lambda}\coth(\sqrt{\mu^2-\lambda}) )+\beta\mu \cosh(\mu)\big]=0. 
\end{split}
\end{equation}
Let $f_\mu(\lambda)=\cosh(\mu) [\lambda - \beta \sqrt{\mu^2-\lambda}\tanh(\sqrt{\mu^2-\lambda})]+\beta \mu \sinh(\mu)$ for $\lambda\in(0,\mu^2)$; computing 
$$
f_\mu'(\lambda)=\cosh(\mu)\bigg(1+\dfrac{\beta\tanh(\sqrt{\mu^2-\lambda})}{2\sqrt{\mu^2-\lambda}}+\dfrac{\beta}{2[\cosh{(\sqrt{\mu^2-\lambda)}}]^2}\bigg),
$$
we observe $f_\mu'(\lambda)>0$ for all $\lambda\in(0,\mu^2)$. \\
Similarly, we consider $g_\mu(\lambda)=\sinh(\mu) \big[\lambda - \beta \sqrt{\mu^2-\lambda}\coth(\sqrt{\mu^2-\lambda})\big]+\beta \mu \cosh(\mu)$ for $\lambda\in(0,\mu^2)$; computing 
$$
g_\mu'(\lambda)=\sinh(\mu) \bigg(1+\dfrac{\beta\coth(\sqrt{\mu^2-\lambda})}{2\sqrt{\mu^2-\lambda}}-\dfrac{\beta}{2(\sinh{\sqrt{\mu^2-\lambda}})^2}\bigg),
$$
and, observing that $\sinh(x)\geq x$ and $\coth{(x)}\ge 1$ for all $x\geq0$, we obtain $g_\mu'(\lambda)>0$ for all $\lambda\in(0,\mu^2)$. Since 
$
\lim\limits_{\lambda\rightarrow 0^+}f_\mu(\lambda)=\lim\limits_{\lambda\rightarrow 0^+}g_\mu(\lambda)=0,
$
we infer that $f_\mu(\lambda)$ and $g_\mu(\lambda)$ have no zeros for $\lambda\in(0,\mu^2)$, so that~\eqref{eqcaso2}
does not admit solutions. 

The cases $c_1=0$ and $c_2\neq 0$ or $c_1\neq0$ and $c_2=0$ lead to discuss only one of the factors in~\eqref{eqcaso2}, yielding to the same conclusion.

$\bullet$ \textit{Case} $m=0$, $ n\in\mathbb{N}_+.$ The system~\eqref{sist23d} gives nontrivial solutions if only if~\eqref{eq00} holds; not possible. Hence $U(z)\equiv 0$ and we obtain~\eqref{sist3d2} with $\mu=n$, giving trivial solution.

$\bullet$ \textit{Case} $n=0$, $ m\in\mathbb{N}_+.$  We find the same result as before as long as the first and the second components of the velocity field are exchanged.\\ \\
\noindent
\textbf{3) Case $\lambda > \mu^2$.}  The solution to~\eqref{Wmnp} is
\begin{equation}\label{WW}
    W(z) = a_0 \sin(\sqrt{\lambda-\mu^2}z)+ b_0\cos(\sqrt{\lambda-\mu^2}z) + \frac{c_1\mu}{\lambda}\cosh (\mu z)+ \frac{c_2\mu}{\lambda} \sinh (\mu z),
\end{equation}
with $a_0$, $b_0\in\mathbb R$.
The boundary conditions~\eqref{stokes23d}$_8$  yield to the system
\begin{equation}\label{sist333}
    \begin{cases}
      \displaystyle  +a_0 \sin (\sqrt{\lambda-\mu^2})+b_0\cos (\sqrt{\lambda-\mu^2})+\frac{c_1\mu}{\lambda}\cosh(\mu)+\frac{c_2\mu}{\lambda}\sinh(\mu)=0
      \vspace{.2cm}
      \\
       \displaystyle -a_0 \sin (\sqrt{\lambda-\mu^2})+b_0\cos (\sqrt{\lambda-\mu^2})+\frac{c_1\mu}{\lambda}\cosh(\mu)-\frac{c_2\mu}{\lambda}\sinh(\mu)=0,
    \end{cases}
\end{equation}
that to be discussed needs the following lemma.
\begin{lemma}\label{lemma}
    Let $\beta, \mu>0$ and $p\in\mathbb N_+$, then $\lambda=\mu^2+p^2\dfrac{\pi^2}{4}$ is not an eigenvalue of~\eqref{stokes}.
\end{lemma}
\begin{proof}
If $\lambda=\mu^2+p^2\frac{\pi^2}{4}$ with $p$ even the system~\eqref{sist333} may have nontrivial solutions only if $c_1=c_2=0$ or $c_1\neq 0\wedge c_2=0$; note that if $c_1,c_2\neq 0$ and $c_1=0\wedge c_2\neq0$ the system~\eqref{sist333} does not admit solution. 

If $c_1=c_2=0$ we have $W(z)=a_0\sin(\frac{p\pi}{2} z)$ ($p$ even) and $Q(z)\equiv 0$.
Hence, we compute $U(z)=a_1\sin(\frac{p\pi}{2} z)+b_1\cos(\frac{p\pi}{2} z)$, with $p$ even, $a_1$, $b_1\in\mathbb C$ (similarly $V(z)$) and applying~\eqref{stokes23d}$_6$ we find
$$
\begin{cases}
    +a_1 p\frac{\pi}{2}+ b_1\beta=0\\
    -a_1p\frac{\pi}{2}+ b_1\beta =0
\end{cases}\qquad \iff\qquad a_1=b_1=0.
$$
Therefore, through the divergence condition~\eqref{stokes23d}$_4$ we
obtain $a_0=0$ and  only the trivial solution.

If $c_1\neq 0\wedge c_2=0$ we find $W(z)=a_0\sin(\frac{p\pi}{2} z)-\frac{c_1\mu}{\mu^2+p^2\pi^2/4}\frac{\cosh(\mu)}{(-1)^{p/2}}\cos(\frac{p\pi}{2} z)+\frac{c_1\mu}{\mu^2+p^2\pi^2/4}\cosh (\mu z)$ and $Q(z)=c_1\sinh(\mu z)$; computing $U(z)$, $V(z)$, we do not find $a_1,b_1,a_2,b_2\in\mathbb{C}$ compatible with the boundary conditions~\eqref{stokes23d}$_{6,7}$.
Similarly if $\lambda=\mu^2+p\dfrac{\pi^2}{4}$  with $p$ odd, where we should discuss only the cases  $c_1=c_2=0$ and $c_1=0\wedge c_2\neq0$. 
\end{proof}
Thanks to Lemma~\ref{lemma} we infer
\begin{equation*}
    \begin{cases}
a_0 = -\frac{c_2\mu}{\lambda}\frac{\sinh(\mu)}{\sin(\sqrt{\lambda-\mu^2})}
\vspace{.2cm}\\
    b_0=-\frac{c_1\mu}{\lambda}\frac{\cosh(\mu)}{\cos(\sqrt{\lambda-\mu^2})}.
    \end{cases}
\end{equation*}
From~\eqref{stokes23d}$_1$-\eqref{stokes23d}$_2$ we get
\begin{equation*}
\begin{split}
        U(z) = a_1\sin(\sqrt{\lambda-\mu^2}z)+b_1\cos(\sqrt{\lambda-\mu^2}z)+ \frac{c_1im}{\lambda}\sinh(\mu z) +\dfrac{c_2im}{\lambda} \cosh(\mu z)\\
        V(z) = a_2\sin(\sqrt{\lambda-\mu^2}z)+b_2\cos(\sqrt{\lambda-\mu^2}z)+\frac{c_1in}{\lambda}\sinh(\mu z) +\dfrac{c_2in}{\lambda} \cosh(\mu z),
\end{split}
\end{equation*}
with $a_1$, $b_1$, $a_2$, $b_2\in\mathbb C$.
Imposing the divergence-free condition~\eqref{stokes23d}$_4$ we find
\begin{equation}\label{div}
\begin{cases}
na_2=-ma_1+i\dfrac{c_1\mu}{\lambda}\dfrac{\cosh(\mu)}{\cos(\sqrt{\lambda-\mu^2})}\sqrt{\lambda-\mu^2}
\vspace{.2cm}
\\
nb_2=-mb_1-i\dfrac{c_2\mu}{\lambda}\dfrac{\sinh(\mu)}{\sin(\sqrt{\lambda-\mu^2})}\sqrt{\lambda-\mu^2},
\end{cases}
\end{equation}
hence, through~\eqref{stokes23d}$_6$, we obtain the system 
\begin{equation}\label{sist33}
    \begin{cases}
        \displaystyle a_1\big[\sqrt{\lambda-\mu^2}\cos(\sqrt{\lambda-\mu^2})+\beta\sin(\sqrt{\lambda-\mu^2})\big]+b_1\big[-\sqrt{\lambda-\mu^2}\sin(\sqrt{\lambda-\mu^2})+\beta\cos(\sqrt{\lambda-\mu^2})\big]+\\ c_1\frac{im}{\lambda}\big[\mu\cosh(\mu)+\beta\sinh(\mu)\big]+c_2\frac{im}{\lambda}\big[\mu\sinh(\mu)+\beta\cosh(\mu)\big]=0
        \vspace{.2cm}
        \\
        \displaystyle -a_1\big[\sqrt{\lambda-\mu^2}\cos(\sqrt{\lambda-\mu^2})+\beta\sin(\sqrt{\lambda-\mu^2})\big]+b_1\big[-\sqrt{\lambda-\mu^2}\sin(\sqrt{\lambda-\mu^2})+\beta\cos(\sqrt{\lambda-\mu^2})\big]+\\ -c_1\frac{im}{\lambda}\big[\mu\cosh(\mu)+\beta\sinh(\mu)\big]+c_2\frac{im}{\lambda}\big[\mu\sinh(\mu)+\beta\cosh(\mu)\big]=0.
    \end{cases}
\end{equation}
The system gives different results with respect to the possibility that $c_1$ and $c_2$ are zero, hence we distinguish these cases as in the statement of the theorem.

$i)$ \underline{$c_1=c_2=0$, i.e. $q=\nabla q\equiv 0$.}\\
\\
We always obtain $W(z)\equiv 0$ and $a_1$, $b_1$, $a_2$, $b_2\in\mathbb R$. The eigenvalues $\lambda>\mu^2$ exist if and only if 
\begin{equation}\label{eq020}
   [\sqrt{\lambda-\mu^2}\cos(\sqrt{\lambda-\mu^2})+\beta\sin(\sqrt{\lambda-\mu^2})][-\sqrt{\lambda-\mu^2}\sin(\sqrt{\lambda-\mu^2})+\beta\cos(\sqrt{\lambda-\mu^2})]=0.
\end{equation}
 Thanks to Lemma~\ref{lemma} to solve~\eqref{eq020} is equivalent to find the solutions of~\eqref{eig3d2}.
Therefore, we obtain two increasing and positively divergent sequences of eigenvalues 
such that
\begin{equation*}
	\lambda^{(1)}_{m,n,p}:=\lambda^{(1)}_{m,n,p}(\beta) \in \bigg(m^2+n^2+\pi^2\bigg(\frac{1}{2}+p\bigg)^2,\,\,m^2+n^2+\pi^2(1+p)^2\bigg)\, \qquad \forall m,n \in \mathbb N_+,\,\,p\in \mathbb N
\end{equation*}
and
\begin{equation*}
	\lambda^{(2)}_{m,n,p}:=\lambda^{(2)}_{m,n,p}(\beta) \in \bigg(m^2+n^2+\pi^2p^2\,,\,m^2+n^2+\pi^2\bigg(\frac{1}{2}+p\bigg)^2\bigg)\, \qquad \forall m,n \in \mathbb N_+,\,\,p\in \mathbb N.
\end{equation*}
so that, combined together, the spectrum is given by 
\begin{equation*}
  \lambda_{m,n,p}\in \bigg(\mu^2+\frac{\pi^2}{4}p^2\,,\,\mu^2+\frac{\pi^2}{4}(1+p)^2\bigg) \qquad \forall\,p\in \mathbb N.
\end{equation*}
Different things may happen with respect to $m,n$.

$\bullet$ \textit{Case} $m,n\in\mathbb N_+.$
We obtain 
\begin{equation*}
    U_{m,n,p}(z)\!=\!\begin{cases}
        b_1\cos(\sqrt{\lambda-\mu^2}z)\hspace{2mm} \text{if } p=0\,\, \text{or even}\vspace{3mm}\\
        a_1\sin(\sqrt{\lambda-\mu^2}z)\hspace{2mm} \text{if } p \text{ odd}
    \end{cases}\qquad
        V_{m,n,p}(z)\!=\!-\dfrac{m}{n}\!\begin{cases}
         b_1\cos(\sqrt{\lambda-\mu^2}z)\hspace{2mm} \text{if } p=0\,\, \text{or even}\vspace{3mm}\\
        a_1\sin(\sqrt{\lambda-\mu^2}z)\hspace{2mm} \text{if } p \text{ odd},
    \end{cases}
\end{equation*}
and, in turn~\eqref{eig_pnullaappendix} with $a_1=b_1=1$.

$\bullet$ \textit{Case} $m=0$, $ n\in\mathbb{N}_+.$ From~\eqref{div} we find $a_2=b_2=0$, so that $V(z)\equiv 0$, while 
\begin{equation}\label{U(z)}
    U(z)=\begin{cases}
        b_1\cos(\sqrt{\lambda-n^2}z)\hspace{3mm} \text{if }  p=0\,\, \text{or even}\vspace{3mm}\\
        a_1\sin(\sqrt{\lambda-n^2}z)\hspace{3mm} \text{if } p \text{ odd},
    \end{cases}
\end{equation}
and, in turn,~\eqref{eig_pnullaappendix} with $m=0$, $\mu=n$ and $a_1=b_1=1$.

$\bullet$ \textit{Case} $n=0$, $ m\in\mathbb{N}_+.$ We obtain the same result as before as long as the first and the second components of the velocity field are exchanged, i.e.
$a_1=b_1=0$, so that $U(z)\equiv 0$, while 
\begin{equation}\label{V(z)}
    V(z)=\begin{cases}
        b_2\cos(\sqrt{\lambda-m^2}z)\hspace{3mm} \text{if }  p=0\,\, \text{or even}\vspace{3mm}\\
        a_2\sin(\sqrt{\lambda-m^2}z)\hspace{3mm} \text{if } p \text{ odd}.
    \end{cases}
\end{equation}
and, in turn,~\eqref{eigm0pappendix} with $a_2=b_2=1$.

$ii)$ \underline{$c_1\neq 0 \vee c_2\neq 0$, i.e. $\nabla q\not\equiv 0$.}\\
\\
We distinguish the following subcases. 

$\bullet$ \textit{Case} $m,n \in\mathbb N_+.$
We observe that~\eqref{eq020} does not hold. More precisely, if  $\tan(\sqrt{\Lambda-\mu^2})=-\sqrt{\Lambda-\mu^2}/\beta$, the system~\eqref{sist33} may have solutions only if $c_1=0\wedge c_2\neq0$ (and not if $c_1,c_2\neq 0$ and $c_1\neq0\wedge c_2=0$); in this case, the boundary conditions~\eqref{stokes23d}$_7$ are satisfied if and only if 
$$
\beta\mu\cosh(\mu)+\Lambda\sinh(\mu)+\beta^2\sinh(\mu)=0,
$$
that is not possible, being all the terms strictly positive. Similarly if $\cot(\sqrt{\Lambda-\mu^2})=\sqrt{\Lambda-\mu^2}/\beta$.

Therefore, we obtain
\begin{equation*}
    \begin{cases}
    \displaystyle a_1=- c_1\frac{im\big[\mu\cosh(\mu)+\beta\sinh(\mu)\big]}{\Lambda[\sqrt{\Lambda-\mu^2}\cos(\sqrt{\Lambda-\mu^2})+\beta\sin(\sqrt{\Lambda-\mu^2})]}
    \vspace{.2cm}
    \\
        \displaystyle b_1=-c_2\frac{im\big[\mu\sinh(\mu)+\beta\cosh(\mu)\big]}{\Lambda[-\sqrt{\Lambda-\mu^2}\sin(\sqrt{\Lambda-\mu^2})+\beta\cos(\sqrt{\Lambda-\mu^2})]}.
    \end{cases}
\end{equation*}
From the boundary condition~\eqref{stokes23d}$_7$ we get
\begin{equation*}
    \begin{cases}
        +c_1\big[\cosh(\mu) (\Lambda + \beta \sqrt{\Lambda-\mu^2}\tan (\sqrt{\Lambda-\mu^2}))+\beta \mu \sinh(\mu)\big]+ \\\,\,\,\,c_2\big[\sinh(\mu) (\Lambda - \beta \sqrt{\Lambda-\mu^2}\cot( \sqrt{\Lambda-\mu^2}))+\beta \mu \cosh(\mu)\big]=0
        \vspace{.2cm}
        \\
        -c_1\big[\cosh(\mu) (\Lambda + \beta \sqrt{\Lambda-\mu^2}\tan(\sqrt{\Lambda-\mu^2}))+\beta \mu \sinh(\mu)\big]+ 
        \\
        \,\,\,\,c_2\big[\sinh(\mu) (\Lambda - \beta \sqrt{\Lambda-\mu^2}\cot (\sqrt{\Lambda-\mu^2}))+\beta \mu \cosh(\mu)\big]=0,
    \end{cases}
\end{equation*}
having nontrivial solutions if and only if
\begin{equation*}
 \big[\cosh(\mu) (\Lambda + \beta \sqrt{\Lambda-\mu^2}\tan (\sqrt{\Lambda-\mu^2}))+\beta \mu \sinh(\mu)\big]\big[\sinh(\mu) (\Lambda - \beta \sqrt{\Lambda-\mu^2}\cot (\sqrt{\Lambda-\mu^2}))+\beta \mu \cosh(\mu)\big]=0.
\end{equation*}
This is equivalent to find solutions of one of the two equations in~\eqref{eig3d} for some $\La>\mu^2$. 
Hence, we obtain two increasing and positively divergent sequences of eigenvalues 
\begin{equation*}
	\Lambda^{(1)}_{m,n,p}:=\Lambda^{(1)}_{m,n,p}(\beta) \in \bigg(m^2+n^2+\pi^2\bigg(\frac{1}{2}+p\bigg)^2,\,\,m^2+n^2+\pi^2(1+p)^2\bigg)\, \qquad \forall m,n \in \mathbb N_+,\,\,p\in \mathbb N.
\end{equation*}
and
\begin{equation*}
	\Lambda^{(2)}_{m,n,p}:=\Lambda^{(2)}_{m,n,p}(\beta) \in \bigg(m^2+n^2+\pi^2(1+p)^2\,,\,m^2+n^2+\pi^2\bigg(\frac{3}{2}+p\bigg)^2\bigg)\, \qquad \forall m,n \in \mathbb N_+,\,\,p\in \mathbb N.
\end{equation*}
so that, combined together, the spectrum is given by
\begin{equation*}
  \La_{m,n,p}\in \bigg(\mu^2+\frac{\pi^2}{4}(1+p)^2\,,\,\mu^2+\frac{\pi^2}{4}(2+p)^2\bigg) \qquad \forall\, p\in \mathbb N.
\end{equation*}

We compute the corresponding eigenfunctions, observing that $U(z)$ and $V(z)$ are purely imaginary; then, we put 
$\mathcal{U}(z):=\frac{U(z)}{i}$ and $\mathcal{V}(z):=\frac{V(z)}{i}$ where
{\small
\begin{equation*}
\begin{split}
    &\mathcal U(z)=\frac{m}{\Lambda}\cdot\\&\begin{cases}
         c_1\dfrac{[\beta\sin(\sqrt{\La-\mu^2})+\sqrt{\La-\mu^2}\cos(\sqrt{\La-\mu^2})]\sinh(\mu z)-[\mu\cosh(\mu)+\beta\sinh(\mu)]\sin(\sqrt{\La-\mu^2}z)}{\beta\sin(\sqrt{\La-\mu^2})+\sqrt{\La-\mu^2}\cos(\sqrt{\La-\mu^2})}\hspace{2mm} \text{if } p=0\,\, \text{or even}\\
 c_2\dfrac{[\beta\cos(\sqrt{\La-\mu^2})-\sqrt{\La-\mu^2}\sin(\sqrt{\La-\mu^2})]\cosh(\mu z)-[\mu\sinh(\mu)+\beta\cosh(\mu)]\cos(\sqrt{\La-\mu^2}z)}{\beta\cos\sqrt{\La-\mu^2}-\sqrt{\La-\mu^2}\sin(\sqrt{\La-\mu^2})}\hspace{2mm} \text{if } p \text{ odd},
    \end{cases}\\
      &\mathcal V(z)\!=\!\begin{cases}
       \dfrac{c_1}{\La}\bigg[ n
     \sinh(\mu z)+\\\dfrac{[\mu(\La-n^2)\cosh(\mu)+\beta(m^2\sinh(\mu)+\mu\sqrt{\La-\mu^2}\cosh(\mu)\tan(\sqrt{\La-\mu^2}))]\sin(\sqrt{\La-\mu^2}z)}{n(\beta\sin(\sqrt{\La-\mu^2})+\sqrt{\La-\mu^2}\cos(\sqrt{\La-\mu^2}))}\bigg]\hspace{1mm} \text{if } p=0\, \text{or even}
     \vspace{.2cm}
     \\
 \dfrac{c_2}{\La}\bigg[n\cosh(\mu z)+\\\dfrac{[\mu(\La-n^2)\sinh(\mu)+\beta(m^2\cosh(\mu)-\mu\sqrt{\La-\mu^2}\sinh(\mu)\cot(\sqrt{\La-\mu^2}))]\cos(\sqrt{\La-\mu^2}z)}{n(\beta\cos(\sqrt{\La-\mu^2})-\sqrt{\La-\mu^2}\sin(\sqrt{\La-\mu^2}))}\bigg]\hspace{1mm} \text{if } p \text{ odd}.
    \end{cases}
    \end{split}
\end{equation*}
}
For coherence, we also put $\mathcal{W}(z):=W(z)$, so that
\begin{equation*}
   \mathcal W(z)=\frac{\mu}{\Lambda}\begin{cases}c_1
        \dfrac{\cos(\sqrt{\La-\mu^2})\cosh(\mu z)-\cosh(\mu)\cos(\sqrt{\La-\mu^2}z)}{\cos(\sqrt{\La-\mu^2})}\hspace{3mm} \text{if } p=0\,\, \text{or even}
        \vspace{.2cm}\\
        c_2\dfrac{\sin(\sqrt{\La-\mu^2})\sinh(\mu z)-\sinh(\mu)\sin(\sqrt{\La-\mu^2}z)}{\sin(\sqrt{\La-\mu^2})}\hspace{5mm} \text{if } p \text{ odd},
    \end{cases}
\end{equation*}
\normalsize
associating the pressure
\begin{equation*}
   Q(z)=\begin{cases}
        c_1\sinh(\mu z)\qquad \text{if } p=0\,\, \text{or even}\\
        c_2\cosh(\mu z)\qquad \text{if } p \text{ odd}.
    \end{cases}
\end{equation*}
Choosing e.g. $c_1=c_2=\La$, we find~\eqref{eigvappendix} with the pressure in~\eqref{press2appendix}.

We observe that the cases $c_1=0$ and $c_2\neq 0$ or $c_1\neq0$ and $c_2=0$ lead to discuss only one among~\eqref{eig3d}, yielding to the same conclusion. 

$\bullet$ \textit{Case} $m=0$, $ n\in\mathbb{N}_+.$ 
We get the decoupling of the $U(z)$ component that solves
$$
\begin{cases}
    U''(z)+(\Lambda-n^2)U(z)=0\quad z\in(-1,1)\\
    \pm U'(\pm 1)+\beta U(\pm1)=0.
\end{cases}
$$
Therefore two scenarios are possible: $U(z)$ is as in~\eqref{U(z)} or $U(z)\equiv 0$.

In the first case the eigenvalues have to solve~\eqref{eq020} with $\mu=n$.
 If  $\tan(\sqrt{\Lambda-n^2})=-\sqrt{\Lambda-n^2}/\beta$,  the boundary conditions~\eqref{stokes23d}$_7$ read 
\begin{equation*}
    \begin{cases}
        +c_1n[\beta  \sinh(n) +n\cosh(n)]+c_2[(\Lambda+\beta^2)\sinh(n) +\beta n \cosh(n)]=0\\
        -c_1n[\beta  \sinh(n) +n\cosh(n)]+c_2[(\Lambda+\beta^2)\sinh(n) +\beta n \cosh(n)]=0,
    \end{cases}
\end{equation*}
that is not possible for $c_1\neq 0 \vee c_2\neq 0$. Similarly if $\cot(\sqrt{\Lambda-n^2})=\sqrt{\Lambda-n^2}/\beta$.

Hence, it remains to discuss the case $U(z)\equiv 0$. Repeating similar computations as when $m,n\in\mathbb N_+$, the eigenvalues satisfy~\eqref{eig3d} with $\mu=n$ and the eigenfunctions are given in~\eqref{eigvappendix} with the pressure in~\eqref{press2appendix} taking $m=0$ and $\mu=n$.

$\bullet$ \textit{Case} $n=0$, $ m\in\mathbb{N}_+.$ We observe the decoupling of the $V(z)$ component that solves
$$
\begin{cases}
    V''(z)+(\Lambda-m^2)V(z)=0\quad z\in(-1,1)\\
    \pm V'(\pm 1)+\beta V(\pm1)=0.
\end{cases}
$$
Therefore two scenarios are possible: either $V(z)$ is as in~\eqref{V(z)} or $V(z)\equiv 0$.

In the first case the eigenvalues have to solve~\eqref{eq020} with $\mu=m$.
 If  $\tan(\sqrt{\Lambda-m^2})=-\sqrt{\Lambda-m^2}/\beta$,  the boundary conditions~\eqref{stokes23d}$_6$ read 
\begin{equation*}
    \begin{cases}
        +c_1m[\beta  \sinh(m) +m\cosh(m)]+c_2[(\Lambda+\beta^2)\sinh(m) +\beta m \cosh(m)]=0\\
        -c_1m[\beta  \sinh(m) +m\cosh(m)]+c_2[(\Lambda+\beta^2)\sinh(m) +\beta m\cosh(m)]=0,
    \end{cases}
\end{equation*}
that is not possible for $c_1\neq 0 \vee c_2\neq 0$. Similarly if $\cot(\sqrt{\Lambda-m^2})=\sqrt{\Lambda-m^2}/\beta$.

In the case $V(z)\equiv 0$ we have as eigenfunctions $(\mathcal U_{m,0,p}(z)\mathcal{P}^u_{m,0}(x,y)\,,\,0\,,\,\mathcal W_{m,0,p}(z)\mathcal{P}_{m,0}(x,y))$, where the $z-$component can be inferred by~\eqref{eigvappendix} with $n=0$ and $\mu=m$, see~\eqref{eigvappendix_m0p}. The corresponding pressure $q_{m,0,p}$ comes from~\eqref{press2appendix} with $n=0$ and $\mu=m$, see~\eqref{press.m0p.appendix}.

Since we are interested in real eigenfunctions we combine the exponential components as follows. First of all we observe that problem~\eqref{stokes23d} is solved by four families of linear independent solutions: 
\begin{equation*}
    \begin{split}
        \uu_{1}&=\big(+U_{m,n,p}(z)e^{+i(mx+ny)},+V_{m,n,p}(z)e^{+i(mx+ny)},W_{m,n,p}(z)e^{+i(mx+ny)}\big)\qquad q_{1}=Q_{m,n,p}(z)e^{ +i(mx+ny)}\\
\uu_{2}&=\big(- U_{m,n,p}(z)e^{-i(mx+ny)},- V_{m,n,p}(z)e^{- i(mx+ny)},W_{m,n,p}(z)e^{- i(mx+ny)}\big)\qquad q_{2}=Q_{m,n,p}(z)e^{- i(mx+ny)}\\
\uu_{3}&=\big( +U_{m,n,p}(z)e^{+i(mx- ny)},- V_{m,n,p}(z)e^{+i(mx- ny)},W_{m,n,p}(z)e^{+i(mx- ny)}\big)\qquad q_{3}=Q_{m,n,p}(z)e^{ +i(mx- ny)}\\
\uu_{4}&=\big(- U_{m,n,p}(z)e^{-i(mx- ny)},+ V_{m,n,p}(z)e^{- i(mx- ny)},W_{m,n,p}(z)e^{- i(mx- ny)}\big)\qquad q_{4}=Q_{m,n,p}(z)e^{- i(mx- ny)}.
    \end{split}
\end{equation*}
If $m,n\in\mathbb N_+$, recalling that $\mathcal{U}(z)=\frac{U(z)}{i}$, $\mathcal{V}(z)=\frac{V(z)}{i}$ and $\mathcal{W}(z)=W(z)$ are real functions, we compute
\begin{equation*}
\begin{split}
    \widetilde \uu_1=\dfrac{\uu_1+\uu_2}{2}&=\big(-\mathcal{U}(z)\sin(mx+ny),-\mathcal{V}(z)\sin(mx+ny),+\mathcal W(z)\cos(mx+ny)\big)\\
    \widetilde \uu_2=\dfrac{\uu_1-\uu_2}{2i}&=\big(+\mathcal{U}(z)\cos(mx+ny),+\mathcal{V}(z)\cos(mx+ny),+\mathcal W(z)\sin(mx+ny)\big)\\
    \widetilde \uu_3=\dfrac{\uu_3+\uu_4}{2}&=\big(-\mathcal{U}(z)\sin(mx-ny),+\mathcal{V}(z)\sin(mx-ny),+\mathcal W(z)\cos(mx-ny)\big)\\
   \widetilde \uu_4 =\dfrac{\uu_3-\uu_4}{2i}&=\big(+\mathcal{U}(z)\cos(mx-ny),-\mathcal{V}(z)\cos(mx-ny),+\mathcal W(z)\sin(mx-ny)\big).
    \end{split}
\end{equation*}
Using trigonometric formula we can combine them as follows
\begin{equation*}
\begin{split}
 \dfrac{\widetilde\uu_1+\widetilde\uu_3}{2}&=\big(-\mathcal{U}(z)\sin(mx)\cos(ny),-\mathcal{V}(z)\cos(mx)\sin(ny),+\mathcal W(z)\cos(mx)\cos(ny)\big)\\
 \dfrac{\widetilde \uu_1-\widetilde\uu_3}{2}&=\big(+\mathcal{U}(z)\cos(mx)\sin(ny),+\mathcal{V}(z)\sin(mx)\cos(ny),+\mathcal W(z)\sin(mx)\sin(ny)\big)\\
\dfrac{\widetilde\uu_2+\widetilde\uu_4}{2}&=\big(+\mathcal{U}(z)\cos(mx)\cos(ny),-\mathcal{V}(z)\sin(mx)\sin(ny),+\mathcal W(z)\sin(mx)\cos(ny)\big)\\
   -\dfrac{\widetilde\uu_2-\widetilde\uu_4}{2}&=\big(-\mathcal{U}(z)\sin(mx)\sin(ny),+\mathcal{V}(z)\cos(mx)\cos(ny),+\mathcal W(z)\cos(mx)\sin(ny)\big),
    \end{split}
\end{equation*}
getting, in turn, the linear combination functions in~\eqref{periodic}.
\\ \\
We conclude the proof showing that the set of eigenfunctions forms a basis of $H_\P$.
We begin with dividing the eigenfunctions into two classes depending on the
associated pressure. We denote by $\{\uu\}$ and $\{\uu^q\}$ the eigenfunctions with constant and non constant pressure, respectively.

It is straightforward to see that, they are orthogonal. We prove that $\mathcal{X}_{\P}:=\text{span }\{\textbf{u},\textbf{u}^q\}$ 
 is complete in $H_\P$. To do this we show that any function $\textbf{f}\in
 H_{\P}$ which is orthogonal to every function in $\mathcal{X}_{\P}$ must be identically zero in
 $\Omega$. By the periodicity conditions in $x$ and $y$, we can write $\textbf{f}=(f_1,f_2,f_3)$ as follows
    \begin{equation*}
        f_k(x,y,z)= \sum_{m,n\in \mathbb N} A_{m,n}^{(k)}(z)\sin{(mx+ny)}+B_{m,n}^{(k)}(z)\cos{(mx+ny)}, \qquad k=1,2,3,
    \end{equation*}
    for some scalar functions $A_{m,n}^{(k)}, B_{m,n}^{(k)}: (-1,1)\to \mathbb R$, $k=1,2,3$.  The divergence condition implies that for every $m,n \in \mathbb N$
    \begin{equation}\label{divcondf1}
        mA_{m,n}^{(1)}(z)+nA_{m,n}^{(2)}(z) + \frac{d}{dz}{A_{m,n}^{(3)}}(z) =0\qquad z\in (-1,1),
    \end{equation}
    and
    \begin{equation}\label{divcondf2}
        mB_{m,n}^{(1)}(z)+nB_{m,n}^{(2)}(z) + \frac{d}{dz}{B_{m,n}^{(3)}}(z) =0\qquad z\in (-1,1).
    \end{equation}
     When the pressure is constant, the eigenfunctions are of the type 
     $$
     (\mathcal{U}_{m,n,p}(z)\mathcal{P}^u_{m,n}(x,y),\mathcal{V}_{m,n,p}(z)\mathcal{P}^u_{m,n}(x,y),0),
     $$ see Corollary~\ref{corollary:press_cost}. In particular, the $z$-components solve~\eqref{stokes23d}$_1$-\eqref{stokes23d}$_2$, and, by the Sturm-Liouville theory, we have that these eigenfunctions form a complete orthogonal system in $L^2(-1,1)$.  Since the variables are separated, the orthogonality condition plays a role only in the $z$-component, which means that 
    \begin{equation*}
        \int_{-1}^{1} A_{m,n}^{(k)}(z)\,\mathcal{U}_{m,n,p}(z)\,dz = \int_{-1}^{1} B_{m,n}^{(k)}(z)\,\mathcal{V}_{m,n,p}(z)\,dz =0 \qquad k=1,2.
    \end{equation*}
    By  the completeness in $L^2(-1,1)$ we get that $A_{m,n}^{(k)}=B_{m,n}^{(k)}=0$ a.e. in $(-1,1)$, with $k=1,2$, then the first two components of $\textbf{f}$ are zero in $\Omega$. Hence from~\eqref{divcondf1}-\eqref{divcondf2}, we infer $A_{m,n}^{(3)}=a_{m,n}$ and $B_{m,n}^{(3)}=b_{m,n}$ constants.\\
    Let us now consider the scalar product with every function in $\uu^q$, where only the third component of the eigenfunctions has a role in this case. When $p$ is odd, $\mathcal{W}_{m,n,p}$ defined in~\eqref{eigvappendix} is a linear combination of trigonometric and hyperbolic sine functions, so that the orthogonality condition is immediately satisfied. If $p$ is zero or even,
    $\mathcal{W}_{m,n,p}$ is a linear combination of trigonometric and hyperbolic cosine functions so that the orthogonality condition becomes
    \begin{equation*}
        \cos(\sqrt{\Lambda-\mu^2})\int_{-1}^1 \cosh (\mu z)\,dz = \cosh(\mu) \int_{-1}^1\cos(\sqrt{\Lambda-\mu^2})\,dz.
    \end{equation*}
    Integrating and rearranging the term we get
    \begin{equation*}
        \frac{\tan(\sqrt{\Lambda-\mu^2})}{\sqrt{\Lambda-\mu^2}}= \frac{\tanh(\mu)}{\mu}.
    \end{equation*}
    Since~\eqref{eig3d} holds true, we have
    \begin{equation*}
        \frac{1}{\Lambda-\mu^2}\bigg(-\frac{\Lambda}{\beta}-\mu\tanh(\mu)\bigg)= \frac{\tanh(\mu)}{\mu},
    \end{equation*}
    which is not possible for $\Lambda>\mu^2$, $\mu>0$, since the left and right hand side have opposite signs. Therefore, the only possibility is that $a_{m,n}=b_{m,n}=0$ for every $m,n \in \mathbb N$, which means that $f_3$ must be zero,  implying $\textbf{f}=\textbf{0}$.

\subsection{Proof of Corollary~\ref{corollary}}\label{proof:corollaryN}
If $\beta=0$ the boundary conditions in~\eqref{stokes} become $w=u_{z}=v_{z}=0$ on $\Gamma$, since we have a flat boundary; hence, writing the differential equations up to the boundary we find
$q_z=w_{zz}$ on $\Gamma$, being $w=0$ on $\Gamma$ and, consequently, $w_{xx}=w_{yy}=0$ on $\Gamma$. Applying the divergence condition up to $\Gamma$ and using again the conditions  $u_{z}=v_{z}=0$ on $\Gamma$, we end up that $w_{zz}=0$ on $\Gamma$, implying
$$
q_z=0\quad \text{on }\Gamma\quad \Rightarrow \quad (Q^N)'(\pm 1)=0.
$$
Therefore, from~\eqref{stokes23d} we get
\begin{equation*}
\left\{
\begin{array}{ll}
(Q^N)''(z)-\mu^2Q^N(z)=0\quad & z\in(-1,1)\,\\
(Q^N)'(\pm 1)=0,
\end{array}\right.
\end{equation*}
giving $Q^N(z)$ constant if $\mu=0$ and $Q^N(z)\equiv 0$ if $\mu>0$.
This produces some simplifications with respect to the case $\beta>0$, indeed~\eqref{stokes23d} uncouples into
\begin{equation}\label{decoup}
    \begin{cases}
        (U^N)''(z)+(\lambda-\mu^2)U^N(z)=0\\
        (U^N)'(\pm1)=0
    \end{cases} \begin{cases}
        (V^N)''(z)+(\lambda-\mu^2)V^N(z)=0\\
        (V^N)'(\pm1)=0
    \end{cases}
      \begin{cases}
        (W^N)''(z)+(\lambda-\mu^2)W^N(z)=0\\
        W^N(\pm1)=0,
    \end{cases}
\end{equation}
with $imU^N(z)+inV^N(z)+(W^N)'(z)=0$ for all $z\in(-1,1)$.

If $m=n=0$ ($\mu=0$), from equation~\eqref{eq:m=0} with $\beta=0$ we find that the eigenvalues satisfy
$\sin({2\sqrt{\lambda^N}})=0$, i.e. we obtain the increasing and diverging sequence of eigenvalues 
\begin{equation*}
    \lambda_{0,0,p}^N = \frac{\pi^2}{4}p^2 \qquad p \in \mathbb N.
\end{equation*}
The corresponding pressure is $q^N_{0,0,p}\equiv Q^N_0$ with $Q_0^N\in \mathbb R$ and the eigenfunctions $(u_{0,0,p}^N,v_{0,0,p}^N,0)$ are written explicitly in~\eqref{eig_0_N}. In this case we also have that the first eigenvalue of the spectrum is $\lambda_{0,0,0}^N=0$, see Proposition~\ref{lemma0}.

If $m\in\mathbb{N}_+\vee n\in\mathbb{N}_+$ ($\mu>0$) from~\eqref{decoup} it is easy to see that in the case $0<\lambda<\mu^2$ we have only the trivial solution.

If $\lambda= \mu^2$, from~\eqref{decoup} we have $W^N(z)\equiv0$,  while  $U^N(z)$ and $V^N(z)$ are constants satisfying the divergence condition. If $m,n\in\mathbb{N}_+$ or $m=0, n\in\mathbb{N}_+$ the eigenfunctions are given in~\eqref{final_N} (case $p=0$), while if $m\in\mathbb{N}_+, n=0$ the eigenfunctions are in~\eqref{final_N2} (case $p=0$).

If $\lambda>\mu^2$ we solve~\eqref{decoup}$_3$ finding 
$
W^N(z)=a_0\sin(\sqrt{\lambda^N-\mu^2}z)+b_0\cos(\sqrt{\lambda^N-\mu^2}z)
$
with $a_0$, $b_0\in\mathbb R$.
The conditions $W^N(\pm1)=0$ yield to the system
\begin{equation*}
    \begin{cases}
      \displaystyle  +a_0 \sin (\sqrt{\lambda^N-\mu^2})+b_0\cos (\sqrt{\lambda^N-\mu^2})=0
      \vspace{.2cm}
      \\
       \displaystyle -a_0 \sin (\sqrt{\lambda^N-\mu^2})+b_0\cos (\sqrt{\lambda^N-\mu^2})=0,
    \end{cases}
\end{equation*}
having nontrivial solutions if and only if
\begin{equation*}
   \sin({2\sqrt{\lambda^N-\mu^2}})=0\quad \iff\quad \lambda_{m,n,p}^N=\mu^2+\frac{\pi^2}{4}p^2\qquad p\in\mathbb{N}_+.
\end{equation*}
The corresponding eigenfunctions are 
$$
 W_{m,n,p}^N(z)=\begin{cases}
a_0\sin(\frac{\pi}{2}p z)\qquad \text{if } p \text{ even}\\
b_0\cos(\frac{\pi}{2}p z)\qquad \text{if }p \text{ odd}
\end{cases}\qquad a_0,b_0\in\mathbb{R}.
$$
From~\eqref{decoup}$_1$-\eqref{decoup}$_2$  we find  
\begin{equation*}
    U_{m,n,p}^N(z)=i\begin{cases}
b_1\cos(\frac{\pi}{2}p z)\hspace{5mm} \text{if } p \text{ even}\\
a_1\sin(\frac{\pi}{2}p z)\hspace{5mm} \text{if }p \text{ odd}
\end{cases}\quad\quad V_{m,n,p}^N(z)=i\begin{cases}
b_2\cos(\frac{\pi}{2}p z)\hspace{5mm} \text{if } p \text{ even}\\
a_2\sin(\frac{\pi}{2}p z)\hspace{5mm} \,\text{if }p \text{ odd},
\end{cases}
\end{equation*}
with $a_1,b_1,a_2,b_2\in \mathbb R$, where by convenience, instead of considering the usual real solutions combination we consider the imaginary one. Through the divergence condition we have
\begin{equation}\label{divN}
\begin{cases}
nb_2=-mb_1+\frac{\pi p}{2}a_0\qquad \text{if } p \text{ even}\\
na_2=-ma_1-\frac{\pi p}{2}b_0\qquad \text{if }p \text{ odd}.
\end{cases}
\end{equation}
We put 
$\mathcal{U}^N(z):=\frac{U^N(z)}{i}$, $\mathcal{V}^N(z):=\frac{V^N(z)}{i}$, $\mathcal{W}^N(z):=W^N(z)$ and we analyze what happens with respect to $m,n$.

$\bullet$ \textit{Case} $m,n\in\mathbb N_+$ or $m=0$, $ n\in\mathbb{N}_+.$ From~\eqref{divN} we obtain 
\begin{equation*}
    \mathcal{U}^N_{m,n,p}(z)=\begin{cases}
       b_1\cos(\frac{\pi}{2}pz)\qquad \text{if } p\, \text{even}\vspace{3mm}\\
        a_1\sin(\frac{\pi}{2}pz)\qquad\text{if } p \text{ odd}
    \end{cases}\quad\quad
        \mathcal{V}^N_{m,n,p}(z)=\begin{cases}
         (-\frac{m}{n}b_1+\frac{\pi p}{2n}a_0)\cos(\frac{\pi}{2}pz)\qquad \text{if } p\, \text{even}\vspace{3mm}\\
        -(\frac{m}{n}a_1+\frac{\pi p}{2n}b_0)\sin(\frac{\pi}{2}pz)\hspace{8.5mm} \text{if } p \text{ odd},
    \end{cases}
\end{equation*}
and, in turn~\eqref{final_N} with $a_0=b_0=a_1=b_1=1$.

$\bullet$ \textit{Case} $n=0$, $ m\in\mathbb{N}_+.$ From~\eqref{divN} we find
\begin{equation*}
    \mathcal{U}_{m,0,p}^N(z)=\begin{cases}
\frac{\pi p}{2m}a_0\cos(\frac{\pi}{2}p z)\qquad \text{if } p \text{ even}\\
-\frac{\pi p}{2m}b_0\sin(\frac{\pi}{2}p z)\hspace{6mm} \text{if }p \text{ odd}
\end{cases}\quad\quad \mathcal{V}_{m,0,p}^N(z)=\begin{cases}
b_2\cos(\frac{\pi}{2}p z)\qquad \text{if } p \text{ even}\\
a_2\sin(\frac{\pi}{2}p z)\hspace{7mm} \,\text{if }p \text{ odd}.
\end{cases}
\end{equation*}
and, in turn~\eqref{final_N2} with $a_0=b_0=a_2=b_2=1$.

\subsection{Proof of Corollary~\ref{corollaryD}}\label{proof:corollaryD}
If $\beta \to \infty$, by Proposition~\ref{lemma0} we know $\lambda^D>0$. We retrace the main steps in the proof of Theorem~\ref{eigcube}.
\\ \\
\textbf{\large{$\star$ The case $m=n=0$ ($\mu=0$).}}\\
\noindent
\\
Referring to~\eqref{stokesmnzero}, the only difference is in the boundary conditions on $U^D(z)$ and $V^D(z)$, that becomes purely Dirichlet. Hence,
from~\eqref{stokesmnzero} we have $W^D(z)\equiv 0$ and the pressure $Q^D(z) \equiv Q^D_0$ with $Q_0^D\in\R$. 
The problems for $U^D(z)$ and $V^D(z)$ decouple, getting for instance $U^D(z)=a_1\sin(\sqrt{\lambda} z)+b_1 \cos(\sqrt{\lambda}z)$  for some $a_1,b_1\in\R$ satisfying the Dirichlet boundary conditions, i.e.
\begin{equation*}
\begin{cases}
+a_1\sin(\sqrt{\lambda}) +b_1\cos(\sqrt{\lambda})=0\\
-a_1\sin(\sqrt{\lambda}) +b_1\cos(\sqrt{\lambda})=0.
\end{cases}
\end{equation*}
The system admits nontrivial solutions if and only if
the eigenvalues satisfy 
\begin{equation*}
\sin{(2\sqrt{\lambda^D})}=0\quad \iff\quad \lambda_{0,0,p}^D=\dfrac{\pi^2}{4}p^2\qquad\text{with } p\in\mathbb N_+.
\end{equation*}
The corresponding eigenfunctions $(u_{0,0,p}^D,v_{0,0,p}^D,0)$ are written explictly in~\eqref{eigD3d}.
\\ \\
\textbf{\large{$\star$ The case $m\in\mathbb N_+ \vee n\in\mathbb N_+$ ($\mu>0$).}}\\
\noindent
\\
From~\eqref{stokes23d} the pressure is of the form~\eqref{pressure3d} and $W^D(z)$ solves~\eqref{Wmnp}.
From now on we distinguish three different subcases, from which we will see that only the third case $\lambda^D>\mu^2$ gives nontrivial solutions. In the sequel for brevity we denote the eigenvalue  $\lambda^D$ (resp. $\Lambda^D$) simply by $\lambda$ (resp. $\Lambda$).\\
\noindent
\\
\textbf{\normalsize{1) Case $\lambda= \mu^2$.}} 
From~\eqref{stokes23d}$_1$-\eqref{stokes23d}$_2$ we obtain
\begin{equation*}
\begin{split}
        U^D(z) = a_1z+b_1+\frac{c_1im}{\mu^2}\sinh(\mu z) +\dfrac{c_2im}{\mu^2} \cosh(\mu z)\quad\quad
        V^D(z) = a_2z+b_2+\frac{c_1in}{\mu^2}\sinh(\mu z) +\dfrac{c_2in}{\mu^2} \cosh(\mu z),
\end{split}
\end{equation*}
with $a_1$, $b_1$, $a_2$, $b_2\in\mathbb C$.  From the boundary conditions $U^D(\pm 1)=0$ we find
\begin{equation}\label{sist3dD}
    \begin{cases}
        \displaystyle +a_1+b_1+ c_1\frac{im}{\mu^2}\sinh(\mu)+c_2\frac{im}{\mu^"}\cosh(\mu)=0\\
        \displaystyle -a_1+b_1- c_1\frac{im}{\mu^2}\sinh(\mu)+c_2\frac{im}{\mu^2}\cosh(\mu)=0,
    \end{cases}
\end{equation}
so that if $c_1=c_2=0$ we find only the trivial solution. If $c_1\neq 0 \vee c_2\neq 0$ we have to distinguish further subcases.

$\bullet$ \textit{Case} $m,n\in\mathbb N_+.$
We obtain
\begin{equation*}
    \begin{cases}
    \displaystyle a_1=- c_1\frac{im}{\mu^2}\sinh(\mu)\\
        \displaystyle b_1=-c_2\frac{im}{\mu^2}\cosh(\mu),
    \end{cases}
\end{equation*}
and, applying the divergence condition~\eqref{stokes23d}$_4$, see~\eqref{divlambamu2}, from the boundary conditions $V^D(\pm 1)=0$ we get
\begin{equation*}
    \begin{cases}
    \displaystyle c_1\mu\sinh(\mu)+c_2\big[\mu\cosh(\mu)-\sinh(\mu)\big]=0\\
        \displaystyle -c_1\mu\sinh(\mu)+c_2\big[\mu\cosh(\mu)-\sinh(\mu)\big]=0,
    \end{cases}
\end{equation*}
having nontrivial solutions if and only if
$\mu\sinh(\mu)\big[\mu\cosh(\mu)-\sinh(\mu)\big]=0$;
since $\mu\cosh(\mu)>\sinh(\mu)$ for all $\mu>0$, the previous factors are strictly positive and we find only the trivial solution.

$\bullet$ \textit{Case} $m=0$, $ n\in\mathbb{N}_+.$ The system~\eqref{sist3dD} gives the trivial solution.

$\bullet$ \textit{Case} $n=0$, $ m\in\mathbb{N}_+.$  We obtain the same result as before considering the homologue of \eqref{sist3dD}.\\ \\
\noindent
\textbf{2) Case $0<\lambda < \mu^2$.} 
From~\eqref{stokes23d}$_1$-\eqref{stokes23d}$_2$ we get
\begin{equation*}
\begin{split}
        U^D(z) = a_1\sinh(\sqrt{\mu^2-\lambda}z)+b_1\cosh(\sqrt{\mu^2-\lambda}z)+\frac{c_1im}{\lambda}\sinh(\mu z) +\dfrac{c_2im}{\lambda} \cosh(\mu z)\\
        V^D(z) = a_2\sinh(\sqrt{\mu^2-\lambda}z)+b_2\cosh(\sqrt{\mu^2-\lambda}z)+\frac{c_1in}{\lambda}\sinh(\mu z) +\dfrac{c_2in}{\lambda} \cosh(\mu z),
\end{split}
\end{equation*}
with $a_1$, $b_1$, $a_2$, $b_2\in\mathbb C$. From the boundary conditions $U^D(\pm 1)=0$ we find
\begin{equation}\label{sist3dD.2}
    \begin{cases}
        \displaystyle +a_1\sinh(\sqrt{\mu^2-\lambda})+b_1\cosh(\sqrt{\mu^2-\lambda})+ c_1\frac{im}{\mu}\sinh(\mu)+c_2\frac{im}{\mu}\cosh(\mu)=0\\
        \displaystyle -a_1\sinh(\sqrt{\mu^2-\lambda})+b_1\cosh(\sqrt{\mu^2-\lambda})- c_1\frac{im}{\mu}\sinh(\mu)+c_2\frac{im}{\mu}\cosh(\mu)=0,
    \end{cases}
\end{equation}
so that if $c_1=c_2=0$ the system has nontrivial solutions if and only if $\sinh(\sqrt{\mu^2-\lambda})\cosh(\sqrt{\mu^2-\lambda})=0$, not possible for $\lambda\in(0,\mu^2)$.
If $c_1\neq 0 \vee c_2\neq 0$ we have to distinguish further subcases.

$\bullet$ \textit{Case} $m,n\in\mathbb N_+.$
We obtain
\begin{equation*}
    \begin{cases}
    \displaystyle a_1=- c_1\frac{im}{\lambda}\frac{\sinh(\mu)}{\sinh(\sqrt{\mu^2-\lambda})}\\
        \displaystyle b_1=-c_2\frac{im}{\lambda}\frac{\cosh(\mu)}{\cosh(\sqrt{\mu^2-\lambda})},
    \end{cases}
\end{equation*}
where the denominators are strictly positive. Applying the divergence condition~\eqref{stokes23d}$_4$, see~\eqref{div2}, from the boundary conditions $V^D(\pm 1)=0$ we get
\begin{equation*}
    \begin{cases}
    \displaystyle +c_1\big[\mu\sinh(\mu)-\sqrt{\mu^2-\lambda}\cosh(\mu)
        \tanh{(\sqrt{\mu^2-\lambda})}\big]+c_2\big[\mu\cosh(\mu)-\sqrt{\mu^2-\lambda}\sinh(\mu)\coth(\sqrt{\mu^2-\lambda}) \big]=0\vspace{1.5mm}\\
        \displaystyle -c_1\big[\mu\sinh(\mu)-\sqrt{\mu^2-\lambda}\cosh(\mu)
        \tanh{(\sqrt{\mu^2-\lambda})}\big]+c_2\big[\mu\cosh(\mu)-\sqrt{\mu^2-\lambda}\sinh(\mu)\coth(\sqrt{\mu^2-\lambda}) \big]=0,
    \end{cases}
\end{equation*}
having nontrivial solutions if and only if
\begin{equation*}
\begin{split}
      &\big[\mu\sinh(\mu)-\sqrt{\mu^2-\lambda}\cosh(\mu)
        \tanh{(\sqrt{\mu^2-\lambda})}\big]\big[\mu\cosh(\mu)-\sqrt{\mu^2-\lambda}\sinh(\mu)\coth(\sqrt{\mu^2-\lambda}) \big]=0.
\end{split}
\end{equation*}
This is not possible, since the two factors are strictly positive for $\lambda\in(0,\mu^2)$.

$\bullet$ \textit{Case} $m=0$, $ n\in\mathbb{N}_+.$ The system~\eqref{sist3dD.2} gives the trivial solution.

$\bullet$ \textit{Case} $n=0$, $ m\in\mathbb{N}_+.$  We obtain the same result as before considering the homologue of \eqref{sist3dD.2}.\\ \\
\noindent
\textbf{3) Case $\lambda > \mu^2$.}  The solution to~\eqref{Wmnp} is \eqref{WW} and the boundary conditions~\eqref{stokes23d}$_8$  yield to the system \eqref{sist333}.
We distinguish the next cases as in the statement of the corollary.

$i)$ \underline{$c_1=c_2=0$, i.e. $q^D=\nabla q^D\equiv 0$}\\
\\
The system~\eqref{sist333} gives nontrivial solutions if and only if 
\begin{equation*}
\sin({2\sqrt{\lambda-\mu^2}})=0\quad \iff\quad \lambda_{m,n,p}^D=\mu^2+\dfrac{\pi^2}{4}p^2\qquad\text{with } p\in\mathbb N_+.
\end{equation*}
Then we have
\begin{equation*}
    W^D(z)=\begin{cases}
    a_0\sin(\frac{\pi p}{2}z)\quad \text{if }p \text{ even}\\
    b_0\cos(\frac{\pi p}{2}z)\quad \text{if }p \text{ odd}.
    \end{cases}
\end{equation*}
From~\eqref{stokes23d}$_1$-\eqref{stokes23d}$_2$ and the Dirichlet conditions we get
\begin{equation*}
\begin{split}
        U^D(z) = i\begin{cases}
        a_1\sin(\frac{\pi p}{2}z)\quad \text{if }p \text{ even}\\
        b_1\cos(\frac{\pi p}{2}z)\quad\text{if }p \text{ odd}
        \end{cases}\quad\quad
         V^D(z) = i\begin{cases}
        a_2\sin(\frac{\pi p}{2}z)\quad \text{if }p \text{ even}\\
        b_2\cos(\frac{\pi p}{2}z)\quad\text{if }p \text{ odd},
        \end{cases}
\end{split}
\end{equation*}
with $a_1$, $b_1$, $a_2$, $b_2\in\mathbb R$, where by convenience, instead of considering the usual real solutions combination we consider the imaginary one. Through the divergence condition we have
\begin{equation}\label{divD}
\begin{cases}
na_2=-ma_1+\frac{\pi p}{2}a_0\hspace{6mm} \text{if } p \text{ even}\\
nb_2=-mb_1-\frac{\pi p}{2}b_0\qquad \text{if }p \text{ odd}.
\end{cases}
\end{equation}
We put 
$\mathcal{U}^D(z):=\frac{U^D(z)}{i}$, $\mathcal{V}^D(z):=\frac{V^D(z)}{i}$, $\mathcal{W}^D(z):=W^D(z)$ and we analyze what happens with respect to $m,n$.

$\bullet$ \textit{Case} $m,n\in\mathbb N_+$ or $m=0$, $ n\in\mathbb{N}_+.$ From~\eqref{divD} we obtain 
\begin{equation*}
    \mathcal{U}^D_{m,n,p}(z)=\begin{cases}
       a_1\sin(\frac{\pi}{2}pz)\qquad \text{if } p\, \text{even}\vspace{3mm}\\
        b_1\cos(\frac{\pi}{2}pz)\qquad\text{if } p \text{ odd}
    \end{cases}\qquad
        \mathcal{V}^D_{m,n,p}(z)=\begin{cases}
         (-\frac{m}{n}a_1+\frac{\pi p}{2n}a_0)\sin(\frac{\pi}{2}pz)\qquad \text{if } p\, \text{even}\vspace{3mm}\\
        -(\frac{m}{n}b_1+\frac{\pi p}{2n}b_0)\cos(\frac{\pi}{2}pz)\hspace{8.5mm} \text{if } p \text{ odd},
    \end{cases}
\end{equation*}
and, in turn~\eqref{final_D} with $a_0=b_0=a_1=b_1=1$.

$\bullet$ \textit{Case} $n=0$, $ m\in\mathbb{N}_+.$ From~\eqref{divD} we find
\begin{equation*}
\mathcal{U}_{m,0,p}^D(z)=\begin{cases}
\frac{\pi p}{2m}a_0\sin(\frac{\pi}{2}p z)\qquad \text{if } p \text{ even}\\
-\frac{\pi p}{2m}b_0\cos(\frac{\pi}{2}p z)\hspace{5mm} \text{if }p \text{ odd}
\end{cases}\quad\quad \mathcal{V}_{m,0,p}^D(z)=\begin{cases}
a_2\sin(\frac{\pi}{2}p z)\qquad \text{if } p \text{ even}\\
b_2\cos(\frac{\pi}{2}p z)\hspace{7mm} \,\text{if }p \text{ odd}.
\end{cases}
\end{equation*}
and, in turn~\eqref{final_D2} with $a_0=b_0=a_2=b_2=1$.

$ii)$ \underline{$c_1\neq 0 \vee c_2\neq 0$, i.e. $\nabla q^D\not\equiv 0$}.\\
\\
We prove at first the following lemma.

\begin{lemma}\label{lemmaD}
	Let $\nabla q^D\not\equiv 0$, $\mu>0$ and $p\in\mathbb N_+$, then $\Lambda^D=\mu^2+p^2\dfrac{\pi^2}{4}$ is not an eigenvalue of~\eqref{stokes}.
\end{lemma}
\begin{proof}
The condition $\nabla q^D\not\equiv 0$ implies to study~\eqref{sist333} with $c_1\neq 0 \vee c_2\neq 0$.
If $\Lambda^D=\mu^2+p^2\frac{\pi^2}{4}$ with $p$ even the system~\eqref{sist333} may have nontrivial solutions only if $c_1\neq 0\wedge c_2=0$; note that if $c_1,c_2\neq 0$ and $c_1=0\wedge c_2\neq0$ the system~\eqref{sist333} does not admit solution. 

If $c_1\neq 0\wedge c_2=0$ we find $W^D(z)=a_0\sin(\frac{p\pi}{2} z)-\frac{c_1\mu}{\mu^2+p^2\pi^2/4}\frac{\cosh(\mu)}{(-1)^{p/2}}\cos(\frac{p\pi}{2} z)+\frac{c_1\mu}{\mu^2+p^2\pi^2/4}\cosh (\mu z)$ and $Q^D(z)=c_1\sinh(\mu z)$; computing $U^D(z)$ we do not find $a_1,b_1\in\mathbb{C}$ compatible with the boundary condition $U^D(\pm1)=0$. Similarly if $\Lambda^D=\mu^2+p^2\frac{\pi^2}{4}$ with $p$ odd, where we should discuss only the case $c_1=0\wedge c_2\neq0$. 
\end{proof}
Thanks to Lemma~\ref{lemmaD} we solve~\eqref{sist333}, getting
\begin{equation*}
    \begin{cases}
a_0 = -\frac{c_2\mu}{\Lambda}\frac{\sinh(\mu)}{\sin(\sqrt{\Lambda-\mu^2})}\\
    b_0=-\frac{c_1\mu}{\Lambda}\frac{\cosh(\mu)}{\cos(\sqrt{\Lambda-\mu^2})}.
    \end{cases}
\end{equation*}
From~\eqref{stokes23d}$_1$-\eqref{stokes23d}$_2$ we compute
\begin{equation*}
\begin{split}
        U^D(z) = a_1\sin(\sqrt{\Lambda-\mu^2}z)+b_1\cos(\sqrt{\Lambda-\mu^2}z)+ \frac{c_1im}{\Lambda}\sinh(\mu z) +\dfrac{c_2im}{\Lambda} \cosh(\mu z)\\
        V^D(z) = a_2\sin(\sqrt{\Lambda-\mu^2}z)+b_2\cos(\sqrt{\Lambda-\mu^2}z)+\frac{c_1in}{\Lambda}\sinh(\mu z) +\dfrac{c_2in}{\Lambda} \cosh(\mu z),
\end{split}
\end{equation*}
with $a_1$, $b_1$, $a_2$, $b_2\in\mathbb C$. Applying the divergence condition~\eqref{stokes23d}$_4$ we obtain~\eqref{div}.
Hence, from the Dirichlet conditions $U^D(\pm 1)=0$, we get
\begin{equation*}
    \begin{cases}
        \displaystyle +a_1\sin(\sqrt{\Lambda-\mu^2})+b_1\cos(\sqrt{\Lambda-\mu^2})+ c_1\frac{im}{\Lambda}\sinh(\mu)+c_2\frac{im}{\Lambda}\cosh(\mu)=0\\
        \displaystyle -a_1\sin(\sqrt{\Lambda-\mu^2})+b_1\cos(\sqrt{\Lambda-\mu^2})-c_1\frac{im}{\Lambda}\sinh(\mu)+c_2\frac{im}{\Lambda}\cosh(\mu)=0,
    \end{cases}
\end{equation*}
giving
$$
\begin{cases}
        \displaystyle a_1=-c_1\frac{im}{\Lambda}\frac{\sinh(\mu)}{\sin(\sqrt{\Lambda-\mu^2})}\\
        \displaystyle b_1=-c_2\frac{im}{\Lambda}\frac{\cosh(\mu)}{\cos(\sqrt{\Lambda-\mu^2})}.
    \end{cases}
$$
We distinguish the following subcases. 

$\bullet$ \textit{Case} $m,n \in\mathbb N_+.$
Imposing $V^D(\pm1)=0$ we obtain
\begin{equation*}
    \begin{cases}
        +c_1[\cosh(\mu) \sqrt{\Lambda-\mu^2}\tan (\sqrt{\Lambda-\mu^2})+ \mu \sinh(\mu)]+ c_2[-\sinh(\mu)  \sqrt{\Lambda-\mu^2}\cot (\sqrt{\Lambda-\mu^2})+\mu \cosh(\mu)]=0
        \vspace{.2cm}
        \\
        -c_1[\cosh(\mu) \sqrt{\Lambda-\mu^2}\tan (\sqrt{\Lambda-\mu^2})+ \mu \sinh(\mu)]+ c_2[-\sinh(\mu)  \sqrt{\Lambda-\mu^2}\cot (\sqrt{\Lambda-\mu^2})+\mu \cosh(\mu)]=0,
    \end{cases}
\end{equation*}
having nontrivial solutions if and only if
\begin{equation*}    
[\cosh(\mu) \sqrt{\Lambda-\mu^2}\tan( \sqrt{\Lambda-\mu^2})+ \mu \sinh(\mu)][-\sinh(\mu)  \sqrt{\Lambda-\mu^2}\cot( \sqrt{\Lambda-\mu^2})+\mu \cosh(\mu)]=0.
\end{equation*}
This is equivalent to search solutions of one of the two equations in~\eqref{eig3dbetainfty} for some $\Lambda>\mu^2$.  Hence,
we obtain two increasing and positively divergent sequences of eigenvalues so that, combined together, the spectrum is given by
\begin{equation*}
  \La^D_{m,n,p}\in \bigg(\mu^2+\frac{\pi^2}{4}(1+p)^2\,,\,\mu^2+\frac{\pi^2}{4}(2+p)^2\bigg) \qquad \forall p\in \mathbb N.
\end{equation*}

We compute the corresponding eigenfunctions, observing that $U^D(z)$ and $V^D(z)$ are imaginary pure; then, we put 
$\mathcal{U}^D(z):=\frac{U^D(z)}{i}$, $\mathcal{V}^D(z):=\frac{V^D(z)}{i}$ and $\mathcal{W}^D(z):=W^D(z)$ where

\small{
\begin{equation*}
\begin{split}
    &\mathcal U^D(z)=\frac{m}{\Lambda}\begin{cases}
         c_1\dfrac{\sin(\sqrt{\La-\mu^2})\sinh(\mu z)-\sinh(\mu)\sin(\sqrt{\La-\mu^2}z)}{\sin(\sqrt{\La-\mu^2})}\hspace{5mm} \text{if } p=0\,\, \text{or even}\\
 c_2\dfrac{\cos(\sqrt{\La-\mu^2})\cosh(\mu z)-\cosh(\mu)\cos(\sqrt{\La-\mu^2}z)}{\cos(\sqrt{\La-\mu^2})}\hspace{5mm} \text{if } p \text{ odd}
    \end{cases}
    \vspace{.2cm}\\
      &\mathcal V^D(z)=\begin{cases}
       \dfrac{c_1}{\La}\bigg[ n
     \sinh(\mu z)+\dfrac{[m^2\sinh(\mu)+\mu\sqrt{\La-\mu^2}\cosh(\mu)\tan(\sqrt{\La-\mu^2})]\sin(\sqrt{\La-\mu^2}z)}{n\sin(\sqrt{\La-\mu^2})}\bigg]\hspace{5mm} \text{if } p=0\,\,\, \text{or }\, \text{even}\\
 \dfrac{c_2}{\La}\bigg[n\cosh(\mu z)\!+\!\dfrac{[m^2\cosh(\mu)-\mu\sqrt{\La-\mu^2}\sinh(\mu)\cot(\sqrt{\La-\mu^2})]\cos(\sqrt{\La-\mu^2}z)}{n\cos(\sqrt{\La-\mu^2})}\bigg]\hspace{6mm} \text{if } p \text{ odd}
    \end{cases}
    \vspace{.2cm}\\
    &\mathcal W^D(z)=\mu\begin{cases}\dfrac{c_1}{\La}
        \dfrac{\cos(\sqrt{\La-\mu^2})\cosh(\mu z)-\cosh(\mu)\cos(\sqrt{\La-\mu^2}z)}{\cos(\sqrt{\La-\mu^2})}\hspace{3mm} \text{if } p=0\,\, \text{or even}\\
        \dfrac{c_2}{\La}\dfrac{\sin(\sqrt{\La-\mu^2})\sinh(\mu z)-\sinh(\mu)\sin(\sqrt{\La-\mu^2}z)}{\sin(\sqrt{\La-\mu^2})}\hspace{5mm} \text{if } p \text{ odd},
    \end{cases}
    \end{split}
\end{equation*}}
\normalsize
associating the pressure
\begin{equation*}
   Q^D(z)=\begin{cases}
        c_1\sinh(\mu z)\qquad \text{if } p=0\,\, \text{or even}\\
        c_2\cosh(\mu z)\qquad \text{if } p \text{ odd}.
    \end{cases}
\end{equation*}
Chosing  e.g. $c_1=c_2=\La$, we find~\eqref{eigvD} with the pressure in~\eqref{pressD}.

We observe that the cases $c_1=0$ and $c_2\neq 0$ or $c_1\neq0$ and $c_2=0$ lead to discuss only one among~\eqref{eig3dbetainfty}, yielding to the same conclusion.

$\bullet$ \textit{Case} $m=0$, $ n\in\mathbb{N}_+.$ 
We observe the decoupling of the $U^D(z)$ component that solves
$$
\begin{cases}
    (U^D)''(z)+(\Lambda-n^2)U^D(z)=0\quad z\in(-1,1)\\
    U^D(\pm1)=0.
\end{cases}
$$
Therefore two scenarios are possible: $U^D(z)=a_1\sin(\sqrt{\Lambda-n^2}z)+b_1\cos(\sqrt{\Lambda-n^2}z)$ or $U^D(z)\equiv 0$. Due to Lemma~\ref{lemmaD} in the first case $a_1=b_1=0$.
Therefore, it remains the case $U^D(z)\equiv 0$ where the eigenfunctions are $(\mathcal U^D_{0,n,p}(z)\mathcal{P}^u_{0,n}(y)\,,\,0\,,\,\mathcal W^D_{0,n,p}(z)\mathcal{P}_{0,n}(y))$; 
here the $z-$component can be inferred by~\eqref{eigvD} with $m=0$ and $\mu=n$. The corresponding pressure $q^D_{0,n,p}$ comes from~\eqref{pressD} with $m=0$ and $\mu=n$.

$\bullet$ \textit{Case} $n=0$, $ m\in\mathbb{N}_+.$ We observe the decoupling of the $V^D(z)$ component that solves
$$
\begin{cases}
    (V^D)''(z)+(\Lambda-m^2)V^D(z)=0\quad z\in(-1,1)\\
    V^D(\pm 1)=0.
\end{cases}
$$
Therefore two scenarios are possible: $V^D(z)=a_2\sin(\sqrt{\Lambda-m^2}z)+b_2\cos(\sqrt{\Lambda-m^2}z)$ or $V^D(z)\equiv 0$. Due to Lemma~\ref{lemmaD} in the first case $a_2=b_2=0$.
Therefore, it remains the case $V^D(z)\equiv 0$ where the eigenfunctions are $(\mathcal U^D_{m,0,p}(z)\mathcal{P}^u_{m,0}(x)\,,\,0\,,\,\mathcal W^D_{m,0,p}(z)\mathcal{P}_{m,0}(x))$; 
here the $z-$component can be inferred by~\eqref{eigvD} with $n=0$ and $\mu=m$, see~\eqref{eigvappendix_m0pD}. The corresponding pressure $q^D_{m,0,p}$ comes from~\eqref{pressD} with $n=0$ and $\mu=m$, see~\eqref{press.m0p.appendixD}.
\subsection{Proof of
  Corollary~\ref{cor:monotonicitybeta}}\label{proof:corollarymono}
We know that for any fixed $m,n,p \in \mathbb N$, the eigenvalues $\lambda_{m,n,p}(\beta)$ and $\Lambda_{m,n,p}(\beta)$ satisfy respectively one among equations~\eqref{eig3d2}$_1$-\eqref{eig3d2}$_2$ and~\eqref{eig3d}$_1$-\eqref{eig3d}$_2$, see Theorem~\ref{eigcube}.

We focus at first on $\Lambda$ eigenvalues satisfying~\eqref{eig3d}$_1$. 
Solutions are given by the intersections between the function $f(\Lambda)=\sqrt{\Lambda-\mu^2}\tan(\sqrt{\Lambda-\mu^2})$ and the negative decreasing half-line $r^{-}(\Lambda,\beta)=-\Lambda/\beta-\mu\tanh(\mu)$ with starting point in $(\mu^2,-\mu[\mu/\beta+\tanh(\mu)])$. If $0<\beta_1<\beta_2<\infty$, then 
     \begin{equation*}
       r^-(\Lambda,\beta_1)< r^-(\Lambda,\beta_2).
     \end{equation*}
     If $\La(\beta_i)$ ($i=1,2$) is the unique solution of equation~\eqref{eig3d}$_1$  with the two different $\beta_i$, then we have
     \begin{equation*}
         f(\La(\beta_1))< f(\La(\beta_2)),
     \end{equation*}
     and since $f(\Lambda)$ is strictly increasing, we have $\La(\beta_1)<\La(\beta_2)$.
     
    Solutions to equation~\eqref{eig3d}$_2$ are intersections between the function $g(\Lambda)=\sqrt{\Lambda-\mu^2}\cot({\sqrt{\Lambda-\mu^2}})$ and the increasing positive half-line $r^+(\Lambda,\beta)= \Lambda/\beta+\mu\coth(\mu)$. In this case for any fixed $\Lambda$, $r^+(\Lambda,\beta)$ is decreasing in $\beta$, so that for any $0<\beta_1<\beta_2<\infty$
     \begin{equation*}
         g(\La(\beta_1))> g(\La(\beta_2)).
     \end{equation*}
     Since $g(\Lambda)$ is strictly decreasing we have again $\La(\beta_1)<\La(\beta_2)$.
     
     We observe that for $\beta \to \infty$, equations~\eqref{eig3d}$_1$ and~\eqref{eig3d}$_2$, become respectively~\eqref{eig3dbetainfty}$_1$ and~\eqref{eig3dbetainfty}$_2$, so that
     \begin{equation}\label{limbetainfty}
         \lim_{\beta \to \infty} \La(\beta) = \La^{D}.
     \end{equation}
     The same arguments follow for $\lambda(\beta)$, solving equations~\eqref{eig3d2}$_1$ and~\eqref{eig3d2}$_2$. Indeed, in these cases, solutions are given by the intersection between the functions $ \tan(\sqrt{\lambda-\mu^2})$ and $\cot(\sqrt{\lambda-\mu^2})$, respectively with  $s^-(\lambda,\beta)=-\sqrt{\lambda-\mu^2}/\beta$ and $s^+(\lambda,\beta)=\sqrt{\lambda-\mu^2}/\beta$. When $\beta\to \infty$, the eigenvalues are $\lambda^D=\mu^2+\frac{\pi^2}{4}p^2$ with $p\in\mathbb{N}_+$ and consequently~\eqref{limbetainfty} hold also for $\lambda$.
\subsection{Proof of
  Proposition~\ref{prop-mono}}\label{proof:prop-mono}
We observe that the eigenfunctions $\uu_{0,n_k,p_k}=(u^k,v^k,w^k)$
have components such that $\frac{\partial u^k}{\partial{x}}=v^k=w^k=0$ for all $k\in \mathbb{N}_+$, implying that   $(\uu_{0,n_k,p_k}\cdot\nabla)\,\uu_{0,n_j,p_j}=0$ for all $k,j\in \mathbb{N}_+$. Hence, $(\vv,q)$ in~\eqref{sol-mono} solves 
	$$
	\vv_t-\Delta \vv=0=-\nabla q\, ,
	$$
	since $\uu_{0,n_k,p_k}$ are eigenfunctions of~\eqref{stokes} corresponding to the eigenvalues $\lambda_{0,n_k,p_k}$.
Therefore, $(\vv,q)$ in~\eqref{sol-mono} solve~\eqref{ns}.
\subsection{Proof of Proposition~\ref{prop-single}}\label{proof:prop-single}
The eigenfunction $\uu_{m,n,p}(x,y,z)$ satisfies $$(\uu_{m,n,p}\cdot\nabla)\,\uu_{m,n,p}=-\dfrac{m}{2}\mathcal{Z}(z)^2\left(\begin{matrix}
    2(ab+cd)\cos(2mx)+(a^2-b^2-c^2+d^2)\sin (2mx)\\
     \dfrac{m}{n}\big[2(bc+ad)\cos(2ny)+(a^2+b^2-c^2-d^2)\sin (2ny)\big]\\
    0
\end{matrix}\right).$$
Therefore, assuming~\eqref{hp} we obtain  that $(\uu_{m,n,p}\cdot\nabla)\uu_{m,n,p}=0$.
Hence, $(\vv,q)$ in~\eqref{sol-single} solves 
	$$
	\vv_t-\Delta \vv=0=-\nabla q\, ,
	$$
	since $\uu_{m,n,p}$ is an eigenfunction of~\eqref{stokes} corresponding to the eigenvalue $\lambda_{m,n,p}$.
Therefore, $(\vv,q)$ in~\eqref{sol-single} solve~\eqref{ns}.
\subsection{Proof of Theorem~\ref{teo_p_nulla}}\label{proof:theorem_p_nulla}
The proof relies on the next crucial lemma, saying that the projection of the nonlinearity onto the space $H_\P$ is zero. 
The projection on $H_\P$ requires to apply the so-called Helmholtz-Weyl decomposition to the nonlinear term~\cite{helmholtz,weyl}.
This operation is explicit but quite lengthy in this context, hence we recall it in the Appendix~\ref{nonlin_proj} and we use it directly in the proof of the following lemma.
\begin{lemma}\label{lemma3}
    Let $m_1,m_2, m_3, p_1,p_2, p_3\in\mathbb{N}$, $n_1,n_2,n_3
    \in\mathbb{N}_+$, and let $\uu_{m_1,n_1,p_1}(x,y,z)$,
    $\uu_{m_2,n_2,p_2}(x,y,z)$, and $\uu_{m_3,n_3,p_3}(x,y,z)$ be three eigenfunctions in Corollary~\ref{corollary:press_cost}. If the coefficients in 
     $\mathcal{P}^u_{m_k,n_k}$ and $\mathcal{P}^v_{m_k,n_k}$, see~\eqref{periodic}, satisfy~\eqref{hp3} for all $k$, then
    \begin{equation}\label{ts}
        	\int_\Omega\big(\uu_{m_1,n_1,p_1}\cdot\nabla\big)\,\uu_{m_2,n_2,p_2}\cdot \uu_{m_3,n_3,p_3}=0.
    \end{equation}
\end{lemma}
\begin{proof}
    We prove the first case in~\eqref{hp3}, where $\mathcal{P}^u_{m_k,n_k}$ and $\mathcal{P}^v_{m_k,n_k}$ have $a,b\in\mathbb{R}$ and $c=d=0$ for all $k$; the other cases are similar.
    Applying the Helmholtz-Weyl decomposition to $(\uu_{m_1,n_1,p_1}\cdot\nabla)\,\uu_{m_2,n_2,p_2}$, see Appendix~\ref{nonlin_proj}, by Computation~\ref{computation3} we find 
     \begin{equation*}
    \begin{split}
    &(\uu_{m_1,n_1,p_1}\cdot\nabla)\,\uu_{m_2,n_2,p_2}=G+\\&\left(\begin{matrix}
            abZ^u_{1}(z)\cos[(m_1+m_2)x]\cos[(n_1-n_2)y]+(a^2+b^2)Z^u_{3}(z)\sin[(m_1-m_2)x]\cos[(n_1+n_2)y]\\+(a^2-b^2)Z^u_{4}(z)\sin[(m_1+m_2)x]\cos[(n_1-n_2)y]+abZ^u_{5}(z)\cos[(m_1+m_2)x]\cos[(n_1+n_2)y]\\
          +(a^2+b^2)Z^u_{9}(z)\sin[(m_1-m_2)x]\cos[(n_1-n_2)y]+(a^2-b^2)Z^u_{10}(z)\sin[(m_1+m_2)x]\cos[(n_1+n_2)y]\\
            \\
            abZ^v_1(z)\sin[(m_1+m_2)x]\sin[(n_1-n_2)y]+(a^2+b^2)Z^v_3(z)\cos[(m_1-m_2)x]\sin[(n_1+n_2)y]\\(a^2-b^2)Z^v_4(z)\cos[(m_1+m_2)x]\sin[(n_1-n_2)y]         +abZ^v_5(z)\sin[(m_1+m_2)x]\sin[(n_1+n_2)y]\\+
          (a^2+b^2)Z^v_9(z)\cos[(m_1-m_2)x]\sin[(n_1-n_2)y]+(a^2-b^2)Z^v_{10}(z)\cos[(m_1+m_2)x]\sin[(n_1+n_2)y]\\
            \\abZ^w_{1}(z)\sin[(m_1+m_2)x]\cos[(n_1-n_2)y]+(a^2+b^2)Z^w_{3}(z)\cos[(m_1-m_2)x]\cos[(n_1+n_2)y]\\+(a^2-b^2)Z^w_{4}(z)\cos[(m_1+m_2)x]\cos[(n_1-n_2)y]+abZ^w_{5}(z)\sin[(m_1+m_2)x]\cos[(n_1+n_2)y]\\
          +(a^2+b^2)Z^w_{9}(z)\cos[(m_1-m_2)x]\cos[(n_1-n_2)y]+(a^2-b^2)Z^w_{10}(z)\cos[(m_1+m_2)x]\cos[(n_1+n_2)y]
        \end{matrix}\right);
     \end{split}
    \end{equation*}
the precise expression of the $z-$functions is not necessary for this proof.
   
    Assuming~\eqref{hp3}$_1$ we find that the vector $\big(\P^u_{m_3,n_3}(x,y),\P^v_{m_3,n_3}(x,y),0\big)$ is $L^2(\To^2)-$orthogonal to the projection on $H_\P$ of $(\uu_{m_1,n_1,p_1}\cdot\nabla)\,\uu_{m_2,n_2,p_2}$ for any choice of $m_3\in\mathbb{N}$, $n_3\in\mathbb{N}_+$, implying that 
    $$
\uu_{m_3,n_3,p_3}=\mathcal{Z}_{m_3,n_3,p_3}(z)\big(\P^u_{m_3,n_3}(x,y),-\tfrac{m_3}{n_3}\P^v_{m_3,n_3}(x,y),0\big),
$$ is $L^2(\Omega)-$orthogonal to the projection onto $H_\P$ of $(\uu_{m_1,n_1,p_1}\cdot\nabla)\,\uu_{m_2,n_2,p_2}$ for any choice of $m_3\in\mathbb{N}$, $n_3\in\mathbb{N}_+$, giving~\eqref{ts}.

 \end{proof}

Let us conclude the proof of Theorem~\ref{teo_p_nulla} as follows.
For $\vv_0\in V_\P$ stated in the assumptions, let $T^*\le T$ be as in Proposition~\ref{existence} and let $\vv=\vv(t)$ be the unique local solution of~\eqref{ns} in $(0,T^*)$, which then satisfies~\eqref{regularitya1}. We write it in the form
$$
\vv(t,x,y,z)=\sum_{k=1}^\infty A_k(t)\, \uu_{m_k,n_k,p_k}(x,y,z)\, \quad\quad \vv(0,x,y,z)=\vv_0(x,y,z)=\sum_{k=1}^\infty\overbrace{ A_k(0)}^{=\gamma_k}\, \uu_{m_k,n_k,p_k}(x,y,z)\in V_\P\, .
$$
If~\eqref{hp3} holds then by Lemma~\ref{lemma3} we have
$$
\int_\Omega\big(\uu_{m_i,n_i,p_i}\cdot\nabla\big)\,\uu_{m_j,n_j,p_j}\cdot\Delta \uu_{m_k,n_k,p_k}=-\lambda_k\int_\Omega\big(\uu_{m_i,n_i,p_i}\cdot\nabla\big)\,\uu_{m_j,n_j,p_j}\cdot \uu_{m_k,n_k,p_k}=0\qquad\forall\,i,j,k\in\mathbb{N}_+,
$$
being $\lambda_k$ the eigenvalue corresponding to $\uu_{m_k,n_k,p_k}$.

Therefore,
\begin{equation*}
	\int_\Omega\big(\vv(t)\cdot\nabla\big)\,\vv(t)\cdot\Delta \vv(t)=-\sum_{k=1}^\infty\lambda_kA_k(t)\int_\Omega\big(\vv(t)\cdot\nabla\big)\,\vv(t)\cdot\uu_{m_k,n_k,p_k}=0\, \quad\mbox{for a.e. }t\in(0,T^*).
\end{equation*}

Thanks to~\eqref{regularitya1} and the fact that the gradient of the pressure is zero in this case, we test~\eqref{ns} by $-\Delta \vv$ getting
\begin{equation*}
\begin{split}
\dfrac{1}{2}\dfrac{d}{dt}\|\nabla \vv(t)\|^2_{L^2(\Omega)}+ \|\Delta \vv(t)\|^2_{L^2(\Omega)} &= \int_\Omega(\vv(t)\cdot\nabla)\,\vv(t)\cdot\Delta \vv(t)=0\,\quad\mbox{for a.e. }t\in(0,T^*)\, .
\end{split}
\end{equation*}
After integration over $(0,t)$, we obtain
\begin{equation}\label{en1}
\|\nabla \vv(t)\|^2_{L^2(\Omega)}\le\|\nabla \vv_0\|^2_{L^2(\Omega)}<\infty\quad\mbox{for a.e. }t\in(0,T^*)\, .
\end{equation}
Since the upper bound
\eqref{en1} is uniform and independent on $T^*$, if it was $T^*<T$ we could extend the solution beyond $T^*$; hence, $T^*=T$. Being $T>0$ arbitrary, we get a global and unique solution with the regularity stated.

\bigskip

\noindent
{\bf Acknowledgements.}
All authors are members of INdAM-GNAMPA. 

 LCB and RS are supported by MIUR
within project PRIN20204NT8W ``Nonlinear evolution PDEs, fluid
dynamics and transport equations: theoretical foundations and
applications'' and MIUR Excellence, Department of Mathematics,
University of Pisa, CUP I57G22000700001, 2023-2027.

LCB is partially supported by National Research Center in High
Performance Computing, Big Data and Quantum Computing (CN1 -- Spoke
6),  

AF is partially supported by the INdAM - GNAMPA project 2023 ``Modelli
matematici di EDP per fluidi e strutture e proprietà geometriche delle
soluzioni di EDP'',  INdAM - GNAMPA  project 2024 “Problemi
frazionari:  propriet\`a quantitative ottimali, simmetria,
regolarit\`a, and he is supported by the MUR (Italy) grant
Dipartimento di Eccellenza 2023-2027, Dipartimento di Matematica,
Politecnico di Milano.

\bigskip

\noindent
{\bf Data availability statement.} Data sharing not applicable to this article as no datasets were generated or analysed during the current study. There are no conflicts of interest.

\appendix 

\section{Details on the spectrum of the Stokes problem} \label{eigenfunctions1}
In this appendix we provide further details on the spectrum of~\eqref{stokes}. Let us begin with the extended version of Theorem~\ref{eigcube}, coming from its proof in Section~\ref{prooftheorem3d}.
\begin{theorem}\label{eigcube1}
Let $\beta>0$ and $\mathcal{P}^u_{m,n}$, $\mathcal{P}^v_{m,n}$, $\mathcal{P}_{m,n}$,  as in~\eqref{periodic} for some $a,b,c,d\in\R$, not contemporary all zero. For any $m,n\in \mathbb{N}$ there exist two positive and increasing sequences of diverging eigenvalues of~\eqref{stokes} $\lambda:=\lambda_{m,n,p}(\beta)$  and $\Lambda:=\Lambda_{m,n,p}(\beta)$ ($p\in\mathbb{N}$), counted with their multiplicity, such that:
\begin{itemize}
    \item[i)] If $\nabla q\equiv 0$ the eigenvalues $\lambda\in \big(\mu^2+\frac{\pi^2}{4}p^2,\mu^2+\frac{\pi^2}{4}(1+p)^2\big)$ satisfy \eqref{eig3d2}.

$\bullet$ If $m=n=0$ and $p\in\mathbb N$ the pressure is $q_{0,0,p}=Q_0\in \R$ and the eigenfunctions are $(u_{0,0,p},v_{0,0,p},0)$ 
\begin{equation}\label{eig03d}
u_{0,0,p}(z)=\begin{cases}
d\cos(\sqrt{\lambda} z)\qquad \text{if }p=0\,\, \text{or even}\\
d\sin(\sqrt{\lambda} z)\qquad \text{if }p \text{ odd}
\end{cases} \,
v_{0,0,p}(z)=\begin{cases}
b\cos(\sqrt{\lambda} z)\qquad \text{if }p=0\,\, \text{or even}\\
b\sin(\sqrt{\lambda} z)\qquad \text{if }p \text{ odd}
\end{cases}
\end{equation}
for some $b,d\in\R\setminus\{0\}$.

$\bullet$ If $n=0$, $m\in\mathbb{N}_+$ and $p\in\mathbb N$ the pressure is $q_{m,0,p}\equiv 0$ and the eigenfunctions are $(0,v_{m,0,p},0)$ where
\begin{equation}\label{eigm0pappendix}
v_{m,0,p}(x,z)=\mathcal{P}^v_{m,0}(x)\begin{cases}
\cos(\sqrt{\lambda-m^2} z)\qquad \text{if }p=0\,\, \text{or even}\\
\sin(\sqrt{\lambda-m^2} z)\qquad \text{if }p \text{ odd}.
\end{cases} 
\end{equation}

 $\bullet$ If $m\in\mathbb{N}$, $n\in\mathbb N_+$ and $p\in\mathbb N$ the pressure is $q_{m,n,p}\equiv 0$ and the eigenfunctions are\\ $(\mathcal U_{m,n,p}(z)\mathcal{P}^u_{m,n}(x,y)\,,\,\mathcal V_{m,n,p}(z)\mathcal{P}^v_{m,n}(x,y)\,,\,0)$ where
\begin{equation}\label{eig_pnullaappendix}
   \mathcal U_{m,n,p}(z)\!=\!\begin{cases}
        \cos(\sqrt{\lambda-\mu^2}z)\hspace{3mm} \text{if } p=0\,\, \text{or ev.}\vspace{3mm}\\
        \sin(\sqrt{\lambda-\mu^2}z)\hspace{3mm} \text{if } p \text{ odd}
    \end{cases}
        \mathcal V_{m,n,p}(z)\!=\!-\dfrac{m}{n}\!\begin{cases}
         \cos(\sqrt{\lambda-\mu^2}z)\hspace{3mm} \text{if } p=0\,\, \text{or ev.}\vspace{3mm}\\
        \sin(\sqrt{\lambda-\mu^2}z)\hspace{3mm} \text{if } p \text{ odd}.
    \end{cases}
\end{equation}
   \item[ii)] if $\nabla q\not\equiv 0$ the eigenvalues $\La\in \big(\mu^2+\frac{\pi^2}{4}(1+p)^2,\mu^2+\frac{\pi^2}{4}(2+p)^2\big)$ with $\mu>0$ satisfy \eqref{eig3d}.

$\bullet$ If $n=0$, $m\in\mathbb{N}_+$ and $p\in\mathbb N$ the pressure is 
\begin{equation}\label{press.m0p.appendix}
   q_{m,0,p}(x,z)=\mathcal{P}_{m,0}(x)\begin{cases}
        \La\sinh(m z)\qquad \text{if } p=0\,\, \text{or even}\\
        \La\cosh(m z)\qquad \text{if } p \text{ odd}
    \end{cases}
\end{equation}
 and the eigenfunctions are $(\mathcal U_{m,0,p}(z)\mathcal{P}^v_{m,0}(x)\,,\,0\,,\,\mathcal{W}_{m,0,p}(z)\mathcal{P}_{m,0}(x))$,  where
\begin{equation}\label{eigvappendix_m0p}
\begin{split}
    &\mathcal U_{m,0,p}(z)=m\begin{cases}
        \sinh(m z)-\dfrac{[m\cosh(m)+\beta\sinh(m)]\sin(\sqrt{\La- m^2}z)}{\beta\sin(\sqrt{\La- m^2})+\sqrt{\La-m^2}\cos(\sqrt{\La-m^2})}\hspace{5mm} \text{if } p=0\,\, \text{or even}\\
\cosh(m z)-\dfrac{[m\sinh(m)+\beta \cosh(m)]\cos(\sqrt{\La-m^2}z)}{\beta\cos(\sqrt{\La-m^2})-\sqrt{\La-m^2}\sin(\sqrt{\La-m^2})}\hspace{5mm} \text{if } p \text{ odd},
    \end{cases}\\
    &\mathcal W_{m,0,p}(z)=m\begin{cases}
        \dfrac{\cos(\sqrt{\La-m^2})\cosh(m z)-\cosh(m) \cos(\sqrt{\La-m^2}z)}{\cos(\sqrt{\La-m^2})}\hspace{3mm} \text{if }p=0\,\, \text{or even}\\
        \dfrac{\sin\sqrt{\La-m^2}\sinh(m z)-\sinh(m)\sin(\sqrt{\La-m^2}z)}{\sin(\sqrt{\La-m^2})}\hspace{7mm} \text{if } p \text{ odd}.
    \end{cases}
    \end{split}
\end{equation}

 $\bullet$ If $m\in\mathbb{N}$, $n\in\mathbb N_+$ and $p\in\mathbb N$ the pressure is 
\begin{equation}\label{press2appendix}
   q_{m,n,p}(x,y,z)=\mathcal{P}_{m,n}(x,y)\begin{cases}
        \La\sinh(\mu z)\qquad \text{if } p=0\,\, \text{or even}\\
        \La\cosh(\mu z)\qquad \text{if } p \text{ odd}
    \end{cases}
\end{equation}
and the eigenfunctions are 
$\big(\mathcal U_{m,n,p}(z)\mathcal{P}^u_{m,n}(x,y), \mathcal V_{m,n,p}(z)\mathcal{P}^v_{m,n}(x,y),\mathcal W_{m,n,p}(z)\mathcal{P}_{m,n}(x,y)\big)$ where
\small{
\begin{equation}\label{eigvappendix}
\begin{split}
    &\mathcal U_{m,n,p}(z)=m\cdot\\&\begin{cases}
        \dfrac{[\beta\sin(\sqrt{\La-\mu^2})+\sqrt{\La-\mu^2}\cos(\sqrt{\La-\mu^2})]\sinh(\mu z)-[\mu\cosh(\mu)+\beta\sinh(\mu)]\sin(\sqrt{\La-\mu^2}z)}{\beta\sin(\sqrt{\La-\mu^2})+\sqrt{\La-\mu^2}\cos(\sqrt{\La-\mu^2})}\hspace{2mm} \text{if } p=0\,\, \text{or ev.}\\
\dfrac{[\beta\cos(\sqrt{\La-\mu^2})-\sqrt{\La-\mu^2}\sin(\sqrt{\La-\mu^2})]\cosh(\mu z)-[\mu\sinh(\mu)+\beta\cosh(\mu)]\cos(\sqrt{\La-\mu^2}z)}{\beta\cos(\sqrt{\La-\mu^2})-\sqrt{\La-\mu^2}\sin(\sqrt{\La-\mu^2})}\hspace{2mm} \text{if } p \text{ odd},
    \end{cases}\vspace{.2cm}
    \\
      &\mathcal V_{m,n,p}(z)=\\&\begin{cases}
       n
     \sinh(\mu z)\!+\\\dfrac{[\mu(\La-n^2)\cosh(\mu)+\beta(m^2\sinh(\mu)+\mu\sqrt{\La-\mu^2}\cosh(\mu)\tan(\sqrt{\La-\mu^2}))]\sin(\sqrt{\La-\mu^2}z)}{n[\beta\sin(\sqrt{\La-\mu^2})+\sqrt{\La-\mu^2}\cos(\sqrt{\La-\mu^2})]}\hspace{1mm} \text{if } p=0\,\text{or ev.}\\
n\cosh(\mu z)+\\\dfrac{[\mu(\La-n^2)\sinh(\mu)+\beta(m^2\cosh(\mu)-\mu\sqrt{\La-\mu^2}\sinh(\mu)\cot(\sqrt{\La-\mu^2}))]\cos(\sqrt{\La-\mu^2}z)}{n[\beta\cos(\sqrt{\La-\mu^2})-\sqrt{\La-\mu^2}\sin(\sqrt{\La-\mu^2})]}\hspace{1mm} \text{if } p \text{ odd},
    \end{cases}
    \vspace{.2cm}\\
    &\mathcal W_{m,n,p}(z)=\mu\begin{cases}
        \dfrac{\cos(\sqrt{\La-\mu^2})\cosh(\mu z)-\cosh(\mu)\cos(\sqrt{\La-\mu^2}z)}{\cos(\sqrt{\La-\mu^2})}\hspace{3mm} \text{if }p=0\,\, \text{or even}\\
        \dfrac{\sin(\sqrt{\La-\mu^2})\sinh(\mu z)-\sinh(\mu)\sin(\sqrt{\La-\mu^2}z)}{\sin(\sqrt{\La-\mu^2})}\hspace{5mm} \text{if } p \text{ odd}.
    \end{cases}
    \end{split}
\end{equation}}
\normalsize
\end{itemize}  
 \end{theorem}
In the next proposition we give some further details on the first eigenvalue of~\eqref{stokes}.
\begin{proposition}\label{primo}
    Let $\beta>0$. The first eigenvalue of~\eqref{stokes} is $\lambda_{0,0,0}\in (0,\frac{\pi^2}{4}) $.
\end{proposition}
\begin{proof}
Fixed $m,n,p\in\mathbb N$, if $m$ and $n$ are not contemporary zero, then $\lambda_{m,n,p}< \La_{m,n,p}$, since their interval of definition is translated by $1$, see Theorem~\ref{eigcube}. Thus we only need to prove that 
 $\lambda_{0,0,0}<\min\{\lambda_{0,1,0},\lambda_{1,1,0}\}$. Let us begin proving that $\lambda_{0,0,0}<\lambda_{0,1,0}=\lambda_{1,0,0}$. We know that 
 $\lambda_{0,0,0}$ and $\lambda_{0,1,0}$ are respectively zeros of the functions
 \begin{equation*}
     \begin{split}
         &f(s) = \cot(\sqrt{s})-\frac{\sqrt{s}}{\beta}    \quad  s \in \big(0,\tfrac{\pi^2}{4}\big),\qquad
        g(s) = \cot(\sqrt{s-1})-\frac{\sqrt{s-1}}{\beta} \quad s \in \big(1,1+\tfrac{\pi^2}{4}\big).
     \end{split}
 \end{equation*}
In particular, if $\lambda_{0,0,0}\in (0,1]$ or $\lambda_{0,1,0}\in [\frac{\pi^2}{4},1+\frac{\pi^2}{4})$, there is nothing to prove. Hence, we prove the inequality in the interval  $(1,\frac{\pi^2}{4})$. 

Let us consider the function 
$ h(s)= f(s)-g(s)$, for all $s \in \big(1,\tfrac{\pi^2}{4}\big)$. For any $s>1$, it is obvious that $\sqrt{s-1}<\sqrt{s}$, so that $h(s) < \cot(\sqrt{s})-\cot(\sqrt{s-1})$. If $h(s)<0$ we have finished. Remembering the formula
$    \cot(\alpha-\beta) = \frac{\cot(\alpha)\cot(\beta)+1}{\cot(\beta)-\cot(\alpha)}$, 
if we chose $\alpha= \sqrt{s-1}$ and $\beta= \sqrt{s}$, then 
\begin{equation*}
    h(s)<\cot(\sqrt{s})-\cot(\sqrt{s-1})= \frac{\cot(\sqrt{s-1})\cot{(\sqrt{s})}+1}{\cot(\sqrt{s-1}-\sqrt{s})}.
\end{equation*}
Since $s \in (1,\frac{\pi^2}{4})$, we have  $\cot(\sqrt{s})\cot(\sqrt{s-1})+1>0$ and, since $-\frac{\pi}{2}<\sqrt{s-1}-\sqrt{s}<0$, we have $\cot(\sqrt{s-1}-\sqrt{s})<0$. This proves that $h(s)<0$. \\
The proof for $\lambda_{0,0,0}<\lambda_{1,1,0}$ follows with similar arguments, considering the function $\Bar{h}(s)=f(s)-\Bar{g}(s)$, with $\Bar{g}(s) = g(s-1)$, in the interval $(2,\frac{\pi^2}{4})$.
\end{proof}

We give here the extended version of Corollary~\ref{corollary}, whose proof is in Section~\ref{proof:corollaryN}.
\begin{corollary}[Navier frictionless]
\label{corollaryN1}
Let $\beta=0$ and $\mathcal{P}^u_{m,n}$, $\mathcal{P}^v_{m,n}$, $\mathcal{P}_{m,n}$,  as in~\eqref{periodic} for some $a,b,c,d\in\R$, not contemporary all zero. For any $m,n\in \mathbb{N}$ there exists a non-negative sequence of diverging eigenvalues of~\eqref{stokes} $\lambda^N:=\lambda_{m,n,p}^N$ ($p\in\mathbb{N}$), counted with their multiplicity, given by \eqref{lambdaN}.

\noindent$\bullet$ If $m=n=0$ and $p\in\mathbb N$, the pressure is $q^N_{0,0,p}\equiv Q^N_0\in \mathbb R$ and the eigenfunctions  $(u^N_{0,0,p}(z),v^N_{0,0,p}(z),0)$ are
\begin{equation}\label{eig_0_N}
    u^N_{0,0,p}(z)=\begin{cases}
d\cos(\frac{\pi}{2}p z)\qquad \text{if }p=0\;\text{or even}\\
d\sin(\frac{\pi}{2}p z)\qquad \text{if }p \text{ odd},
\end{cases}\quad v^N_{0,0,p}(z)=\begin{cases}
b\cos(\frac{\pi}{2}p z)\qquad \text{if }p=0\;\text{or even}\\
b\sin(\frac{\pi}{2}p z)\qquad \text{if }p \text{ odd}.
\end{cases}
\end{equation}
for some $b,d\in\mathbb{R}\setminus\{0\}$.

In the remaining cases $q^N_{m,n,p}\equiv 0$ and
\begin{equation}\label{WN}
\mathcal W^N_{m,n,p}(z)=
\begin{cases}
\sin(\frac{\pi}{2}p z)\hspace{5mm} \text{if } p=0\;\text{or even}\\
\cos(\frac{\pi}{2}p z)\hspace{5mm} \text{if }p \text{ odd}.    
\end{cases}
\end{equation}
\noindent$\bullet$ If $n=0$, $m\in \mathbb N$ and $p\in \mathbb N$ the  eigenfunctions are 
$$(\mathcal{U}^N_{m,0,p}(z)\mathcal P^u_{m,0}(x,y),\mathcal{V}^N_{m,0,p}(z)\mathcal P^v_{m,0}(x,y),\mathcal W^N_{m,0,p}(z)\P_{m,0}(x,y)\big),$$ where $\mathcal W^N_{m,0,p}$ is given in~\eqref{WN} and
\begin{equation}\label{final_N2}
\begin{split}
    &  \mathcal{U}^N_{m,0,p}(z)=
\begin{cases}
\frac{\pi p}{2m}\cos(\frac{\pi}{2}p z)\hspace{5mm} \text{if } p=0\;\text{or even}\\
-\frac{\pi p}{2m}\sin(\frac{\pi}{2}p z)\hspace{2mm} \text{if }p \text{ odd}    
\end{cases}
\,
    \mathcal{V}^N_{m,0,p}(z)=
\begin{cases}
\cos(\frac{\pi}{2}p z)\hspace{5mm} \text{if } p=0\;\text{or even}\\\sin(\frac{\pi}{2}p z)\hspace{7.5mm} \text{if }p \text{ odd}.    
\end{cases}
\end{split}
\end{equation}

\noindent$\bullet$ If $m\in \mathbb N$, $n\in\mathbb{N}_+$ and $p\in \mathbb N$ the eigenfunctions are $$(\mathcal{U}^N_{m,n,p}(z)\mathcal P^u_{m,n}(x,y),\mathcal{V}^N_{m,n,p}(z)\mathcal P^v_{m,n}(x,y),\mathcal W^N_{m,n,p}(z)\P_{m,n}(x,y)\big),$$ where $\mathcal W^N_{m,n,p}$ is given in~\eqref{WN} and
\begin{equation}\label{final_N}
\begin{split}
    &  \mathcal{U}^N_{m,n,p}(z)=
\begin{cases}
\cos(\frac{\pi}{2}p z)\hspace{3mm} \text{if } p=0\;\text{or even}\\
\sin(\frac{\pi}{2}p z)\hspace{3mm} \text{if }p \text{ odd}    
\end{cases}
    \mathcal{V}^N_{m,n,p}(z)=
\begin{cases}
\big(-\frac{m}{n}+\frac{\pi p}{2n}\big)\cos(\frac{\pi}{2}p z)\hspace{3mm} \text{if } p=0\;\text{or even}\\
-\big(\frac{m}{n}+\frac{\pi p}{2n}\big)\sin(\frac{\pi}{2}p z)\hspace{5mm} \text{if }p \text{ odd}.    
\end{cases}
\end{split}
\end{equation}
\end{corollary}

We conclude providing the extended version of Corollary~\ref{corollaryD}, whose proof is in Section~\ref{proof:corollaryD}.
\begin{corollary}[Dirichlet]
\label{corollaryD1}
Let $\beta\rightarrow\infty$ and $\mathcal{P}^u_{m,n}$, $\mathcal{P}^v_{m,n}$, $\mathcal{P}_{m,n}$,  as in~\eqref{periodic} for some $a,b,c,d\in\R$, not contemporary all zero. For any $m,n\in \mathbb{N}$ there exist two positive and increasing sequences of diverging eigenvalues of~\eqref{stokes} $\lambda^D:=\lambda^D_{m,n,p}$ ($p\in\mathbb{N}_+$) and $\Lambda^D:=\Lambda^D_{m,n,p}$ ($p\in\mathbb{N}$), counted with their multiplicity, such that:
\begin{itemize}
    \item[i)] If $\nabla q^D\equiv 0$ the eigenvalues are in \eqref{lambdaD}.
    
$\bullet$ If $m=n=0$ and $p\in\mathbb N_+$ the pressure is $q^D_{0,0,p}=Q^D_0\in \R$ and the eigenfunctions are $(u^D_{0,0,p},v^D_{0,0,p},0)$ 
\begin{equation}\label{eigD3d}
u^D_{0,0,p}(z)=\begin{cases}
d\cos(\tfrac{\pi}{2}p z)\qquad \text{if }p \text{ odd}\\
d\sin(\tfrac{\pi}{2}p z)\qquad \text{if }p\, \text{even}
\end{cases} \,
v^D_{0,0,p}(z)=\begin{cases}
b\cos(\tfrac{\pi}{2}p z)\qquad \text{if }p \text{ odd}\\
b\sin(\tfrac{\pi}{2}p z)\qquad \text{if }p\, \text{even}
\end{cases}
\end{equation}
for some $b,d\in\R\setminus\{0\}$.

$\bullet$ If $n=0$, $m\in\mathbb{N}$ and $p\in\mathbb N_+$ the pressure is $q^D_{m,0,p}\equiv 0$ and  the eigenfunctions are 
$$(\mathcal{U}^D_{m,0,p}(z)\mathcal P^u_{m,0}(x,y),\mathcal{V}^D_{m,0,p}(z)\mathcal P^v_{m,0}(x,y),\mathcal W^D_{m,0,p}(z)\P_{m,0}(x,y)\big),$$ where 
\begin{equation}\label{final_D2}
    \mathcal{U}_{m,0,p}^D(z)=\begin{cases}
\frac{\pi p}{2m}\sin(\frac{\pi}{2}p z)\qquad \text{if } p \text{ even}\\
-\frac{\pi p}{2m}\cos(\frac{\pi}{2}p z)\hspace{6mm} \text{if }p \text{ odd}
\end{cases}\qquad \mathcal{V}_{m,0,p}^D(z)=\begin{cases}
\sin(\frac{\pi}{2}p z)\qquad \text{if } p \text{ even}\\
\cos(\frac{\pi}{2}p z)\hspace{7mm} \,\text{if }p \text{ odd}.
\end{cases}
\end{equation}
and 
\begin{equation}\label{WD}
    \mathcal W^D_{m,0,p}(z)=\begin{cases}
    \sin(\frac{\pi p}{2}z)\quad \text{if }p \text{ even}\\
    \cos(\frac{\pi p}{2}z)\quad \text{if }p \text{ odd}.
    \end{cases}
\end{equation}

 $\bullet$ If $m\in\mathbb N$, $n\in\mathbb N_+$ and $p\in\mathbb N_+$ the pressure is $q^D_{0,n,p}\equiv 0$ and the eigenfunctions are
 $$(\mathcal{U}^D_{m,n,p}(z)\mathcal P^u_{m,n}(x,y),\mathcal{V}^D_{m,n,p}(z)\mathcal P^v_{m,n}(x,y),\mathcal W^D_{m,n,p}(z)\P_{m,n}(x,y)\big),$$ where $\mathcal W^D_{m,n,p}$ is given in~\eqref{WD} and
\begin{equation}\label{final_D}
    \mathcal{U}^D_{m,n,p}(z)=\begin{cases}
       \sin(\frac{\pi}{2}pz)\qquad \text{if } p\, \text{even}\vspace{3mm}\\
       \cos(\frac{\pi}{2}pz)\qquad\text{if } p \text{ odd}
    \end{cases}\quad
        \mathcal{V}^D_{m,n,p}(z)=\begin{cases}
         (-\frac{m}{n}+\frac{\pi p}{2n})\sin(\frac{\pi}{2}pz)\qquad \text{if } p\, \text{even}\vspace{3mm}\\
        -(\frac{m}{n}+\frac{\pi p}{2n})\cos(\frac{\pi}{2}pz)\hspace{8.5mm} \text{if } p \text{ odd}.
    \end{cases}
\end{equation}
   \item[ii)] if $\nabla q\not\equiv 0$ the eigenvalues $\La^D\in \big(\mu^2+\frac{\pi^2}{4}(1+p)^2,\mu^2+\frac{\pi^2}{4}(2+p)^2\big)$ with $\mu>0$ satisfy \eqref{eig3dbetainfty}.

$\bullet$ If $n=0$, $m\in\mathbb{N}_+$ and $p\in\mathbb N$ the pressure is 
\begin{equation}\label{press.m0p.appendixD}
   q^D_{m,0,p}(x,z)=\mathcal{P}_{m,0}(x)\begin{cases}
        \La^D\sinh(m z)\qquad \text{if } p=0\,\, \text{or even}
        \\
        \La^D\cosh(m z)\qquad \text{if } p \text{ odd}
    \end{cases}
\end{equation}
 and the eigenfunctions are $(\mathcal U^D_{m,0,p}(z)\mathcal{P}^v_{m,0}(x)\,,\,0\,,\,\mathcal W^D_{m,0,p}(z)\mathcal{P}_{m,0}(x))$,  where
\begin{equation}\label{eigvappendix_m0pD}
\begin{split}
    &\mathcal U^D_{m,0,p}(z)=m\begin{cases}
        \dfrac{\sin(\sqrt{\La^D- m^2})\sinh(m z)-\sinh(m)\sin(\sqrt{\La^D- m^2}z)}{\sin(\sqrt{\La^D- m^2})}\hspace{5mm} \text{if } p=0\,\, \text{or even}\\
\dfrac{\cos(\sqrt{\La^D-m^2})\cosh(m z)- \cosh(m)\cos(\sqrt{\La^D-m^2}z)}{\cos(\sqrt{\La^D-m^2})}\hspace{5mm} \text{if } p \text{ odd}
    \end{cases}\\
    &\mathcal W^D_{m,0,p}(z)=m\begin{cases}
        \dfrac{\cos(\sqrt{\La^D-m^2})\cosh(m z)-\cosh(m) \cos(\sqrt{\La^D-m^2}z)}{\cos(\sqrt{ \La^D-m^2})}\hspace{3mm} \text{if }p=0\,\, \text{or even}\\
        \dfrac{\sin(\sqrt{\La^D-m^2})\sinh(m z)-\sinh(m)\sin(\sqrt{\La^D-m^2}z)}{\sin(\sqrt{ \La^D-m^2})}\hspace{5mm} \text{if } p \text{ odd}.
    \end{cases}
    \end{split}
\end{equation}

 $\bullet$ If $m\in\mathbb{N}$, $n\in\mathbb N_+$ and $p\in\mathbb N$ the pressure is 
\begin{equation}\label{pressD}
   q^D_{m,n,p}(x,y,z)=\mathcal{P}_{m,n}(x,y)\begin{cases}
        \Lambda^D\sinh(\mu z)\qquad \text{if } p=0\,\, \text{or even}\\
        \Lambda^D\cosh(\mu z)\qquad \text{if } p \text{ odd}.
    \end{cases}
\end{equation}
and the eigenfunctions are 
$\big(\mathcal U^D_{m,n,p}(z)\mathcal{P}^u_{m,n}(x,y), \mathcal V^D_{m,n,p}(z)\mathcal{P}^v_{m,n}(x,y),\mathcal W^D_{m,n,p}(z)\mathcal{P}_{m,n}(x,y)\big)$, where
\small{
\begin{equation}\label{eigvD}
\begin{split}
    &\mathcal U^D_{m,n,p}(z)=m\begin{cases}
         \dfrac{\sin(\sqrt{ \La^D-\mu^2})\sinh(\mu z)-\sinh(\mu)\sin(\sqrt{ \La^D-\mu^2}z)}{\sin(\sqrt{ \La^D-\mu^2})}\hspace{5mm} \text{if } p=0\,\, \text{or even}\\
\dfrac{\cos(\sqrt{ \La^D-\mu^2})\cosh(\mu z)-\cosh(\mu)\cos(\sqrt{ \La^D-\mu^2}z)}{\cos(\sqrt{ \La^D-\mu^2})}\hspace{5mm} \text{if } p \text{ odd}
    \end{cases}\\
      &\mathcal V_{m,n,p}^D(z)=\\&\begin{cases}
       n
     \sinh(\mu z)+\dfrac{[m^2\sinh(\mu)+\mu\sqrt{ \La^D-\mu^2}\cosh(\mu)\tan(\sqrt{ \La^D-\mu^2})]\sin(\sqrt{ \La^D-\mu^2}z)}{n\sin(\sqrt{ \La^D-\mu^2})}\hspace{5mm} \text{if } p=0\text{ or even}\\
n\cosh(\mu z)\!+\!\dfrac{[m^2\cosh(\mu)-\mu\sqrt{ \La^D-\mu^2}\sinh(\mu)\cot(\sqrt{ \La^D-\mu^2})]\cos(\sqrt{ \La^D-\mu^2}z)}{n\cos(\sqrt{ \La^D-\mu^2})}\hspace{6mm} \text{if } p \text{ odd}
    \end{cases}\\
    &\mathcal W_{m,n,p}^D(z)=\mu\begin{cases}
        \dfrac{\cos(\sqrt{ \La^D-\mu^2})\cosh(\mu z)-\cosh(\mu)\cos(\sqrt{ \La^D-\mu^2}z)}{\cos(\sqrt{ \La^D-\mu^2})}\hspace{3mm} \text{if } p=0\,\, \text{or even}\\
        \dfrac{\sin(\sqrt{ \La^D-\mu^2})\sinh(\mu z)-\sinh(\mu)\sin(\sqrt{ \La^D-\mu^2}z)}{\sin(\sqrt{ \La^D-\mu^2})}\hspace{5mm} \text{if } p \text{ odd}.
    \end{cases}
    \end{split}
\end{equation}}
\normalsize
\end{itemize}  
 \end{corollary}

\begin{remark}
    We point out that when $\beta\to\infty$ we recover the eigenvalues and the eigenfunctions that were already computed by Rummler in~\cite{Rum1997b}, making consistent our computations.
\end{remark}

\noindent

\section{Helmholtz-Weyl decomposition}\label{nonlin_proj}
 
Let us recall briefly what we mean by the Helmholtz-Weyl~\cite{helmholtz,weyl} orthogonal decomposition of $L^2_\P(\Omega)$. We introduce the space
$$
G:=\{\vv\in L^2_\P(\Omega):\, \exists\, g\in H^1_\P(\Omega),\, \vv=\nabla g\}\, ,
$$
and, recalling the definition of $H_\P$ in~\eqref{HP}, it holds that
\begin{equation}\label{HW}
L^2_\P(\Omega)=H_\P\oplus G\ ,\qquad H_\P\perp G\ .
\end{equation}
Here the orthogonality is intended with respect to the scalar product in $L^2(\Omega)$. The space $G$ is made of weakly irrotational (conservative)
vector fields, namely
$$
G=\{\w\in L^2_\P(\Omega):\ {\rm curl}\, \w={\bf 0}\mbox{ in distributional sense}\}\, .
$$
We use the notation
$$
\forall\,\vv,\w\in L^2_\P(\Omega)\qquad\qquad \vv=\w+G\quad
\begin{array}{c}
\mbox{{\tiny notation}}\\
\Longleftrightarrow\\
\
\end{array}
\quad
\vv-\w\in G\, ,
$$
which means that $\vv$ and $\w$ have the same projection onto $H_\P$. 
To determine the components of a vector field $\Phi\in L^2_\P(\Omega)$
according to~\eqref{HW}, one proceeds as follows. Let $\varphi\in H^1_\P(\Omega)$ be a (scalar) weak solution of the Neumann problem
\neweq{scalNeu}
\Delta\varphi=\nabla\cdot\Phi\quad\mbox{in }\Omega\, ,\qquad \partial_\nu\varphi=\Phi\cdot\nu\quad\mbox{on }\Gamma\, .
\endeq
Since the compatibility condition is satisfied, up to the addition of constants, such $\varphi$ exists and is unique.
Then we observe that $\nabla\cdot(\Phi-\nabla\varphi)=0$ in $\Omega$ and $(\Phi-\nabla\varphi)\cdot\nu=0$
on $\Gamma$. Therefore $(\Phi-\nabla\varphi)\in H_\P$ and we can write
\neweq{HWscomp}
\Phi=(\Phi-\nabla\varphi)+\nabla\varphi=(\Phi-\nabla\varphi)+G\, .
\endeq

Our aim is to apply this procedure to the term $(\uu_{1}\cdot\nabla)\,\uu_{2}$, being $\uu_1$ and $\uu_2$ two eigenfunctions given by Corollary~\ref{corollary:press_cost}. Hence, the first step is to compute $(\uu_{1}\cdot \nabla)\,\uu_{2}$.
\begin{computation}\label{computation1}
    Let $m_1,m_2, p_1,p_2\in\mathbb{N}$, $n_1,n_2\in\mathbb{N}_+$, $\uu_{m_1,n_1,p_1}(x,y,z)$ and $\uu_{m_2,n_2,p_2}(x,y,z)$ two eigenfunctions in Corollary~\ref{corollary:press_cost}. Then, it holds 
    \begin{equation}
    \label{compu}
    \begin{split}
        (\uu_{m_1,n_1,p_1}\cdot\nabla)\,&\uu_{m_2,n_2,p_2}=\dfrac{\mathcal{Z}_{m_1,n_1,p_1}(z)\mathcal{Z}_{m_2,n_2,p_2}(z)}{4n_1}\cdot\\&\cdot\bigg\{(m_2n_1+m_1n_2)\left(\begin{matrix}
            -2(ab+cd)\cos[(m_1+m_2)x]\cos[(n_1-n_2)y]\\+2(bc+ad)\sin[(m_1-m_2)x]\sin[(n_1+n_2)y]\\-(a^2+b^2-c^2-d^2)\sin[(m_1-m_2)x]\cos[(n_1+n_2)y]\\-(a^2-b^2-c^2+d^2)\sin[(m_1+m_2)x]\cos[(n_1-n_2)y]\\
            \\
            \dfrac{m_2}{n_2}\big\{2(ab+cd)\sin[(m_1+m_2)x]\sin[(n_1-n_2)y]\\-2(bc+ad)\cos[(m_1-m_2)x]\cos[(n_1+n_2)y]\\-(a^2+b^2-c^2-d^2)\cos[(m_1-m_2)x]\sin[(n_1+n_2)y]\\-(a^2-b^2-c^2+d^2)\cos[(m_1+m_2)x]\sin[(n_1-n_2)y]\big\}\\
            \\
            0
        \end{matrix}\right)\\
        &+(m_2n_1-m_1n_2)\left(\begin{matrix}
            2(ab-cd)\cos[(m_1+m_2)x]\cos[(n_1+n_2)y]\\-2(ac-bd)\cos[(m_1-m_2)x]\sin[(n_1-n_2)y]\\
            -2(ac+bd)\cos[(m_1+m_2)x]\sin[(n_1+n_2)y]\\+2(bc-ad)\sin[(m_1+m_2)x]\sin[(n_1+n_2)y]\\+(a^2+b^2+c^2+d^2)\sin[(m_1-m_2)x]\cos[(n_1-n_2)y]\\+(a^2-b^2+c^2-d^2)\sin[(m_1+m_2)x]\cos[(n_1+n_2)y]\\
            \\
            \dfrac{m_2}{n_2}\big\{2(ab-cd)\sin[(m_1+m_2)x]\sin[(n_1+n_2)y]\\+2(ac-bd)\sin[(m_1-m_2)x]\cos[(n_1-n_2)y]\\
            +2(ac+bd)\sin[(m_1+m_2)x]\cos[(n_1+n_2)y]\\+2(bc-ad)\cos[(m_1+m_2)x]\cos[(n_1+n_2)y]\\-(a^2+b^2+c^2+d^2)\cos[(m_1-m_2)x]\sin[(n_1-n_2)y]\\-(a^2-b^2+c^2-d^2)\cos[(m_1+m_2)x]\sin[(n_1+n_2)y]\big\}
            \\
            \\
            0
        \end{matrix}\right)\bigg\}\,.
     \end{split}
    \end{equation}
    \end{computation}
The second step is to decompose~\eqref{compu} onto the space $H_\P$, applying the Helmholtz-Weyl decomposition. This operation is possible, but it is very heavy. 
For sake of clarity we do the decomposition only for a term in~\eqref{compu}, e.g. that multiplied by $(ab-cd)$ assuming $p_1,p_2$ zero or even. The others can be computed in the same way.
 \begin{computation}\label{computation2}
     Under the assumptions of Computation~\ref{computation1} and $p_1,p_2$ zero or even, let $\alpha:=(m_1+m_2)$, $\beta:=(n_1+n_2)$, $\gamma:=\sqrt{\lambda_1-m_1^2-n_1^2}$, $\delta:=\sqrt{\lambda_2-m_2^2-n_2^2}$ and
     $$
     \Phi(x,y,z)=(ab-cd)\dfrac{m_2n_1-m_1n_2}{2n_1}\cos(\gamma z)\cos(\delta z)\left(\begin{matrix}
            \cos(\alpha x)\cos(\beta y)\\
            \dfrac{m_2}{n_2}\sin(\alpha x)\sin(\beta y)\\
            0
        \end{matrix}\right).
     $$
 Then
   	\begin{equation}\label{NN0}
 \begin{split}
    \Phi(x,y,z)=&(ab-cd)\dfrac{m_2n_1-m_1n_2}{2n_1n_2}\cdot\\&\cdot\left(\begin{matrix}
            \big[n_2\cos(\gamma z)\cos(\delta z)-(m_2n_1-m_1n_2)\alpha \upsilon(z)\big]\cos(\alpha x)\cos(\beta y)\\
            \big[m_2\cos(\gamma z)\cos(\delta z)+(m_2n_1-m_1n_2)\beta \upsilon(z)\big]\sin(\alpha x)\sin(\beta y)\\
            -(m_2n_1-m_1n_2)\upsilon'(z)\sin(\alpha x)\cos(\beta y)
        \end{matrix}\right)+G,
 \end{split}
\end{equation}
where 
\begin{equation}\label{ZZ0}
    \begin{split}
         \upsilon(z)=&\dfrac{-1}{[(\gamma-\delta)^2+\alpha^2+\beta^2][(\gamma+\delta)^2+\alpha^2+\beta^2]}\bigg[2\gamma\delta \sin(\gamma z)\sin(\delta z)+(\alpha^2+\beta^2+\gamma^2+\delta^2)\cos(\gamma z)\cos(\delta z)\\&\frac{\gamma \sin (\gamma) \cos (\delta) \left(\alpha^2+\beta^2+\gamma^2-\delta^2\right)+\delta \cos (\gamma) \sin (\delta) \left(\alpha^2+\beta^2-\gamma^2+\delta^2\right)}{\sqrt{\alpha^2+\beta^2} \sinh \left(\sqrt{\alpha^2+\beta^2}\right)}\cosh(\sqrt{\alpha^2+\beta^2}z)\bigg].
    \end{split}
\end{equation}
\end{computation}
\begin{proof}
	We observe that $\Phi\cdot\nu=0$ on $\Gamma$, since the third component of $\Phi$ is zero. We then compute
	\begin{equation*}
	\nabla\cdot \Phi=(ab-cd)\dfrac{(m_2n_1-m_1n_2)^2}{2n_1n_2}\cos(\gamma z)\cos(\delta z)\sin(\alpha x)\cos(\beta y),
	\end{equation*}
	and, we have to find $\varphi$ satisfying~\eqref{scalNeu}. To do this  we search $$\varphi(x,y,z)=(ab-cd)\frac{(m_2n_1-m_1n_2)^2}{2n_1n_2}\upsilon(z)\sin(\alpha x)\cos(\beta y),
 $$
 solving the problem
 $$
 \begin{cases}
 \upsilon''(z)-[\alpha^2+\beta^2]\upsilon(z)=\cos(\gamma z)\cos(\delta z)\\
 \upsilon'(\pm 1)=0,
 \end{cases}
 $$
 since $p_1,p_2$ are zero or even. Hence we obtain~\eqref{ZZ0}.
Computing $\Phi-\nabla\varphi$, by~\eqref{HWscomp} we obtain~\eqref{NN0}.    
\end{proof}

For our aims the key point is that this operation does not affect the coefficients depending on $a,b,c,d$ in front of each line of~\eqref{compu} and produces some complications only in the $z-$part of the field; for instance in Computation~\ref{computation2} the term $(ab-cd)$ is preserved as well $\cos(\alpha x)\cos(\beta y)$ in the first component and $\sin(\alpha x)\sin (\beta y)$ in the second, see~\eqref{NN0}. Since we are not interested in finding the exact analytical expression of the decomposed vector field we state the result  we need in a slightly general form.
\begin{computation}\label{computation3}
    Let $m_1,m_2, p_1,p_2\in\mathbb{N}$, $n_1,n_2\in\mathbb{N}_+$, $\uu_{m_1,n_1,p_1}(x,y,z)$ and $\uu_{m_2,n_2,p_2}(x,y,z)$ two eigenfunctions in Corollary~\ref{corollary:press_cost}. Then there exist $Z_i^u(z)$, $Z_i^v(z)$ and $Z_i^w(z)$ for $i=1,\dots,10$ such that
    \begin{equation*}
    \begin{split}
        &(\uu_{m_1,n_1,p_1}\cdot\nabla)\,\uu_{m_2,n_2,p_2}=G+\\
        \\&\left(\begin{matrix}
            (ab+cd)Z_1^u(z)\cos[(m_1+m_2)x]\cos[(n_1-n_2)y]+(bc+ad)Z_2^u(z)\sin[(m_1-m_2)x]\sin[(n_1+n_2)y]\\+(a^2+b^2-c^2-d^2)Z_3^u(z)\sin[(m_1-m_2)x]\cos[(n_1+n_2)y]\\+(a^2-b^2-c^2+d^2)Z_4^u(z)\sin[(m_1+m_2)x]\cos[(n_1-n_2)y]\\+(ab-cd)Z_5^u(z)\cos[(m_1+m_2)x]\cos[(n_1+n_2)y]+(ac-bd)Z_6^u(z)\cos[(m_1-m_2)x]\sin[(n_1-n_2)y]\\
            +(ac+bd)Z_7^u(z)\cos[(m_1+m_2)x]\sin[(n_1+n_2)y]+(bc-ad)Z_8^u(z)\sin[(m_1+m_2)x]\sin[(n_1+n_2)y]\\+(a^2+b^2+c^2+d^2)Z_9^u(z)\sin[(m_1-m_2)x]\cos[(n_1-n_2)y]\\+(a^2-b^2+c^2-d^2)Z_{10}^u(z)\sin[(m_1+m_2)x]\cos[(n_1+n_2)y]\\
            \\
            (ab+cd)Z_1^v(z)\sin[(m_1+m_2)x]\sin[(n_1-n_2)y]+(bc+ad)Z_2^v(z)\cos[(m_1-m_2)x]\cos[(n_1+n_2)y]\\+(a^2+b^2-c^2-d^2)Z_3^v(z)\cos[(m_1-m_2)x]\sin[(n_1+n_2)y]\\+(a^2-b^2-c^2+d^2)Z_4^v(z)\cos[(m_1+m_2)x]\sin[(n_1-n_2)y]\\+(ab-cd)Z_5^v(z)\sin[(m_1+m_2)x]\sin[(n_1+n_2)y]+(ac-bd)Z_6^v(z)\sin[(m_1-m_2)x]\cos[(n_1-n_2)y]\\
            +(ac+bd)Z_7^v(z)\sin[(m_1+m_2)x]\cos[(n_1+n_2)y]+(bc-ad)Z_8^v(z)\cos[(m_1+m_2)x]\cos[(n_1+n_2)y]\\+(a^2+b^2+c^2+d^2)Z_9^v(z)\cos[(m_1-m_2)x]\sin[(n_1-n_2)y]\\+(a^2-b^2+c^2-d^2)Z_{10}^v(z)\cos[(m_1+m_2)x]\sin[(n_1+n_2)y]\\
            \\
            +(ab+cd)Z_1^w(z)\sin[(m_1+m_2)x]\cos[(n_1-n_2)y]+(bc+ad)Z_2^w(z)\cos[(m_1-m_2)x]\sin[(n_1+n_2)y]\\+(a^2+b^2-c^2-d^2)Z_3^w(z)\cos[(m_1-m_2)x]\cos[(n_1+n_2)y]\\+(a^2-b^2-c^2+d^2)Z_4^w(z)\cos[(m_1+m_2)x]\cos[(n_1-n_2)y]\\+(ab-cd)Z_5^w(z)\sin[(m_1+m_2)x]\cos[(n_1+n_2)y](ac-bd)Z_6^w(z)\sin[(m_1-m_2)x]\sin[(n_1-n_2)y]\\
            +(ac+bd)Z_7^w(z)\sin[(m_1+m_2)x]\sin[(n_1+n_2)y]+(bc-ad)Z_8^w(z)\cos[(m_1+m_2)x]\sin[(n_1+n_2)y]\\+(a^2+b^2+c^2+d^2)Z_9^w(z)\cos[(m_1-m_2)x]\cos[(n_1-n_2)y]\\+(a^2-b^2+c^2-d^2)Z_{10}^w(z)\cos[(m_1+m_2)x]\cos[(n_1+n_2)y]
        \end{matrix}\right).
     \end{split}
    \end{equation*}
    \end{computation}
\begin{proof}
    The proof follows computing the functions $Z_i^u(z)$, $Z_i^v(z)$ and $Z_i^w(z)$ for $i=1,\dots,10$ as we did for $Z_5^u(z)$, $Z_5^v(z)$ and $Z_5^w(z)$ in the Computation~\ref{computation2}. 
\end{proof}
%
\def\ocirc#1{\ifmmode\setbox0=\hbox{$#1$}\dimen0=\ht0 \advance\dimen0
  by1pt\rlap{\hbox to\wd0{\hss\raise\dimen0
  \hbox{\hskip.2em$\scriptscriptstyle\circ$}\hss}}#1\else {\accent"17 #1}\fi}
  \def\polhk#1{\setbox0=\hbox{#1}{\ooalign{\hidewidth
  \lower1.5ex\hbox{`}\hidewidth\crcr\unhbox0}}} \def\cprime{$'$}

\end{document}